\documentclass[10pt,a4paper]{article}
\usepackage[utf8]{inputenc}
\usepackage{amsmath}
\usepackage{amsfonts}
\usepackage{amssymb}
\usepackage{amsthm}
\usepackage{algpseudocode}
\usepackage{mathtools}
\usepackage{nicefrac}
\usepackage{a4wide}
\usepackage[small,it]{caption}
\usepackage{subfig}
\usepackage{url}
\usepackage{color}
\usepackage{soul}
\usepackage{mathtools}
\usepackage{pifont}
\usepackage{algorithm,algorithmicx,algpseudocode}

\newcommand{\vect}[1]{\boldsymbol{#1}}
\newcommand{\avg}[1]{{\{\!\!\{ #1 \}\!\!\}}}
\newcommand{\jmp}[1]{{[\![ #1 ]\!]}}
\DeclarePairedDelimiter{\abs}{\lvert}{\rvert}
\DeclarePairedDelimiter{\norm}{\lVert}{\rVert}

\newcommand{\potdir}{\hat{\vect{\frak{d}}}_K}
\newcommand{\actdir}{\vect{\frak{d}}_K}
\newcommand{\e}{\mathrm{e}}
\newcommand{\real}{\mathbb{R}}
\newcommand{\complex}{\mathbb{C}}
\newcommand{\mesh}[1][h]{\mathcal{T}_{#1}}
\newcommand{\face}[1][]{\mathcal{F}_h^{#1}}
\newcommand{\iin}{\text{in }}
\newcommand{\oon}{\text{on }}
\newcommand{\conj}[1]{\overline{#1}}

\newcommand{\ds}{\, \mathrm{d}s}
\newcommand{\Hessian}{\boldsymbol{\mathcal{H}}}
\newcommand{\cmark}{\ding{51}}
\newcommand{\xmark}{\ding{55}}
\newcommand{\degrees}{^\circ}

\DeclareMathOperator{\Real}{Re}
\DeclareMathOperator{\Imag}{Im}

\theoremstyle{plain}
\newtheorem{theorem}{Theorem}[section]

\theoremstyle{remark}
\newtheorem{remark}[theorem]{Remark}

\newcommand{\configfigure}{\captionsetup[subfloat]{farskip=0pt,captionskip=5pt}\centering}

\title{Adaptive Refinement for $hp$--Version Trefftz Discontinuous
  Galerkin Methods for the Homogeneous Helmholtz Problem}
\author{
Scott Congreve\thanks{Faculty of Mathematics, University of Vienna,
Oskar-Morgenstern-Platz 1,
  1090 Vienna, Austria (\tt{scott.congreve@univie.ac.at})},
Paul Houston\thanks{School of Mathematical Sciences, University of Nottingham,
University Park, Nottingham, NG7 2RD, UK (\tt{Paul.Houston@nottingham.ac.uk})},  
Ilaria Perugia\thanks{Faculty of Mathematics, University of Vienna,
Oskar-Morgenstern-Platz 1,
  1090 Vienna, Austria (\tt{ilaria.perugia@univie.ac.at})}}
\date{}

\begin{document}

\maketitle

\begin{abstract}
In this article we develop an $hp$--adaptive refinement procedure for Trefftz discontinuous
Galerkin methods applied to the homogeneous Helmholtz problem. Our approach combines 
not only mesh subdivision ($h$--refinement) and local basis enrichment ($p$--refinement),
but also incorporates local directional adaptivity, whereby the elementwise plane wave
basis is aligned with the dominant scattering direction. 
Numerical experiments based on employing an empirical {\em a posteriori} error indicator
clearly highlight the efficiency of the proposed approach for various examples. 
\end{abstract}
 
\section{Introduction}
Trefftz discontinuous Galerkin (TDG) methods are finite element schemes which employ discontinuous 
test and trial functions whose restriction to each mesh element belongs to the kernel of the differential 
operator to be discretized. For time-harmonic wave problems, Trefftz discretization spaces 
are made of oscillating functions with the same frequency as that of the underlying
analytical solution.
This results in improved approximation properties, as compared to standard piecewise polynomial 
spaces. Moreover, based on Trefftz spaces, one can construct discontinuous Galerkin methods 
which feature unconditional unique solvability, as well as coercivity of the discrete bilinear forms
in suitable (mesh-dependent) norms. 
We focus here on the case of the Helmholtz problem and refer, e.g., to the survey~\cite{Hiptmair2016}
for a review of the construction, properties, and relevant literature of Trefftz methods
for its numerical approximation.

The purpose of this article is to develop an efficient $hp$--adaptive refinement
algorithm for TDG methods applied to the homogeneous Helmholtz problem; we will specifically 
consider the ultra-weak variational formulation with plane wave basis
functions~\cite{Cessenat1998}.
Within the adaptive procedure, elements will be marked for refinement based on employing an
empirical {\em a posteriori} error indicator, stimulated by the upper
bounds derived in~\cite{Kapita2015} for the $h$--version of the TDG method.
For the $h$--version of the plane wave discontinuous Galerkin method, incorporating Lagrange multipliers, a similar error
indicator has been presented in \cite{Amara2009}.
Once an element has been marked for refinement, a decision must then
be made regarding the type of refinement to be undertaken, i.e., whether the
element should be subdivided ($h$--refinement), or whether the local basis should
be enriched ($p$--refinement). The choice of whether to $h$-- or $p$--refine an element
is typically based on the observation that when the underlying solution is
smooth, then $p$--refinement will be more efficient in terms of reducing the
error, for a given increase in the number of degrees of freedom, than if the
element is subdivided. On the other hand, if the solution is not smooth, then
$h$--refinement should be employed. In general, {\em a posteriori} error estimators
only provide an estimate of the local elementwise error, but do not indicate
which type of refinement should be employed. Within the existing literature
a number of algorithms have been devised for determining the type of
refinement ($h$-- or $p$--) to be undertaken. For a comprehensive review of this
subject, we refer to~\cite{Mitchell2011,Mitchell2014}, and the references
cited therein. In the present context, given the oscillatory nature
of solutions to high-frequency scattering problems, the exploitation of
$hp$--strategies based on local regularity estimation techniques is
not generally applicable. Thereby, we consider an alternative approach based on 
estimating the predicted decay rate of the {\em a posteriori} estimator,
given the refinement history of each element; see, for example,
\cite{Melenk2001}.
For {\em a posteriori} error estimation of conforming
finite element approximations of the Helmholtz problem, we refer,
e.g., to \cite{BabI,BabII} and~\cite{SauterDoerfler}; analogous bounds have been established
for polynomial-based discontinuous
Galerkin finite element methods in \cite{SauterZech, Zech}.

In addition to standard $h$-- and $hp$--adaptivity, we also consider the
issue of directional refinement of the underlying plane wave basis employed 
within our TDG scheme. In particular, we rotate the underlying elementwise
plane wave basis in order that the first basis function is aligned with
the local dominant propagation direction; strategies for determining the
local dominant propagation direction have been proposed 
in~\cite{Amara2014,Betcke2011,Betcke2012,Gittelson2008}, for example.
Stimulated by the work undertaken on
anisotropic mesh adaptation in~\cite{Formaggia2001,Formaggia2003}, cf., 
also,~\cite{ghh-paper,hall-thesis}, we propose an alternative approach
based on studying the properties of the Hessian of the computed TDG solution.
More precisely, the principal eigenvector of the Hessian of the 
solution indicates the dominant direction of wave propagation. However,
since eigenvectors are only unique up to scalar multiples, the
precise wave direction must be fixed, based on exploiting an impedance condition.
In this way, we can locally orientate the elementwise plane wave
basis to reduce the error in the underlying computed TDG solution
in a simple and computationally cheap manner. When combined with 
$hp$--refinement, the resulting adaptive procedure is capable of
generating highly optimized $hp$--refined Trefftz spaces. Indeed,
the efficiency of the proposed strategy is illustrated for a number
of test problems, where we compare the performance between
an $h$-- and $hp$--refinement algorithm, both with and without directional adaptivity.

The outline of this article is as follows: in Section~\ref{sec:pde_tdg}
we introduce the model problem to be studied within this article, together
with its TDG discretization. Then in Section~\ref{section:adaptive_refinement}
we develop an $hp$--refinement algorithm, based on employing both local
mesh subdivision and local basis enrichment, together with directional
adaptivity for the underlying Trefftz space. The performance of
this procedure is studied in Section~\ref{sec:numerical_examples}
through a series of two-- and three--dimensional examples. Finally,
in Section~\ref{sec:conclusions} we summarize the work undertaken
within this article and highlight potential future directions of research.

\section{Model problem and TDG discretization} \label{sec:pde_tdg}

In this section we state the model problem to be studied in this
article, together with its TDG discretization; for further details, we refer
to
\cite{Hiptmair2016},
for example.

\subsection{Model problem}

We study the homogeneous Helmholtz 
equation; to this end, we let 
$\Omega\subset\real^d$, $d=2,3$, be an open bounded, Lipschitz domain with boundary 
$\partial \Omega$. Thereby, we seek $u:\Omega\mapsto {\mathbb C}$ such that
\begin{equation}
\label{eqn:helmholtz}
  \begin{aligned}
    -\Delta u-k^2 u &=0 && \iin \Omega\;,\\
    \frac{\partial u}{\partial \vect{n}} + ik\vartheta u &= g_R && \oon \Gamma_R\;, \\
    u &=g_D && \oon \Gamma_D\;,
  \end{aligned}  
\end{equation}
where $\vect{n}$ denotes the unit outward normal vector on the boundary $\partial\Omega$,
and $\Gamma_R$ and $\Gamma_D$ are non-overlapping open subsets of $\partial\Omega$, such
that $\partial\Omega=\overline{\Gamma}_R\cup\overline{\Gamma}_D$. Furthermore,
$i$ is the imaginary unit, $\vartheta=\pm 1$, $g_R\in L^2(\Gamma_R)$, and we assume, for 
the moment, that the (real-valued) wavenumber $k$ is constant in $\Omega$.

\subsection{Meshes and spaces}
We partition $\Omega$ into computational meshes $\{\mesh\}_{h>0}$ consisting of 
non-overlapping (curvilinear) polygons/polyhedra $K$, which potentially include hanging nodes, 
such that $\overline{\Omega}=\bigcup_{K\in\mesh}{\overline{K}}$. Moreover,
we assume that the family of subdivisions $\{\mesh\}_{h>0}$ is shape regular 
\cite[pp. 61, 114, and 118]{Braess}.
For each element $K\in\mesh$, we write $h_K$ to denote its diameter and
$\vect{n}_K$ signifies the unit outward normal vector to $K$ on~$\partial K$; we set
$h:=\max_{K\in\mesh}h_K$.
Furthermore, we introduce the mesh skeleton $\face$, defined by $\face=\cup_{K\in\mesh}\partial K$;
we write $\face[I]$ and $\face[B]$ to denote the interior and boundary 
skeletons, respectively, defined by
$\face[I]=\face\setminus\partial\Omega$ and $\face[B]=\partial\Omega$.
Implicitly, we assume that the finite element mesh $\mesh$ respects the decomposition 
of the boundary, in the sense that, given an element face $f\subset \partial K$, $K\in\mesh$, 
which lies on the boundary $\partial\Omega$, i.e., $f\subset\partial\Omega$, then 
$f$ is entirely contained within either $\Gamma_R$ or $\Gamma_D$. 

Let $K$ and $K'$ be two adjacent elements of~ $\mesh$, and
$\vect{x}$ an arbitrary point on the interior face $f\subset \face[I]$ given by
$f=(\partial K\cap\partial K')^\circ$.  Furthermore, let $v$
and~$\vect{w}$ be scalar- and vector-valued functions, respectively,
that are sufficiently smooth inside each
element~$K,K'$. Then, the averages of $v$ and
$\vect{w}$ at $\vect{x}\in f$ are given by
\[
\avg{v}=\frac{1}{2}(v|_{K}+v|_{K'}), \qquad \avg{\vect{w}}
=\frac{1}{2}(\vect{w}|_{K}+\vect{w}|_{K'}),
\] 
respectively. Similarly, the jumps of $v$ and $\vect{w}$ at $\vect{x}\in f$ are given by
\[
\jmp{v} =v|_{K}\,\vect{n}_{K}+v|_{K'}\,\vect{n}_{K'},\qquad
\jmp{\vect{w}}=\vect{w}|_{K}\cdot\vect{n}_{K}+\vect{w}|_{K'}\cdot\vect{n}_{K'},
\]
respectively.

Given $K\in\mesh$ the local Trefftz space is defined by
\[
    T(K) \coloneqq \{ v\in H^1(K) : -\Delta v - k^2 v = 0 \};
\]
with this notation, we write
\[
    T(\mathcal{T}_h) \coloneqq \{ v\in L^2(\Omega) : v\vert_K \in T(K), K\in\mesh\}.
\]
Thereby, given a local space $V_{p_K}(K) \subset T(K)$, of finite dimension $p_K\geq 1$,
the corresponding TDG finite element space is defined by
\[
    V_{\vect{p}}(\mathcal{T}_h) \coloneqq \{ v \in T(\mesh) : v\vert_K \in V_{p_K}(K), K\in\mesh \},
\]
where $\vect{p} = \{p_K:K\in\mesh\}$.

\subsection{TDG discretization}
Equipped with the TDG finite element space $V_{\vect{p}}(\mesh)$ defined on the mesh partition 
$\mesh$ of $\Omega$, the TDG approximation of~\eqref{eqn:helmholtz} is given by: find $u_{hp}\in V_{\vect{p}}(\mesh)$ such that
\begin{equation}\label{eqn:bilinear_form}
\mathcal{A}_{h}(u_{hp},v_{hp})=\ell_{h}(v_{hp})
\end{equation}
for all $v_{hp}\in V_{\vect{p}}(\mesh)$, where
\begin{align*}
    \mathcal{A}_{h}(u,v)=
        & \int_{\face[I]} \left( \avg{u}\jmp{\nabla_h \conj{v}}
        -\beta(ik)^{-1} \jmp{\nabla_h u}\jmp{\nabla_h \conj{v}} 
        -\avg{\nabla_h u}\cdot\jmp{\conj{v}}
        +\alpha ik \jmp{u}\cdot\jmp{\conj{v}} \right)\ds\\
        &+\int_{\Gamma_R}( (1-\delta)(u\nabla_h \conj{v}\cdot\vect{n}+ik\vartheta  u\conj{v})
        -\delta((ik\vartheta)^{-1}(\nabla_h u\cdot\vect{n})(\nabla_h \conj{v}\cdot\vect{n})
        +\nabla_h u\cdot\vect{n}\,\conj{v}))\ds \\
        &+\int_{\Gamma_D}(-\nabla_h u\cdot\vect{n}\,\conj{v}+ \alpha ik u\conj{v})\ds, \\
    \ell_{h}(v)=&
        \int_{\Gamma_R}g_R((1-\delta)\conj{v}
        -\delta(ik\vartheta)^{-1}\nabla_h \conj{v}\cdot\vect{n})\ds 
        +\int_{\Gamma_D}g_D(\alpha ik \conj{v} -\nabla_h \conj{v}\cdot\vect{n})\ds,
\end{align*}
and $\nabla_h$ denotes the broken gradient operator, defined elementwise. Here,
$\alpha>0$, $\beta>0$ and $0<\delta\leq \nicefrac12$ are given penalty parameters.
We note that the selection of these penalty parameters has been studied in a number
of different contexts within the literature; in particular, here we mention
the ultra-weak variational formulation (UWVF), cf. \cite{Cessenat1998}, the DG-type scheme studied 
in~\cite{Gittelson2009}, and \cite{Hiptmair2014} which considered their selection
on locally refined meshes; cf.~\cite[Table 1]{Hiptmair2016}. For the purposes of 
this article we consider the UWVF, corresponding to the choice 
$\alpha=\beta=\delta=\nicefrac{1}{2}$.

\subsection{Plane wave basis functions}

Finally, in this section we outline the choice of the underlying discrete space 
$V_{p_K}(K)$, $K\in\mesh$. To this end, we select $V_{p_K}(K)$ to be
a local space consisting of plane waves in $p_K$ different directions, 
all with the same wavenumber $k$. We note that, under suitable assumptions on $K$ 
and the choice of plane wave directions, $V_{p_K}(K)$ approximates smooth Trefftz functions 
with the same order of convergence as polynomials of degree $q_K$, where
\begin{equation}
    p_K =
    \begin{cases}
        2q_K + 1, & d=2, \\
        (q_K+1)^2, & d= 3;
    \end{cases}
    \label{eqn:effective_poly_deg}
\end{equation}
see~\cite{Moiola2011}. Thereby, $q_{K}$ is referred to as the \emph{effective polynomial degree} 
of the discrete Trefftz space; we set $\vect{q} = \{q_K:K\in\mesh\}$.
More precisely, we write
\begin{equation}\label{eqn:pw_basis}
V_{p_K}(K) \coloneqq \left\{ v \in T(K) : v(\vect{x}) = \sum_{\ell=0}^{p_K-1} \alpha_\ell 
\e^{ik\vect{d}_{K,\ell}\cdot(\vect{x}-\vect{x}_K)}, \alpha_\ell\in\complex\right\},
\end{equation}
where $\vect{x}_K$ is the center of mass of element $K$ and $\vect{d}_{K,\ell}$, $\ell=0,\dots,p_K-1$, 
are $p_K$ evenly distributed unit direction vectors (with respect to the unit ball). 
For $d=2$ we can simply define
\begin{equation}
\label{eqn:plane_wave:2d}
\vect{d}_{K,\ell} = (\cos(\nicefrac{2\pi\ell}{p_K}), \sin(\nicefrac{2\pi\ell}{p_K}))^\top, \qquad
\ell=0,\dots,p_K-1;
\end{equation}
for $d=3$ we employ the directions determined by the extremal (maximum determinant) points 
on $S^2$, cf. \cite{Sloan2004,Womersley2007Online}.

\section{Adaptive mesh refinement}\label{section:adaptive_refinement}

In this section we develop an automatic adaptive refinement algorithm which is
capable of not only marking elements for refinement, but also determining the type
of refinement to be undertaken. In particular, here we consider both $h$-- and $p$--refinement,
whereby the local element is subdivided, or the number of elementwise plane wave directions is 
enriched, respectively, as well as 
directional refinement which seeks to rotate the local plane wave basis
in order to align it with the principal scattering direction. 

\subsection{\emph{A posteriori} error indicator}\label{sec:indicator}

In the absence of rigorous {\em a posteriori} error bounds for the numerical approximation of
\eqref{eqn:helmholtz} by the TDG scheme \eqref{eqn:bilinear_form}, which are sharp
with respect to both the local mesh size $h_K$ and the number of local plane waves $p_K$ 
employed on each element $K\in\mesh$, we employ an empirical error estimator stimulated
by the work undertaken in \cite{Kapita2015} in the $h$--version setting. To this end,
we first introduce the dual problem: find $z\in H^1(\Omega)$, such that
\begin{equation*}
  \begin{aligned}
    -\Delta z-k^2 z &=u-u_{hp} && \iin \Omega\;,\\
    \frac{\partial z}{\partial \vect{n}} + ik\vartheta z &= 0 && \oon \Gamma_R\;, \\
    z &=0 && \oon \Gamma_D\;.
  \end{aligned}  
\end{equation*}
Noting that $z\in H^{\nicefrac{2}{3}+s}(\Omega)$, $0<s\leq\nicefrac{1}{2}$,
cf.~\cite{Hiptmair2014}, we recall the following (second) {\em a posteriori}
error bound from \cite{Kapita2015}.

\begin{theorem} \label{thm:apost}
Assume that the mesh $\mesh$ is shape-regular, locally quasi-uniform, in the 
sense that, for two elements $K$ and $K'$ which share a face $f\subset\face[I]$, there 
is a constant $\tau$, independent of $h$, such that
$$
\tau^{-1} \leq \nicefrac{h_K}{h_{K'}} \leq \tau
$$
for all choices of $K$ and $K'$, and that $\mesh$ is quasi-uniform in the vicinity
of $\Gamma_R$, i.e., for all $K\in\mesh$ which lie on the boundary $\Gamma_R$, i.e.,
so that $\partial K\cap \Gamma_R \neq \emptyset$, there exists $\tau_R$ such that
$$
\nicefrac{h}{h_K} \leq \tau_R.
$$
Then, for $g_D\equiv 0$ and fixed $p_K$, $K\in\mesh$, the following {\em a posteriori} bound holds:
	\begin{equation*}
    \norm{u-u_{hp}}_{L^2(\Omega)} \leq {\mathfrak E}(u_{hp},h) \equiv C \left(\sum_{K\in\mesh} \eta_K^2 \right)^{\nicefrac{1}{2}},
\end{equation*}
where
\begin{equation}
\begin{aligned}
\eta_K^2 &=\norm*{\alpha^{\nicefrac{1}{2}} h_K^{s} \jmp{u_{hp}}}_{L^2(\partial K \setminus \partial\Omega)}^2
        +k^{-2}\norm*{\beta^{\nicefrac12} h_K^{s} \jmp{\nabla u_{hp}}}_{L^2(\partial K \setminus \partial\Omega )}^2 \\
        &\quad+ k^{-2}\norm*{\delta^{\nicefrac{1}{2}} h_K^{s} \left( g_R - \nabla u_{hp}\cdot\vect{n}_K + ik u_{hp}\right)}_{L^2(\partial K \cap \Gamma_R)}^2
        +\norm*{\alpha^{\nicefrac{1}{2}} h_K^{s} u_{hp}}_{L^2(\partial K \cap \Gamma_D)}^2,
\end{aligned}
\label{eqn:error_indicator_kmw}
\end{equation}
where $C$ is a positive constant, which is independent of $h$.
\end{theorem}

We stress that the {\em a posteriori} error bound stated in Theorem~\ref{thm:apost}
depends on the regularity index $s$; thereby, {\em a priori} knowledge of $s$ is
required in order to yield a fully computable bound. Moreover, the dependence of
$C$ on $\vect{p}$, or equivalently $\vect{q}$, is unclear; indeed, to the best of our 
knowledge, an $hp$--version 
generalization of Theorem~\ref{thm:apost} is not currently available within the literature. 
Thereby, we propose the following {\em empirical} error estimator, where for simplicity
of notation we also denote it by ${\mathfrak E}$, for 
the $hp$--version TDG method:
\begin{equation}
\label{eqn:error_bound}
    {\mathfrak E}(u_h,h,\vect{p}) =  \left( \sum_{K\in\mesh} \eta_K^2\right)^2,
\end{equation}
where
\begin{equation}
\begin{aligned}
\eta_K^2 &=\norm*{\alpha^{\nicefrac{1}{2}} h_K^{\nicefrac12}q_K^{-\nicefrac12} \jmp{u_{hp}}}_{L^2(\partial K \setminus \partial\Omega )}^2
        +\norm*{\beta^{\nicefrac12}
          h_K^{\nicefrac32}q_K^{-\nicefrac32} \jmp{\nabla
            u_{hp}}}_{L^2(\partial K \setminus \partial\Omega)}^2 \\
        &\quad+ \norm*{\delta^{\nicefrac{1}{2}} h_K^{\nicefrac32}q_K^{-\nicefrac32} \left( g_R - \nabla u_{hp}\cdot\vect{n}_F + ik u_{hp}\right)}_{L^2(\partial K \cap \Gamma_R)}^2
        +\norm*{\alpha^{\nicefrac{1}{2}} h_K^{\nicefrac12}q_K^{-\nicefrac12} (g_D-u_{hp})}_{L^2(\partial K \cap \Gamma_D)}^2.
        \label{eqn:error_indicator}
\end{aligned}
\end{equation}

We stress that the choice of the exponents of $h_K$ and $q_K$ have been selected 
on the basis of numerical experimentation on a problem with a smooth analytical
solution; for details, see Section~\ref{section:aposteriori_effectivity} below.
Compared to the error indicator \eqref{eqn:error_indicator_kmw} from \cite{Kapita2015}, we note that
we have a factor of $h$ instead of $k^{-1}$ in the terms with $\jmp{\nabla u_{hp}}$ and in the Robin boundary terms.
These are both dimensionally correct, but reproducing the numerical experiments conducted in Section~\ref{section:aposteriori_effectivity} with the dependency on $k^{-1}$ results in different effectivities for different wavenumbers $k$.

\subsection{Plane wave directional adaptivity}\label{sec:direction_adapt}

In this section, we discuss the design of a practical algorithm for determining
the direction vectors $\vect{d}_{K,\ell}$, $\ell=0,\dots,p_K-1$, used to define the
plane wave basis within each element
$K$ in the computational mesh $\mesh$. The key observation is that, many wave 
propagation problems typically 
exhibit a dominant direction of propagation of the underlying wave
within each element in $\mesh$. Thereby, by aligning
the plane wave basis in an appropriate fashion, we expect to attain a significant
reduction of the error in the computed TDG solution. Indeed,
in the simple case when the analytical solution is a plane wave, then if
the direction for one of the plane wave basis functions is selected such that it is aligned 
with this plane wave direction,
then the TDG method will exactly recover the analytical solution, subject to rounding errors.

The essential idea here is to simply rotate the element basis according to the
predicted elementwise dominant direction. For simplicity of presentation, let us
consider the two-dimensional case, i.e., $d=2$; we note that $d=3$ follows in 
an analogous manner, cf.~Remarks~\ref{remark:potential_direction} \& \ref{remark:direction} below. 
In two-dimensions, the standard plane wave directions are generally
selected to be evenly spaced, with the first direction $\vect{d}_{K,0} = (1,0)^\top$ always 
pointing along the $x$-axis, cf.~\eqref{eqn:plane_wave:2d} (in the
three-dimensional setting, the first direction vector typically points along the $z$-axis).
Alternatively, assuming that a dominant elementwise direction, denoted by $\actdir$, can be
determined within each $K\in\mesh$, then the direction vectors for the plane wave basis functions 
in $K$ are chosen such that the first plane wave direction is aligned with $\actdir$, 
i.e.,~\eqref{eqn:plane_wave:2d} is replaced by
\begin{equation}
\label{eqn:plane_wave:2d:rotated}
\vect{d}_{K,\ell} = (\cos(\nicefrac{2\pi\ell}{p_K} + \theta_K), \sin(\nicefrac{2\pi\ell}{p_K} + \theta_K))^\top,
\end{equation}
$\ell=0,\ldots,p_K-1$, 
where $\theta_K$ is the angle between $\actdir$ and the $x$-axis.

Clearly, in general, the dominant elementwise direction $\actdir$, $K\in\mesh$, cannot
be determined {\em a priori}, but instead must be numerically estimated as part of the
solution process. In this regard, a number of algorithms have been proposed 
within the literature; here we mention the ray-tracing approach developed 
in~\cite{Betcke2011,Betcke2012}, though this 
includes terms involving integrals over the elements within the underlying TDG formulation. In 
\cite{Amara2014}, the optimal angle of rotation was numerically estimated based on
adding an extra unknown into the problem; however, this leads to a system of nonlinear
equations to be computed. Finally, \cite{Gittelson2008} uses an approximation of
\[
    \frac{\nabla e(\vect{x}_0)}{ike(\vect{x}_0)},
\]
at a given point $\vect{x}_0\in K$, $K\in\mesh$, where $e$ denotes the error. 

Stimulated by the work undertaken in~\cite{Formaggia2001,Formaggia2003}, cf. 
also~\cite{ghh-paper,hall-thesis}, on the design of anisotropically refined
computational meshes, in this section we compute an estimate of $\actdir$, $K\in\mesh$,
based on the properties of the Hessian of the TDG solution $u_{hp}$. Indeed,
we note that the principal eigenvector, i.e., the eigenvector corresponding
to the largest eigenvalue in absolute value, of the Hessian of a given function 
indicates the direction of most rapid variation, and thereby, in our context,
the dominant direction of wave propagation. With this in mind, writing
$\Hessian(\varphi,\vect{x}_0)$ to denote the Hessian matrix of a given 
function $\varphi$, evaluated at the point $\vect{x}_0\in\real^d$, in 
Algorithm~\ref{algo:potential_direction} we outline the steps involved
in computing a {\em potential} dominant plane wave direction
$\potdir$ for a given element $K\in\mesh$.
Table~\ref{table:potential_direction} summarizes how this potential first plane wave 
direction $\potdir$ is selected; for the numerical experiments presented in 
Section~\ref{sec:numerical_examples}, we set $\Lambda =2$.
We note that in the case when no primary propagation direction
is determined, then we leave the first plane wave direction unchanged. 

\begin{algorithm}[t!]
\caption{Computation of the potential first plane wave direction $\potdir$ for element $K$.}
\label{algo:potential_direction}
\begin{algorithmic}[1]
\State Input: the TDG solution $u_{hp}$ of the discrete problem~\eqref{eqn:bilinear_form} and the parameter $\Lambda>1$.
\State Writing $\vect{x}_K$ to denote the centroid of $K$, $K\in\mesh$, evaluate the eigenpairs $(\lambda_1, \vect{v}_1), (\lambda_2, \vect{v}_2)$ of $\Hessian(\Real(u_{hp}\vert_K),\vect{x}_K)$, and $(\mu_1, \vect{w}_1), (\mu_2, \vect{w}_2)$  of $\Hessian(\Imag(u_{hp}\vert_K),\vect{x}_K)$, such that $\abs{\lambda_1} \geq \abs{\lambda_2}$ and $\abs{\mu_1} \geq \abs{\mu_2}$.
\If{$\abs{\lambda_1} \geq \Lambda \abs{\lambda_2}$}
    \If{$\abs{\mu_1} \geq \Lambda \abs{\mu_2}$}
        \If{$\abs{\lambda_1} \geq \Lambda \abs{\mu_1}$}
            \State $\potdir \gets \vect{v}_1$
        \ElsIf{$\abs{\mu_1} \geq \Lambda \abs{\lambda_1}$}
            \State $\potdir \gets \vect{w}_1$
        \Else
            \State $\potdir \gets \frac{\vect{v}_1+\vect{w}_1}{\norm{\vect{v}_1+\vect{w}_1}}$
        \EndIf
    \Else
        \If{$\abs{\lambda_1} \geq \Lambda \abs{\mu_1}$}
            \State $\potdir \gets \vect{v}_1$
        \Else
            \State No primary propagation direction
        \EndIf
    \EndIf
\Else
    \If{$\abs{\mu_1} \geq \Lambda \abs{\mu_2}$}
        \If{$\abs{\mu_1} \geq \Lambda \abs{\lambda_1}$}
            \State $\potdir \gets \vect{w}_1$
        \Else
            \State No primary propagation direction
        \EndIf
    \Else
        \State No primary propagation direction
    \EndIf
\EndIf
\end{algorithmic}
\end{algorithm}

\begin{remark}\label{remark:potential_direction}
We note that in the case when $d=3$, $\Hessian(\Real(u_{hp}\vert_K),\vect{x}_K)$
and $\Hessian(\Imag(u_{hp}\vert_K),\vect{x}_K)$ each have a third eigenpair, $(\lambda_3,\vect{v}_3)$ 
and $(\mu_3, \vect{w}_3)$, respectively. 
However, if the eigenpairs are sorted such that 
$\abs{\lambda_1} \geq \abs{\lambda_2} \geq \abs{\lambda_3}$ and 
$\abs{\mu_1} \geq \abs{\mu_2} \geq \abs{\mu_3}$, the third eigenpairs 
\emph{never} represent a dominant direction, and thereby 
Algorithm~\ref{algo:potential_direction} can be used to identify 
$\potdir$, $K\in\mesh$, without modification.
\end{remark}

\begin{table}[t!]
\centering
\begin{tabular}{c|c|c|c||c}
$\abs{\lambda_1} \geq C \abs{\lambda_2}$ & $\abs{\mu_1} \geq C\abs{\mu_2}$ & $\abs{\lambda_1} \geq C\abs{\mu_1}$ & $\abs{\mu_1}\geq C\abs{\lambda_1}$ & First Plane Wave $\potdir$  \\\hline\hline
\cmark & \cmark & \cmark & \xmark & $\vect{v}_1$ \\
\cmark & \cmark & \xmark & \cmark & $\vect{w}_1$ \\
\cmark & \cmark & \xmark & \xmark & $\frac{(\vect{v}_1+\vect{w}_1)}{\norm{\vect{v}_1+\vect{w}_1}}$ \\
\cmark & \xmark & \cmark & \xmark & $\vect{v}_1$ \\
\cmark & \xmark & \xmark &  ---   & --- \\
\xmark & \cmark & \xmark & \cmark & $\vect{w}_1$ \\
\xmark & \cmark &  ---   & \xmark & --- \\
\xmark & \xmark &  ---   &  ---   & --- 
\end{tabular}
\caption{Summary of selection of first plane wave direction $\potdir$ using Algorithm~\ref{algo:potential_direction}.}
\label{table:potential_direction}
\end{table}

Noting that eigenvectors are only unique up to scalar multiples, 
the vector $\potdir$, $K\in\mesh$, evaluated according to 
Algorithm~\ref{algo:potential_direction} may be pointing in precisely the
opposite direction to the primary wave propagation direction. Thereby,
to ensure that $\potdir$, $K\in\mesh$, is correctly oriented,
we study the impedance trace on the boundary of a ball $B_\delta(\vect{x}_K)$ 
of radius $\delta$, centered  at $\vect{x}_K$, of both the numerical solution and 
a plane wave with (the desired) propagation direction $\actdir$. As we let
$\delta\to 0$, we expect that the numerical solution should be closely approximated by 
the plane wave in the primary propagation direction. 

More precisely, given $K\in\mesh$, the impedance trace of the plane wave
\[
    \tilde{u}_K(\vect{x}) = \e^{ik\actdir \cdot(\vect{x}-\vect{x}_K)}
\]
on $\partial B_\delta(\vect{x}_K)$ is given by
\begin{eqnarray}
	(\nabla\tilde{u}_K(\vect{x})\cdot\vect{n}_{B_\delta} + ik\tilde{u}_K(\vect{x}))
	|_{\partial B_\delta(\vect{x}_K)}
&=& (ik(\actdir \cdot\vect{n}_{B_\delta}+1)\,\e^{ik \actdir \cdot(\vect{x}-\vect{x}_K)}
)|_{\partial B_\delta(\vect{x}_K)},
\label{eqn:impedance_trace}
\end{eqnarray}
where $\vect{n}_{B_\delta}$ denotes the unit outward normal vector on 
$\partial B_\delta(\vect{x}_K)$.
Setting $\vect{x}=\vect{x}_K + \delta \potdir$ in \eqref{eqn:impedance_trace}  
and noting that, at this point of evaluation, $\vect{n}_{B_\delta} = \potdir$, 
we deduce that 
\[
\frac{\nabla\tilde{u}_K(\vect{x}_K + \delta\potdir )\cdot\vect{n}_{B_\delta} 
+ ik\tilde{u}_K(\vect{x}_K + \delta \potdir )}{ik} =
\begin{cases}
    2 \e^{ik\delta}, & \text{if } \potdir = \actdir , \\
    0, & \text{if } \potdir = -\actdir .
\end{cases}
\]
Thereby, the (potential) dominate direction of propagation $\potdir$, $K\in\mesh$, predicted 
according to Algorithm~\ref{algo:potential_direction} may be corrected to yield 
the dominant direction $\actdir$ on the basis of Algorithm~\ref{algo:direction};
this direction will then be selected as the first plane wave direction on element $K$,
$K\in\mesh$. For simplicity, throughout this article we set $\delta=0$.

\begin{algorithm}[t!]
\caption{Evaluation of the first plane wave direction $\actdir$ for element $K$.}
\label{algo:direction}
\begin{algorithmic}[1]
\State Input: the TDG solution $u_{hp}$ of the discrete problem~\eqref{eqn:bilinear_form}, the parameter $0\leq \delta\to 0$, and $\potdir$ computed by Algorithm~\ref{algo:potential_direction}. 
\State The first plane wave direction $\actdir$ on element $K$, $K\in\mesh$, is given by
\[
    \actdir = \begin{cases}
        -\potdir , & ~\mbox{ if } ~\Real\left(\frac{\nabla u_{hp}(\vect{x}_K + \delta\potdir)\cdot\potdir + ik u_{hp}(\vect{x}_K + \delta\potdir)}{ik}\right) < \e^{ik\delta}, \\
        \potdir , & ~\mbox{ if } ~\Real\left(\frac{\nabla u_{hp}(\vect{x}_K + \delta\potdir )\cdot\potdir + ik u_{hp}(\vect{x}_K + \delta\potdir)}{ik}\right) \geq \e^{ik\delta}. \\
    \end{cases}
\]
\end{algorithmic}
\end{algorithm}

\begin{remark}\label{remark:direction}
In the three--dimensional setting, once the selection of the primary wave propagation direction 
$\actdir$ has been computed on the basis of 
Algorithms~\ref{algo:potential_direction} \&~\ref{algo:direction},
we then select the remaining wave directions,
$\vect{d}_{K,\ell}$, $\ell=1,\dots,p_K-1$, by applying a transformation matrix 
$T\in\real^{3\times 3}$ to the original `reference' directions $\vect{\tilde{d}}_{K,\ell}$, 
$\ell=1,\dots,p_K-1$, respectively, where $\vect{\tilde{d}}_{K,0}$ points along the
$z$-axis, cf. above. Thereby, 
$$
\vect{d}_{K,\ell}=T \vect{\tilde{d}}_{K,\ell} ,
$$
$\ell=1,\dots,p_K-1$, where $T$ is selected such that
\[
\actdir \equiv \vect{d}_{K,0} = T \begin{pmatrix} 0 \\ 0 \\ 1 \end{pmatrix} 
\equiv T \vect{\tilde{d}}_{K,0}.
\]
We note that the selection of $T$ is not unique; writing 
$\actdir = (d_x, d_y, d_z)^\top$, we 
define $T$ to be the identity matrix if $d_x=d_y=0$; otherwise, we set
\[
T = \begin{pmatrix}
\frac{d_x d_z}{\sqrt{d_x^2+d_y^2}} & \frac{d_y}{\sqrt{d_x^2+d_y^2}} & d_x \\
\frac{d_y d_z}{\sqrt{d_x^2+d_y^2}} & -\frac{d_x}{\sqrt{d_x^2+d_y^2}} & d_y \\
-\sqrt{d_x^2+d_y^2} & 0 & d_z
\end{pmatrix}.
\]
\end{remark}

\subsection{$hp$--Adaptive mesh refinement}\label{sec:hp}
In this section we discuss the design of an automatic algorithm for generating
sequences of $hp$--adaptively refined TDG finite element spaces in an efficient
manner. This topic has been extensively studied within the finite element
element literature in the case when the local element spaces consist of
polynomial functions; for a comprehensive review, we refer to~\cite{Mitchell2011, Mitchell2014}.
In general, the key underlying principle of most $hp$--refinement strategies is 
to employ local mesh subdivision ($h$--refinement) in regions where the solution is
not smooth, while local enrichment of the finite element space ($p$--refinement)
is undertaken elsewhere. Given that such regularity information is generally
unknown {\em a priori}, several strategies have been developed to {\em a posteriori} estimate
the local smoothness of the analytical solution, based on its numerical approximation;
cf. \cite{Houston2005}, for example. However, in the context of TDG schemes for
the numerical approximation of high-frequency time-harmonic wave problems, the extraction
of such regularity information is expected to be unreliable due to
the oscillatory nature of the computed numerical solution.

Thereby, as an alternative to directly estimating local smoothness of the solution,
we employ the {\em a posteriori} error indicator~\eqref{eqn:error_indicator}
to select the type of refinement to be undertaken on the basis of the refinement
history of the current element, cf.~\cite{Melenk2001}. More precisely, following
\cite{Melenk2001} refinements are selected based on checking if the local error estimate has 
decayed according to the expected rate of convergence based on the last type of refinement
employed. If the expected rate of convergence is achieved, 
then $p$--refinement is performed; otherwise, 
$h$--refinement is undertaken. The variant of \cite[Algorithm 4.4]{Melenk2001} we employ here
is summarized in Algorithm~\ref{algo:refinement}.
Here, we note that $\gamma_h$, $\gamma_p$, 
and $\gamma_n$ are control parameters;
for the purposes of this article, we select $\gamma_h=4$, $\gamma_p=0.4$, and $\gamma_n=1$.
Furthermore, the number of child elements, $N$, cf. step 10. in Algorithm~\ref{algo:refinement},
is dependent on the type of subdivision, i.e., isotropic/anisotropic, undertaken,
as well as the element shape; for isotropic refinement of tensor-product elements, we have
that $N=2^d$. 

\begin{algorithm}[t!]
\begin{algorithmic}[1] 
\caption{$hp$--Adaptive refinement algorithm.}\label{algo:refinement}
\State Input the control parameters $\gamma_h$, $\gamma_p$, and $\gamma_n$.
\State{Choose a coarse initial mesh~$\mesh[h,0]$ of~$\Omega$ and a corresponding low-order starting (effective) polynomial degree vector~$\vect{q}_0$, together with the total dimension vector $\vect{p}_0$ defined as in \eqref{eqn:effective_poly_deg}.}
\State{Set the initial predicted error indicator $\eta_{K,0}^{\mathrm{pred}}=\infty$ for all $K\in\mesh[h,0]$.}
\For{$i=0,1,\ldots,$ until sufficiently many iterations have been performed.}
	\State{Solve \eqref{eqn:bilinear_form} for $u_{hp} \in V_{\vect{p}_i}(\mesh[h,i])$.}
	\State{Compute the {\em a posteriori} error indicators $\eta_{K,i}\equiv \eta_K$, 
		$K\in \mesh[h,i]$, and mark elements for refinement based on their relative magnitude.}
	\For{$K\in\mesh[h,i]$}
    	\If{$K$ is marked for refinement}
        	\If{$\eta_{K,i} > \eta_{K,i}^{\mathrm{pred}}$}
            	\State Perform $h$--refinement: Subdivide $K$ into $N$ children $K_s, s = 1,\dots,N$, and set
            	\State $(\eta_{K_s,i+1}^{\mathrm{pred}})^2\gets \frac{1}{N}\gamma_h \left(\frac{1}{2}\right)^{2q_K} \eta_{K,i}^2$, $1\leq s \leq N$.            
  	      	\Else
    	        \State Perform $p$--refinement: $q_K \gets q_K+1$
        	    \State $(\eta_{K,i+1}^{\mathrm{pred}})^2\gets \gamma_p \eta_{K,i}^2$        
        	\EndIf
    	\Else
        	\State $(\eta_{K,i+1}^{\mathrm{pred}})^2\gets \gamma_n (\eta_{K,i}^{\mathrm{pred}})^2$
    	\EndIf
	\EndFor
	\State Construct the new mesh $\mesh[h,i+1]$ and corresponding Trefftz space $V_{\vect{p}_{i+1}}(\mesh[h,i+1])$.
\EndFor
\end{algorithmic}
\end{algorithm}

\begin{remark}
We note that in \cite{Melenk2001} the initial values of the predicted error indicator 
$\eta_{K,0}^{\mathrm{pred}}$, $K\in\mesh[h,0]$, are set to zero; thereby, this ensures that 
$h$--refinement is undertaken the first time an element is refined. In contrast, in
Algorithm~\ref{algo:refinement} we set $\eta_{K,0}^{\mathrm{pred}}=\infty$ for all 
$K\in\mesh[h,0]$ which instead leads to $p$--enrichment being undertaken as the first refinement
of a given element, since the TDG method for the numerical approximation of
the Helmholtz equation is intrinsically a high-order method.
\end{remark}

\begin{remark}
Plane wave directional adaptivity can be performed at different stages within
Algorithm~\ref{algo:refinement}; for example, the following options are available:
\begin{itemize}
\item undertake directional adaptivity only on elements marked for $p$--refinement,
\item undertake directional adaptivity on all elements marked for refinement, with $h$--refinement performed after plane wave direction adaptivity, or
\item undertake directional adaptivity on every element $K\in\mesh$, even if the element $K$ has not been marked for refinement.
\end{itemize}
In Section~\ref{sec:numerical_examples} we shall numerically investigate each of these approaches in order to assess their relative computational performance in terms of error reduction.
\end{remark}

\begin{remark}
As a final remark, we note that within Algorithm~\ref{algo:refinement} we employ the 
fixed fraction refinement strategy to select elements for refinement, cf. step 6; throughout
this article we set the refinement fraction equal to $25\%$.
\end{remark}

\section{Numerical experiments} \label{sec:numerical_examples}

In this section, we present a series of numerical experiments to highlight the
practical performance of the $hp$--refinement algorithm, with directional adaptivity,
proposed in Algorithm~\ref{algo:refinement}. Throughout this section we shall compare 
the performance of the proposed $hp$--adaptive refinement strategy with the corresponding 
algorithm based on exploiting only local mesh subdivision, i.e., $h$--refinement.
The numerical experiments presented within this section have been undertaken using the
AptoFEM software package~\cite{aptofem}.

\subsection{Plane wave direction adaptivity}\label{section:plane_wave_refine}
\begin{table}[pt]
\centering
\begin{tabular}{c|c|r@{.}l@{$\times$}l|r@{.}l@{$\times$}l|c}
&  & \multicolumn6{c}{Relative $L^2(\Omega)$-Error} & \\
$q$ & No of Dofs & \multicolumn3{c|}{Standard TDG} & \multicolumn3c{Direction Adaptivity} 
& \% Reduction\\ \hline
3 & 112 & \quad 2&015 & $10^{0}$ & \quad 1&959 & $10^{0}$ & 2.7\% \\
4 & 144 & 5&027 & $10^{-1}$ & 3&194 & $10^{-1}$ & 36.5\% \\
5 & 176 & 7&414 & $10^{-2}$ & 2&658 & $10^{-2}$ & 64.1\% \\
6 & 208 & 1&616 & $10^{-2}$ & 6&320 & $10^{-3}$ & 60.9\% \\
7 & 240 & 3&420 & $10^{-3}$ & 1&435 & $10^{-3}$ & 58.0\% \\
8 & 272 & 5&154 & $10^{-4}$ & 3&011 & $10^{-4}$ & 41.6\% \\
9 & 304 & 8&928 & $10^{-5}$ & 6&908 & $10^{-5}$ & 22.6\% \\
\end{tabular}
\caption{Plane Wave Refinement: Comparison of the relative $L^2$-error for uniform $p$--refinement (without direction adaptivity), and $p$--refinement with direction adaptivity (Algorithm~\ref{algo:direction}).}
\label{table:pw_adapt}
\end{table}

\begin{figure}[pt]
    \configfigure
    \subfloat[]{\label{fig:directions:p3}\includegraphics[width=0.4\textwidth]{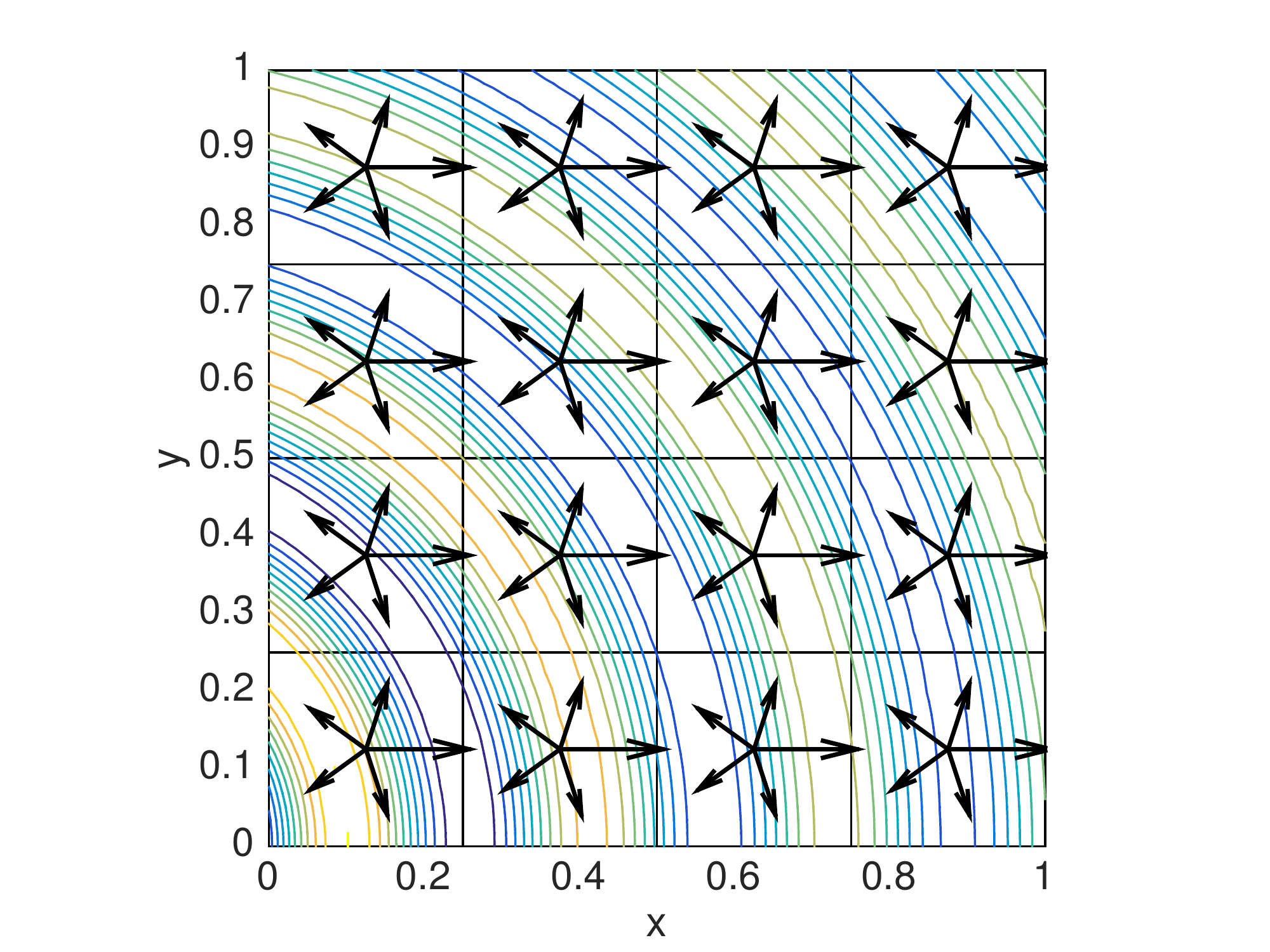}}
    \subfloat[]{\label{fig:directions:p4}\includegraphics[width=0.4\textwidth]{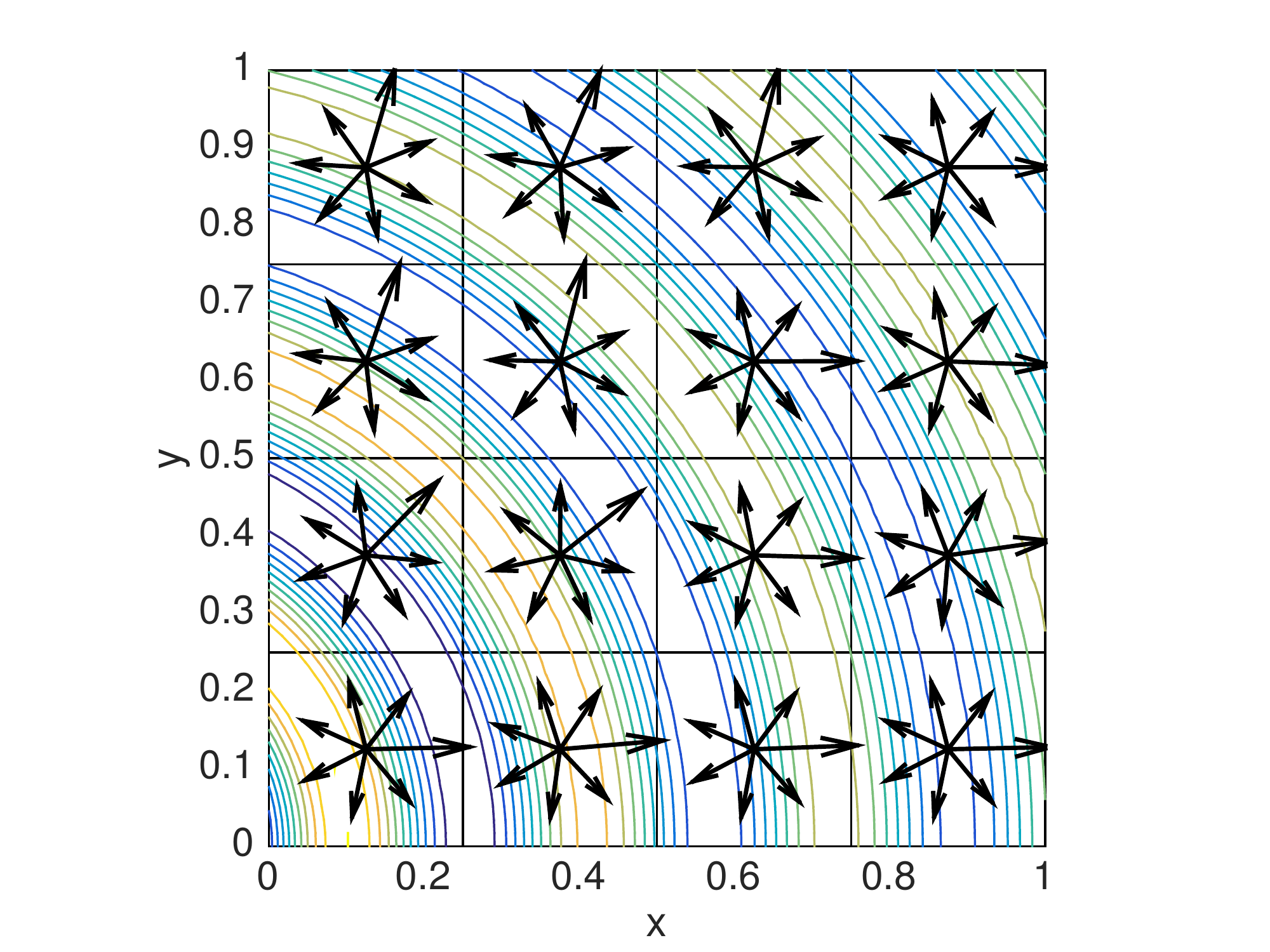}} \\
    \subfloat[]{\label{fig:directions:p5}\includegraphics[width=0.4\textwidth]{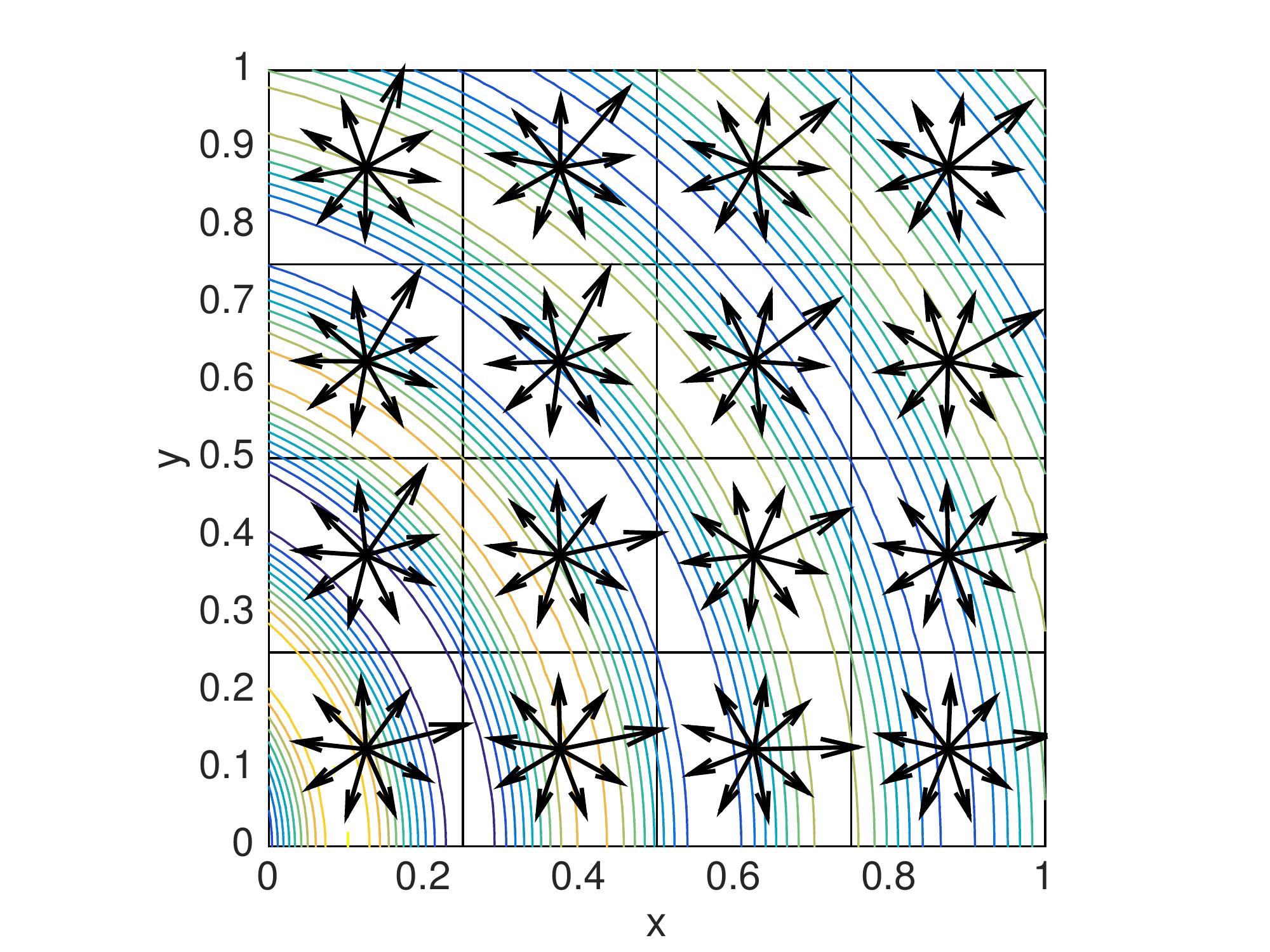}}
    \subfloat[]{\label{fig:directions:p6}\includegraphics[width=0.4\textwidth]{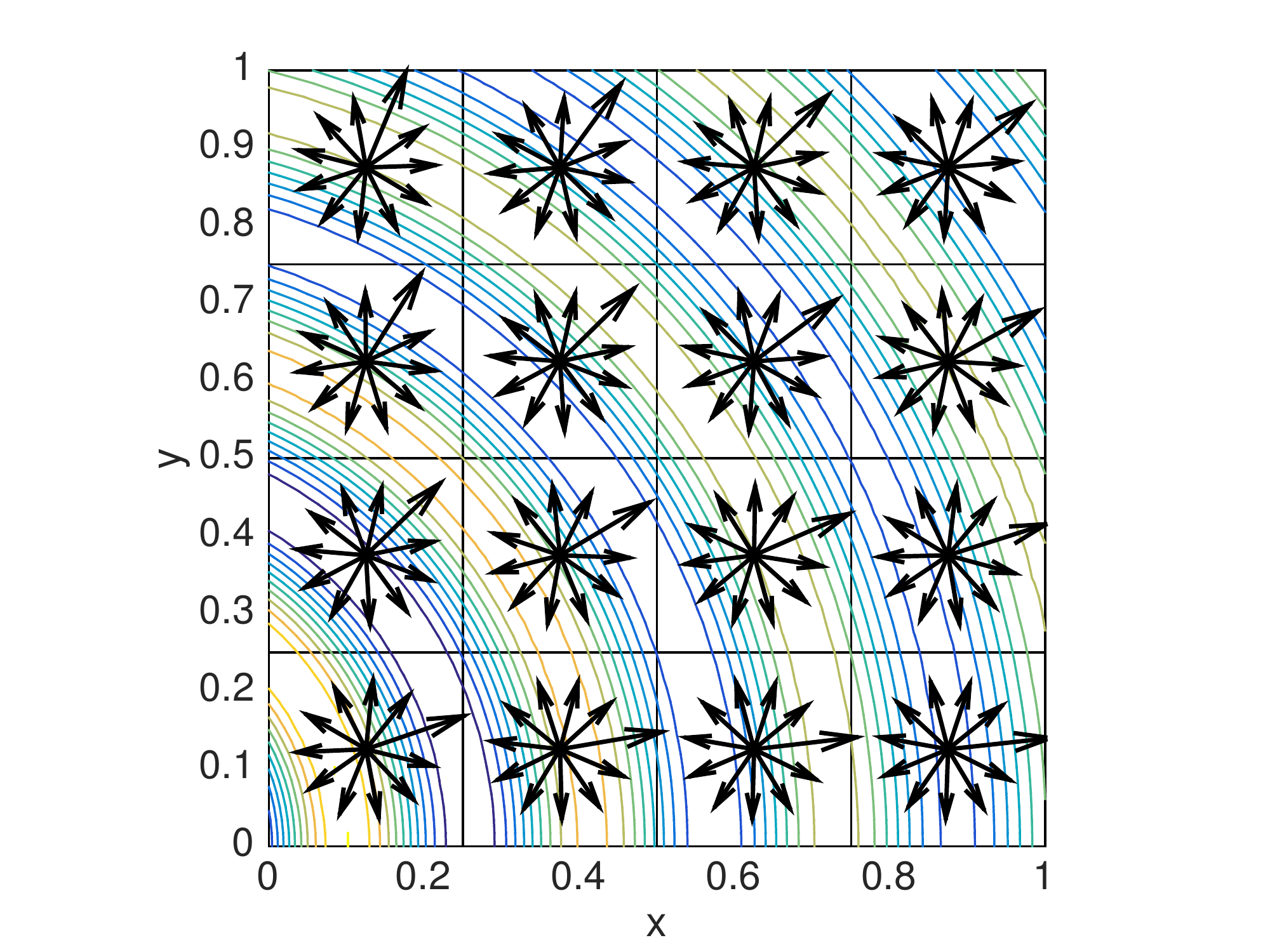}}
    \caption{Plane Wave Refinement: Plane wave directions of \protect\subref{fig:eff:20} initial mesh and after \protect\subref{fig:eff:30} $1$, \protect\subref{fig:eff:40} $2$ and \protect\subref{fig:eff:50} $3$ $p$--refinements with plane wave refinement (Algorithm~\ref{algo:direction})}
    \label{fig:directions}
\end{figure}

In this first example, we study the effect of adjusting the plane wave directions while employing a fixed computational mesh with uniform $p$--refinement. To this end, we consider problem~\eqref{eqn:helmholtz} with $\Omega=(0,1)^2$, $\Gamma_R=\partial\Omega$, and $\Gamma_D\equiv\emptyset$; furthermore, the Robin boundary condition $g_R$ is selected such that the analytical solution $u$ of~\eqref{eqn:helmholtz} is given by
\begin{equation}
\label{eqn:hankal_anal}
u(x,y) = \mathcal{H}_0^{(1)} \left( k \sqrt{ (x+\nicefrac14)^2+y^2 } \right),
\end{equation}
where $\mathcal{H}_0^{(1)}$ denotes the Hankel function of the first kind of order 0. Throughout this section, we set $k=20$; note that for this problem the analytical solution $u$ is smooth in $\Omega$.
 
Here, the underlying computational mesh consists of 16 uniform square elements; on each element we 
initially select the effective polynomial degree $q=2$, i.e., $p=5$. In Table~\ref{table:pw_adapt} we 
compare the computed relative $L^2$-error based on employing uniform $p$--refinement of the underlying TDG
space $V_{\vect{p}}(\mesh)$ in the two cases when the standard TDG scheme is employed, i.e., when the local 
plane wave directions are kept fixed, and when plane wave directional adaptivity is utilised, based on exploiting
Algorithm~\ref{algo:direction} (direction adaptivity). We note that, since uniform $p$--refinement 
is employed in both cases, then at each step of the refinement, both schemes possess the same 
number of degrees of freedom. At each step of the refinement algorithm, we observe that the 
exploitation of directional adaptivity leads to roughly 50\% reduction in the relative 
$L^2$-error when compared to the corresponding quantity computed for the standard TDG
method (without direction adaptivity). We note, however, that in the case when $q=3$, 
the relative $L^2$-error is only reduced by a small amount when directional adaptivity is employed; this
is due to the fact that the local plane wave directions are computed based on the numerical solution evaluated with $q=2$,
which is numerically too inaccurate to reliably predict the correct local direction of wave propagation.
Furthermore, we also note that, as the number of plane waves increases, the improvement in the 
relative $L^2$-error decreases; this is caused by the fact that, as the number of plane waves 
increases for the standard TDG scheme, one of the directions will get closer to the actual 
dominant direction.

In Figure~\ref{fig:directions} we plot, for each element, the initial plane wave directions 
and the plane wave directions computed after $1$, $2$, and $3$ uniform $p$--refinements employing 
directional adaptivity. We emphasize the first plane wave direction with a larger arrow, i.e., 
the dominant wave direction as determined by Algorithm~\ref{algo:direction}. Moreover,
we overlay the directions on top of a contour plot showing the real part of the analytical 
solution~\eqref{eqn:hankal_anal}. From Figure~\ref{fig:directions}, we can clearly observe that
the directional adaptivity algorithm is able to accurately determine the dominant wave direction after 
a few refinements.

\begin{table}[pt]
\centering
\begin{tabular}{c|c|r@{.}l@{$\times$}l|r@{.}l@{$\times$}l|r@{.}l@{$\times$}l}
&  & \multicolumn9{c}{Relative $L^2(\Omega)$-Error} \\
$q$ & No of Dofs & \multicolumn3{c|}{Initial} & \multicolumn3c{One Direction Adapt.} & \multicolumn3c{Two Direction Adapts.} \\ \hline
3 & 112 & 2&015 & $10^{0}$ & \qquad 8&755 & $10^{-1}$ & \qquad 5&856 & $10^{-1}$ \\
4 & 144 & 5&027 & $10^{-1}$ & 1&267 & $10^{-1}$ & 1&149 & $10^{-1}$ \\
5 & 176 & 7&414 & $10^{-2}$ & 2&614 & $10^{-2}$ & 2&584 & $10^{-2}$ \\
6 & 208 & 1&616 & $10^{-2}$ & 6&330 & $10^{-3}$ & 6&327 & $10^{-3}$ \\
7 & 240 & 3&420 & $10^{-3}$ & 1&435 & $10^{-3}$ & 1&435 & $10^{-3}$ \\
8 & 272 & 5&154 & $10^{-4}$ & 3&011 & $10^{-4}$ & 3&011 & $10^{-4}$ \\
\end{tabular}
\caption{Plane Wave Refinement: Comparison of the relative $L^2$-error with fixed effective polynomial degree, $q=3,\dots,8$, and direction adaptivity (Algorithm~\ref{algo:direction}).}
\label{table:pw_adapt_fixed_q}
\end{table}

Finally, in this section we consider performing more than one directional adaptivity step
after each uniform $p$--refinement. To this end, in Table~\ref{table:pw_adapt_fixed_q} we compare 
the relative $L^2$-error for the initial directions, as well as after one and two steps of 
directional adaptivity have been performed, for the case when $q=3,\dots,8$. Here, we observe that
additional application of the direction adaptivity algorithm does not lead to a significant
reduction in the relative $L^2$-error; indeed, most of the reduction, when compared to the standard
TDG scheme, without directional adaptivity, is attained after one step of 
Algorithm~\ref{algo:direction}. Moreover, we emphasise that this first step may be undertaken in 
a very computationally cheap manner.

\begin{figure}[t]
    \configfigure
    \subfloat[$k=20$]{\label{fig:eff:20}\includegraphics[width=0.4\textwidth]{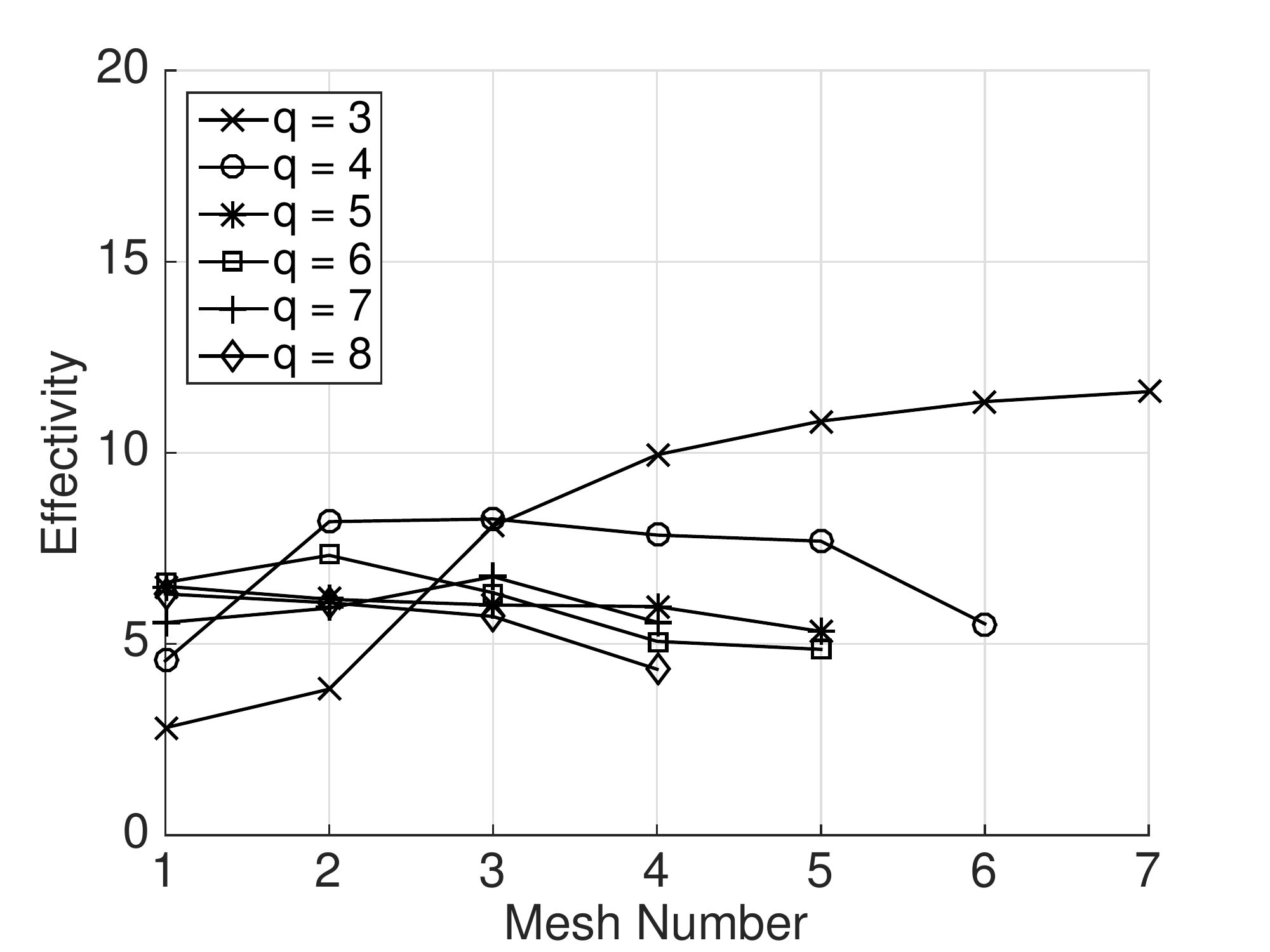}}
    \subfloat[$k=30$]{\label{fig:eff:30}\includegraphics[width=0.4\textwidth]{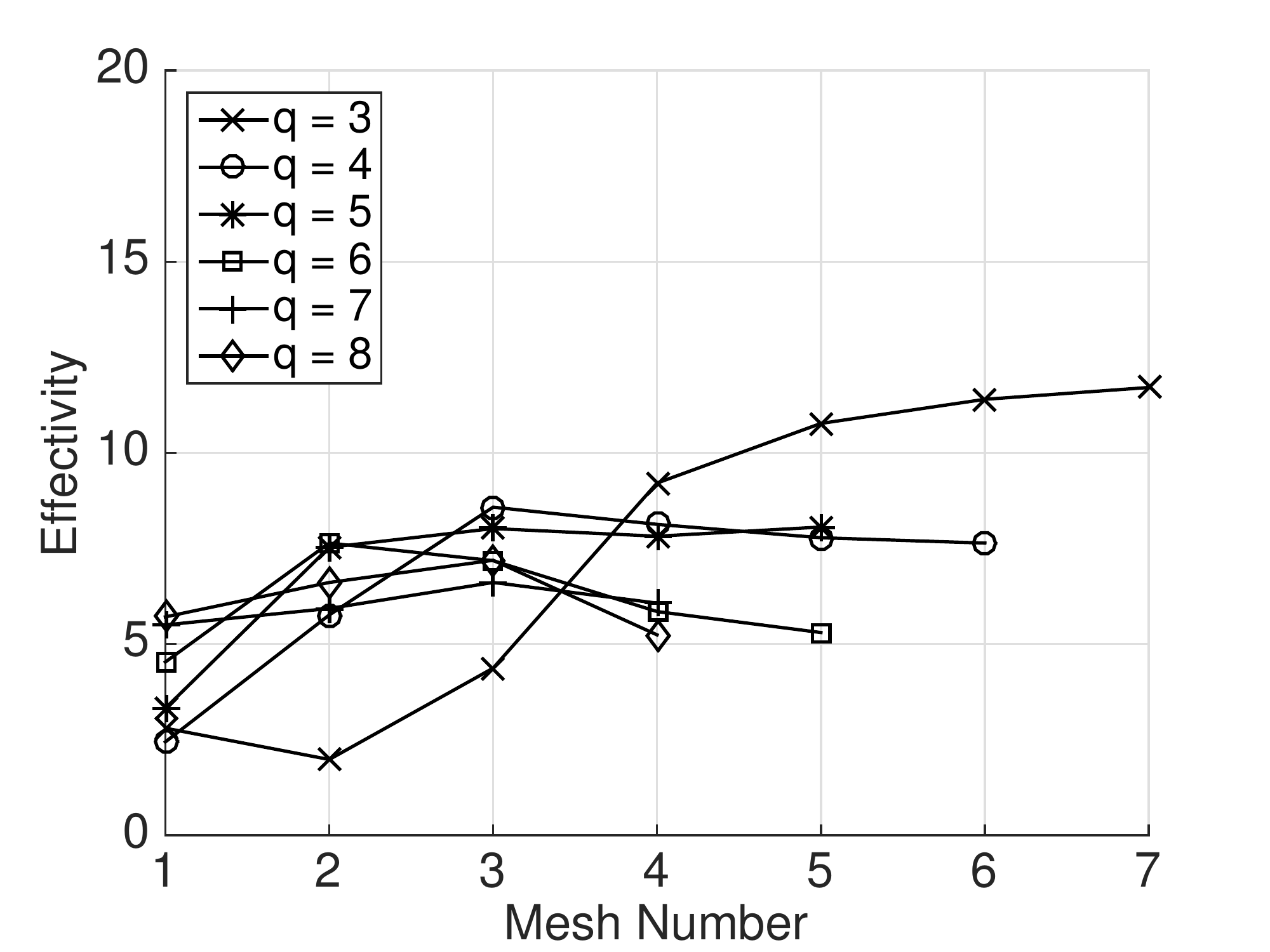}} \\
    \subfloat[$k=40$]{\label{fig:eff:40}\includegraphics[width=0.4\textwidth]{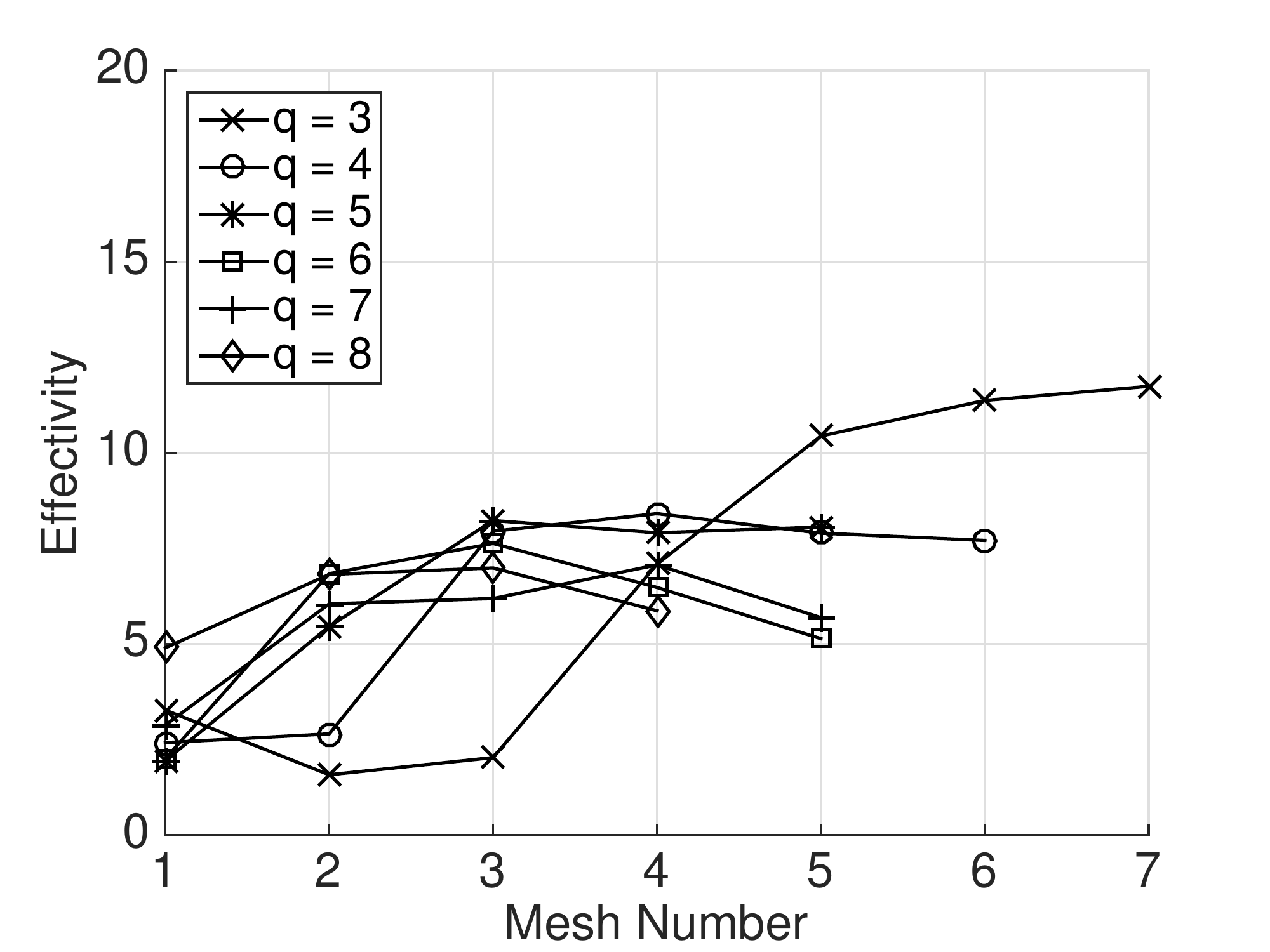}}
    \subfloat[$k=50$]{\label{fig:eff:50}\includegraphics[width=0.4\textwidth]{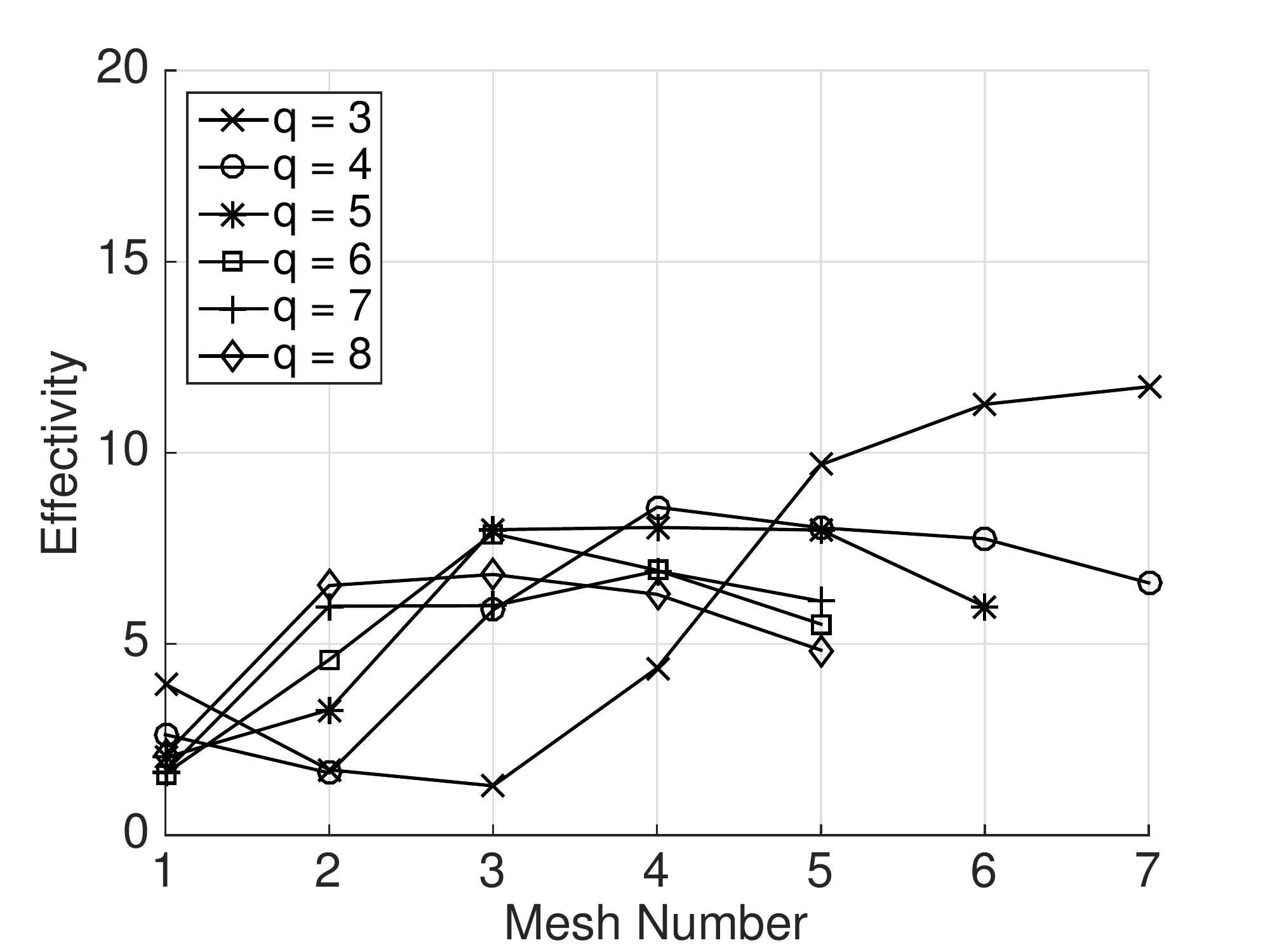}}
    \caption{Effectivities for $h$--refinement with fixed effective polynomial degree of the smooth analytical Hankel solution with different wavenumbers.}
    \label{fig:eff}
\end{figure}
\begin{figure}[tp]
    \configfigure
    \subfloat[$k=20$, $\mathcal{E}_{\jmp{u_{hp}}}$]{\label{fig:eff:u:20}\includegraphics[width=0.33\textwidth]{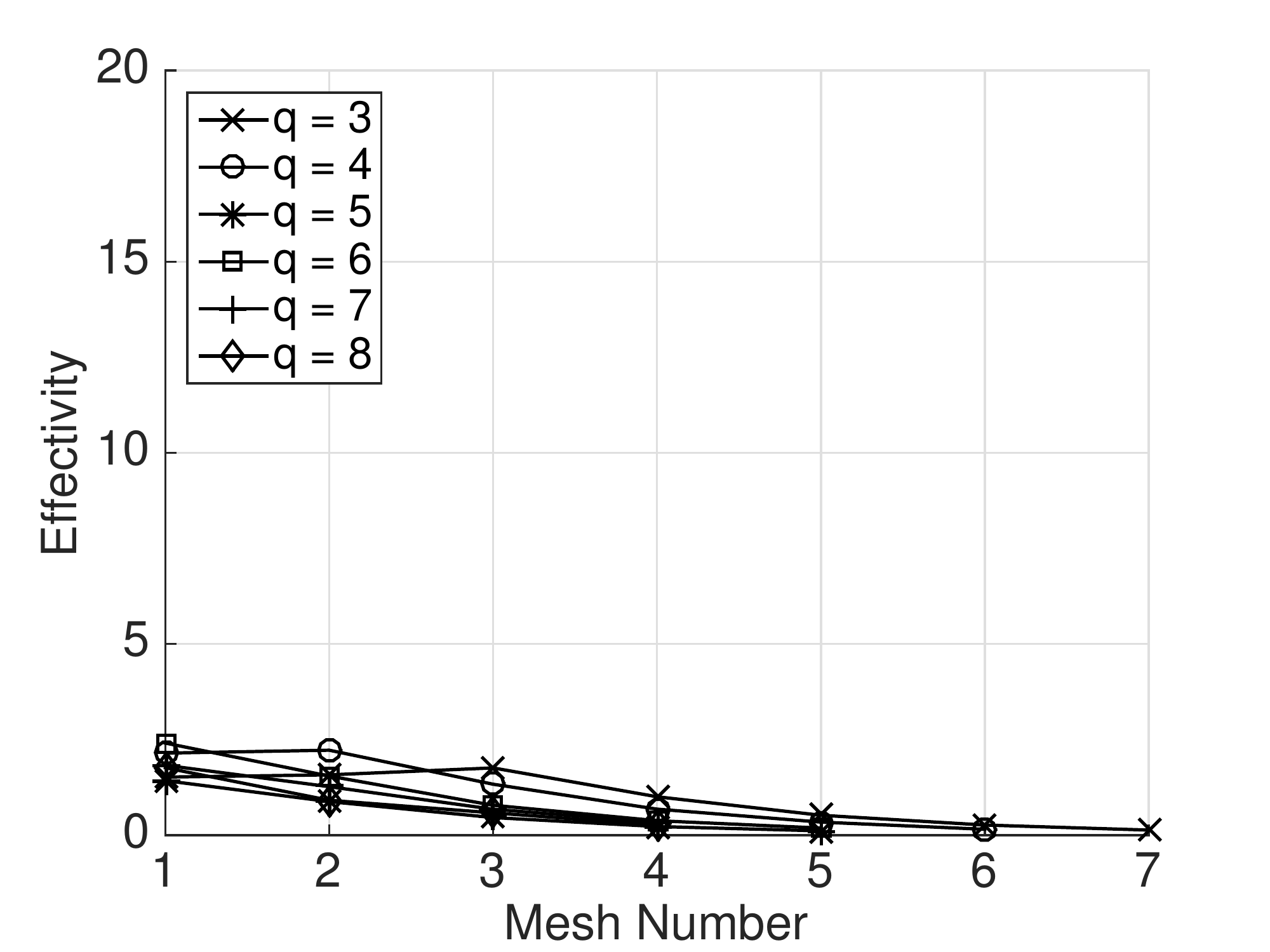}}
    \subfloat[$k=20$, $\mathcal{E}_{\jmp{\nabla u_{hp}}}$]{\label{fig:eff:gradu:20}\includegraphics[width=0.33\textwidth]{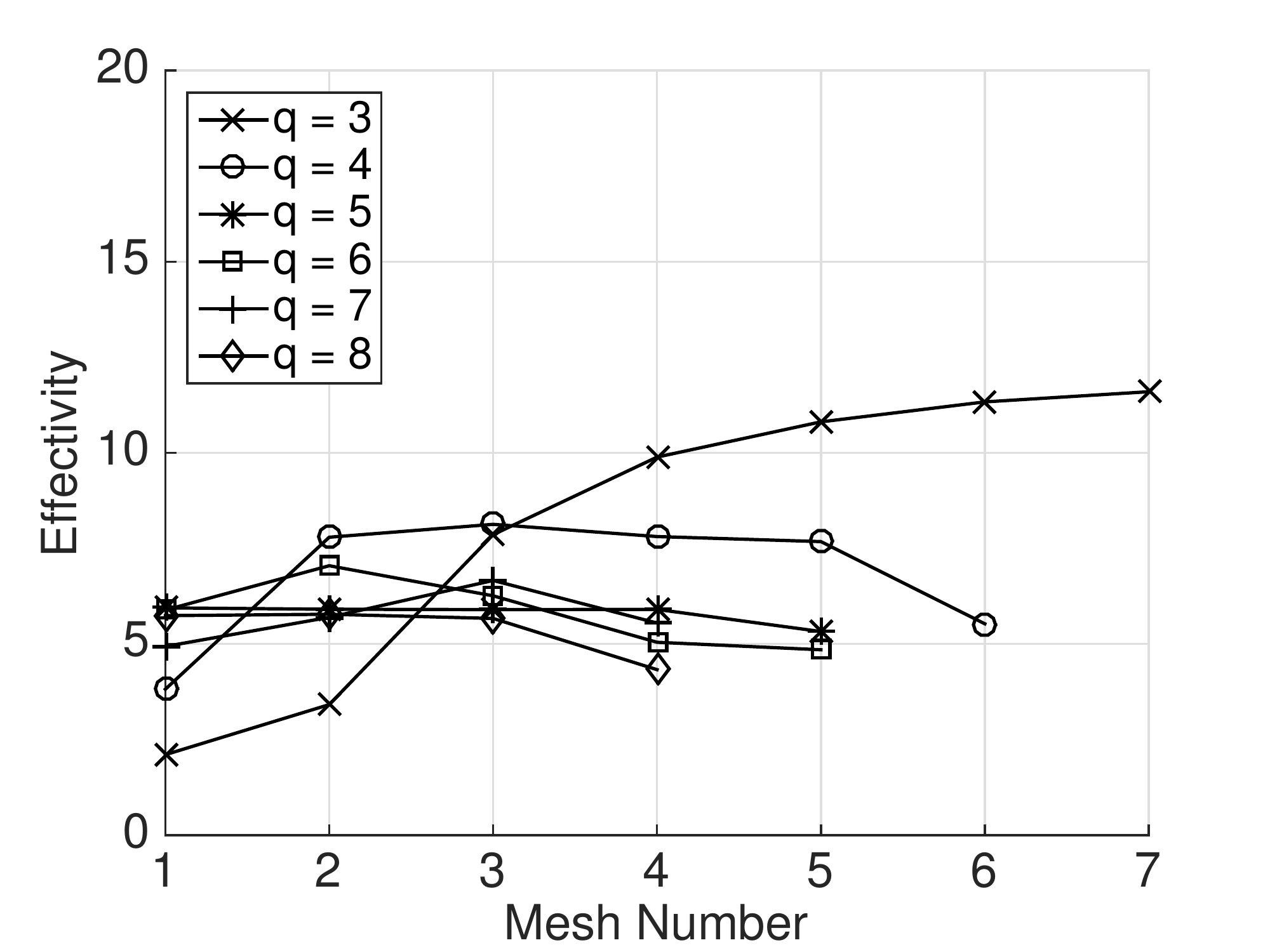}}
    \subfloat[$k=20$, $\mathcal{E}_R$]{\label{fig:eff:robin:20}\includegraphics[width=0.33\textwidth]{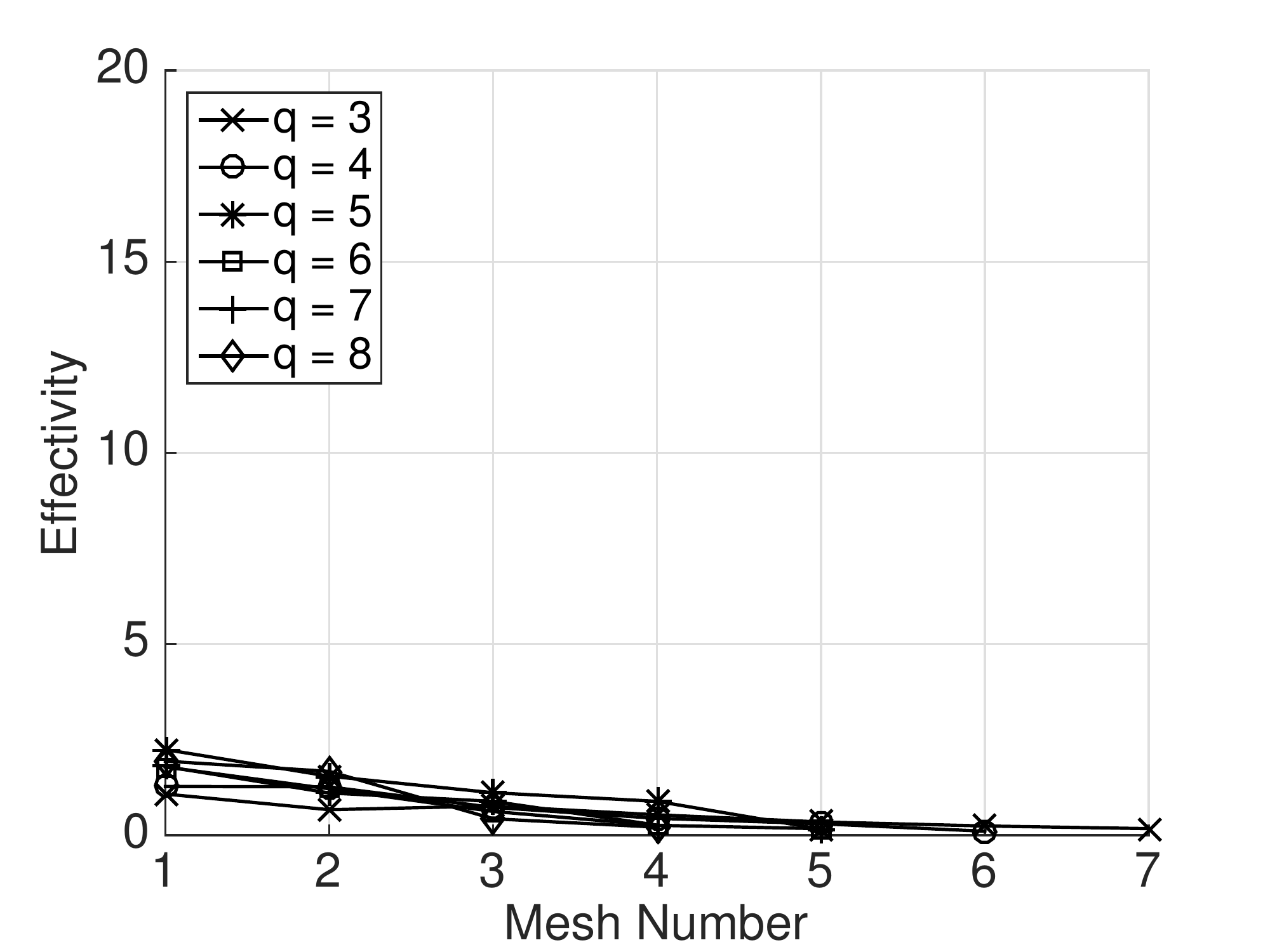}} \\
    \subfloat[$k=30$, $\mathcal{E}_{\jmp{u_{hp}}}$]{\label{fig:eff:u:30}\includegraphics[width=0.33\textwidth]{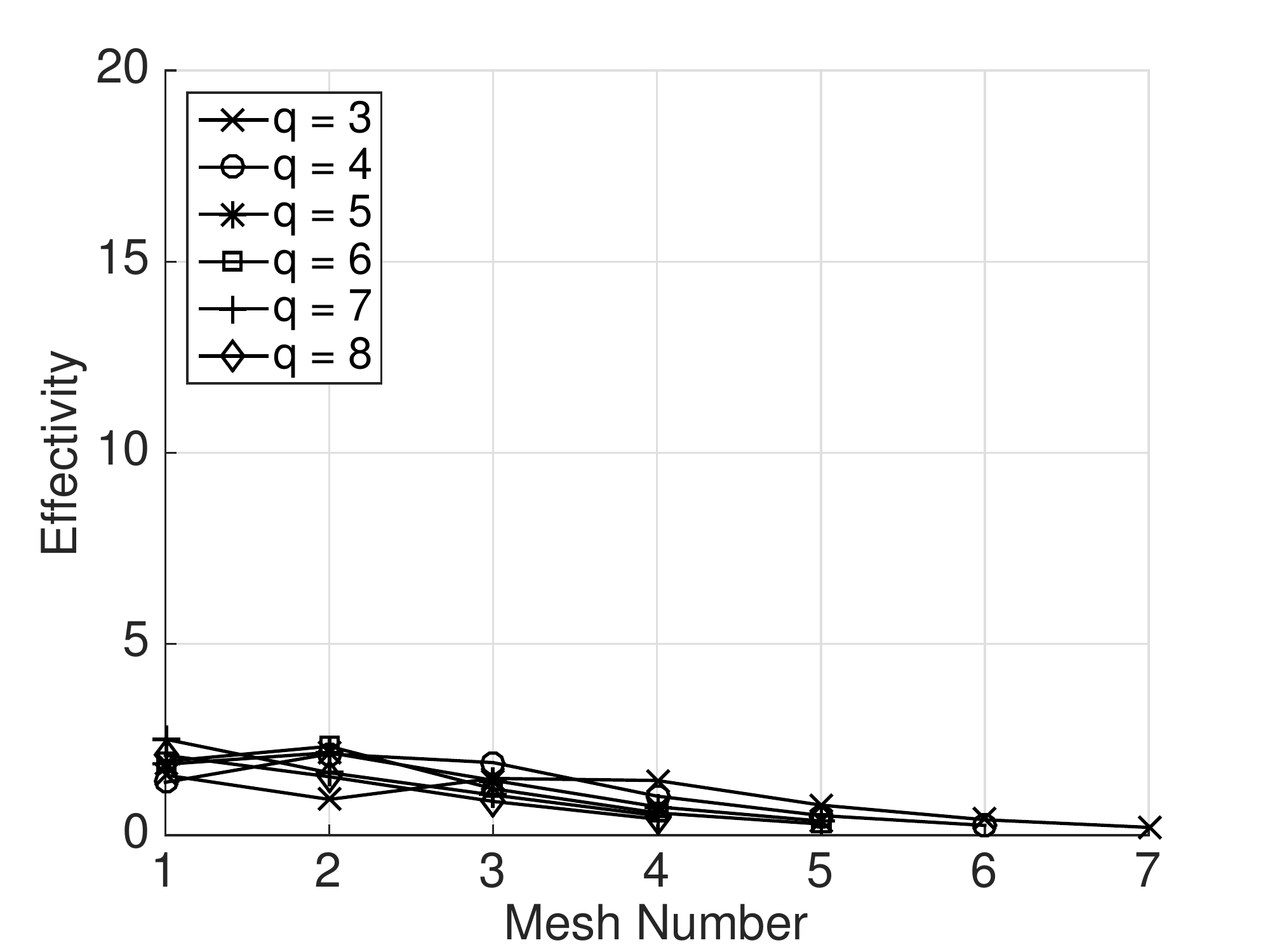}}
    \subfloat[$k=30$, $\mathcal{E}_{\jmp{\nabla u_{hp}}}$]{\label{fig:eff:gradu:30}\includegraphics[width=0.33\textwidth]{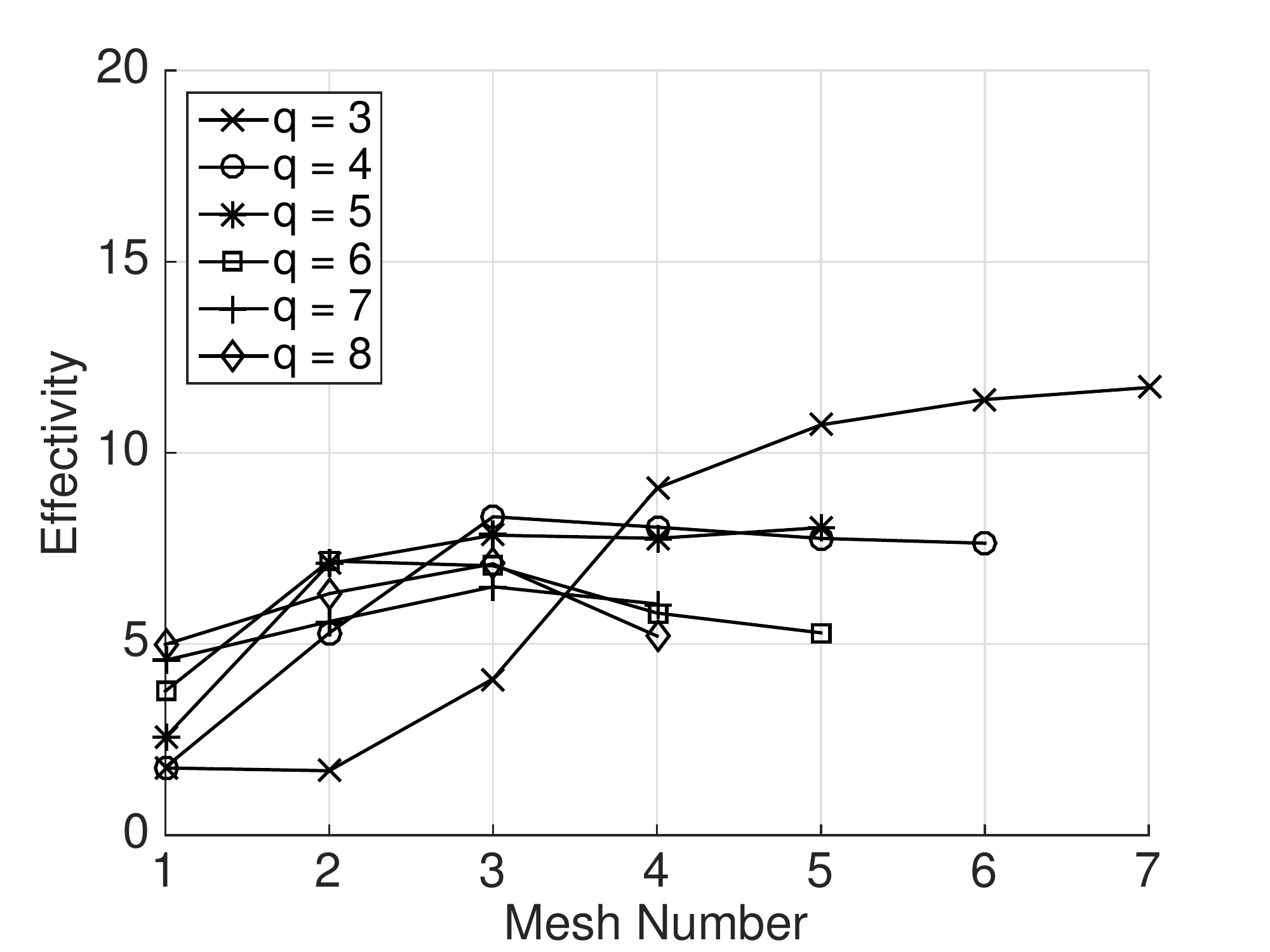}}
    \subfloat[$k=30$, $\mathcal{E}_R$]{\label{fig:eff:robin:30}\includegraphics[width=0.33\textwidth]{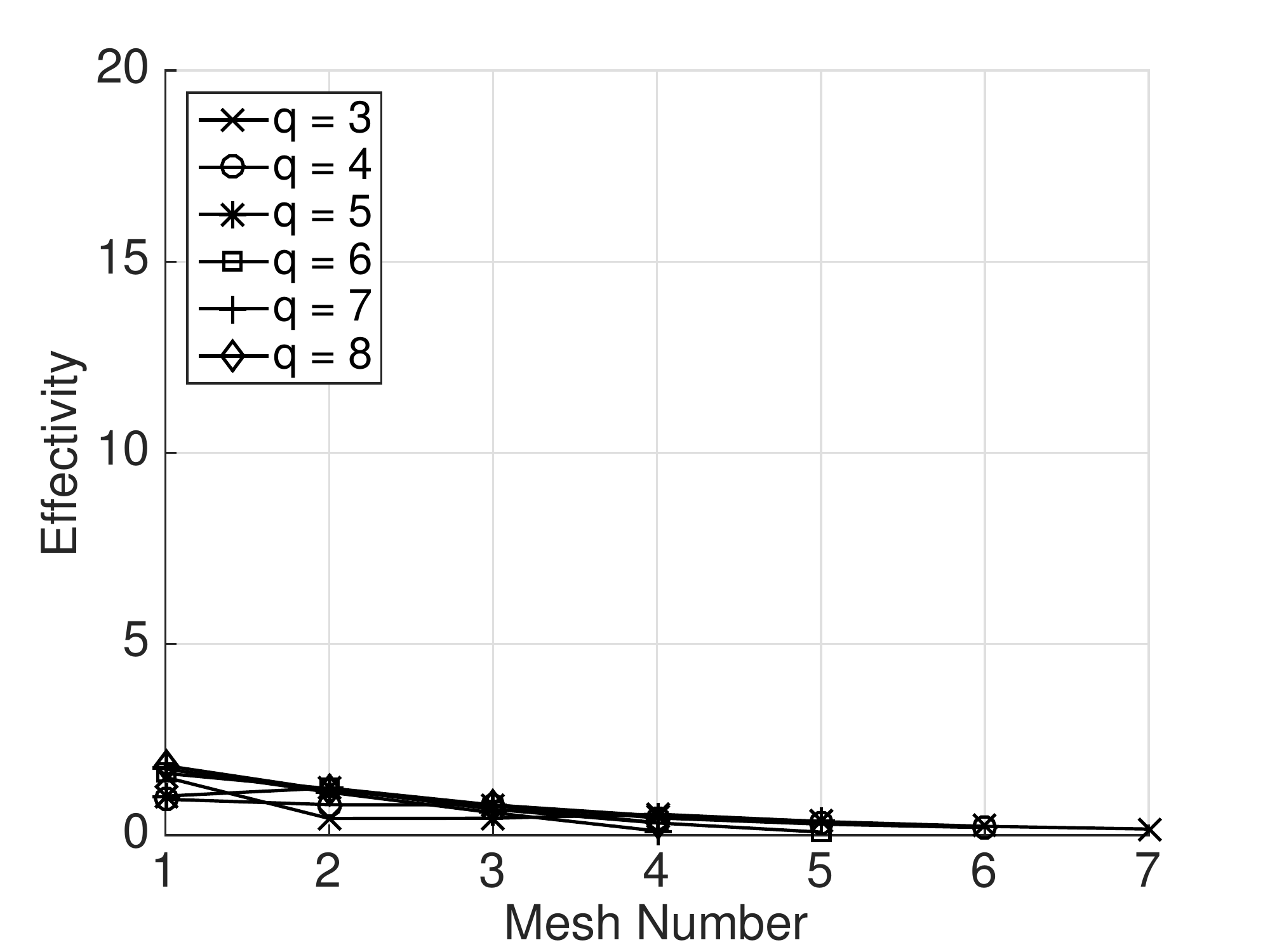}} \\
    \subfloat[$k=40$, $\mathcal{E}_{\jmp{u_{hp}}}$]{\label{fig:eff:u:40}\includegraphics[width=0.33\textwidth]{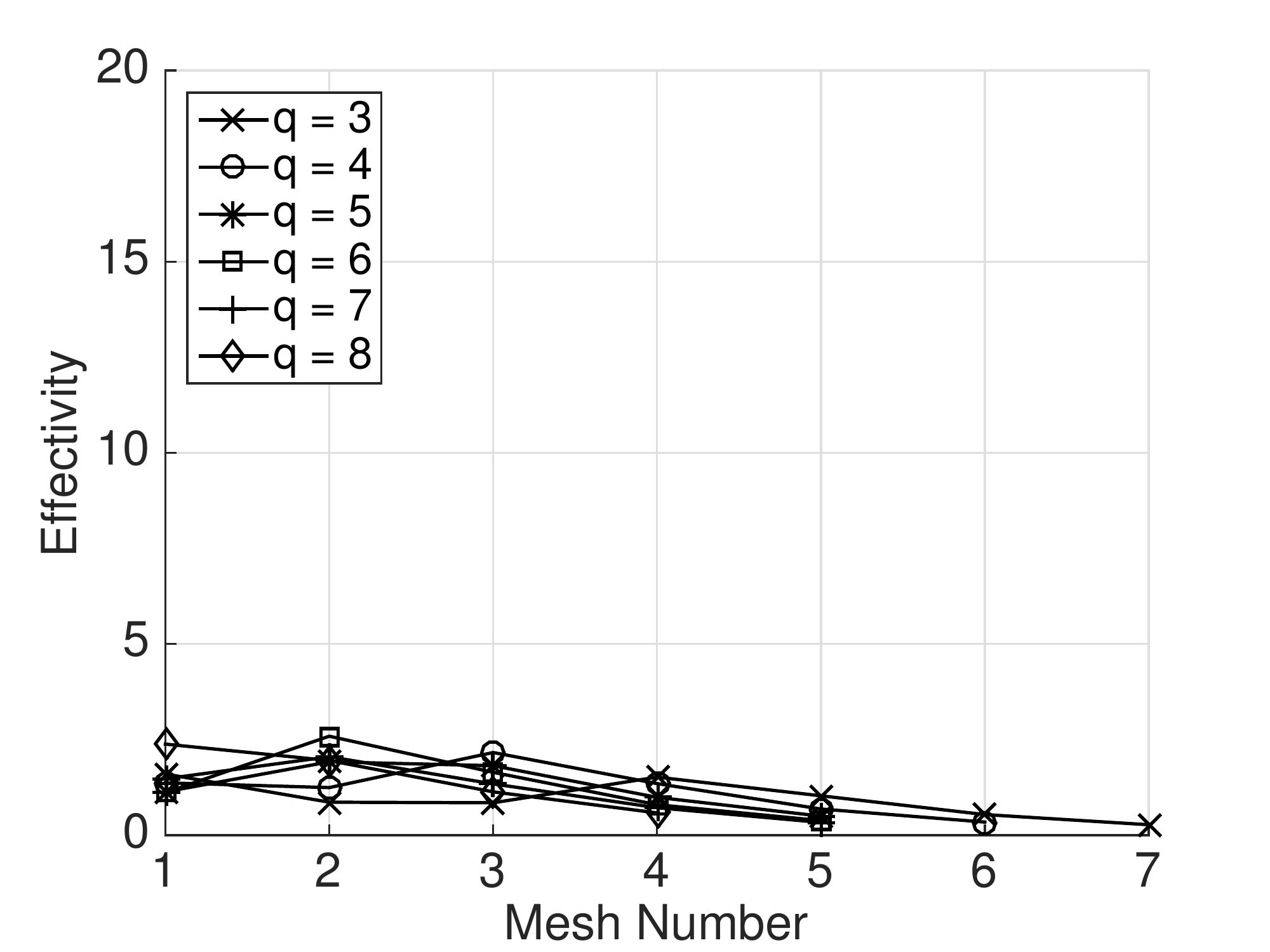}}
    \subfloat[$k=40$, $\mathcal{E}_{\jmp{\nabla u_{hp}}}$]{\label{fig:eff:gradu:40}\includegraphics[width=0.33\textwidth]{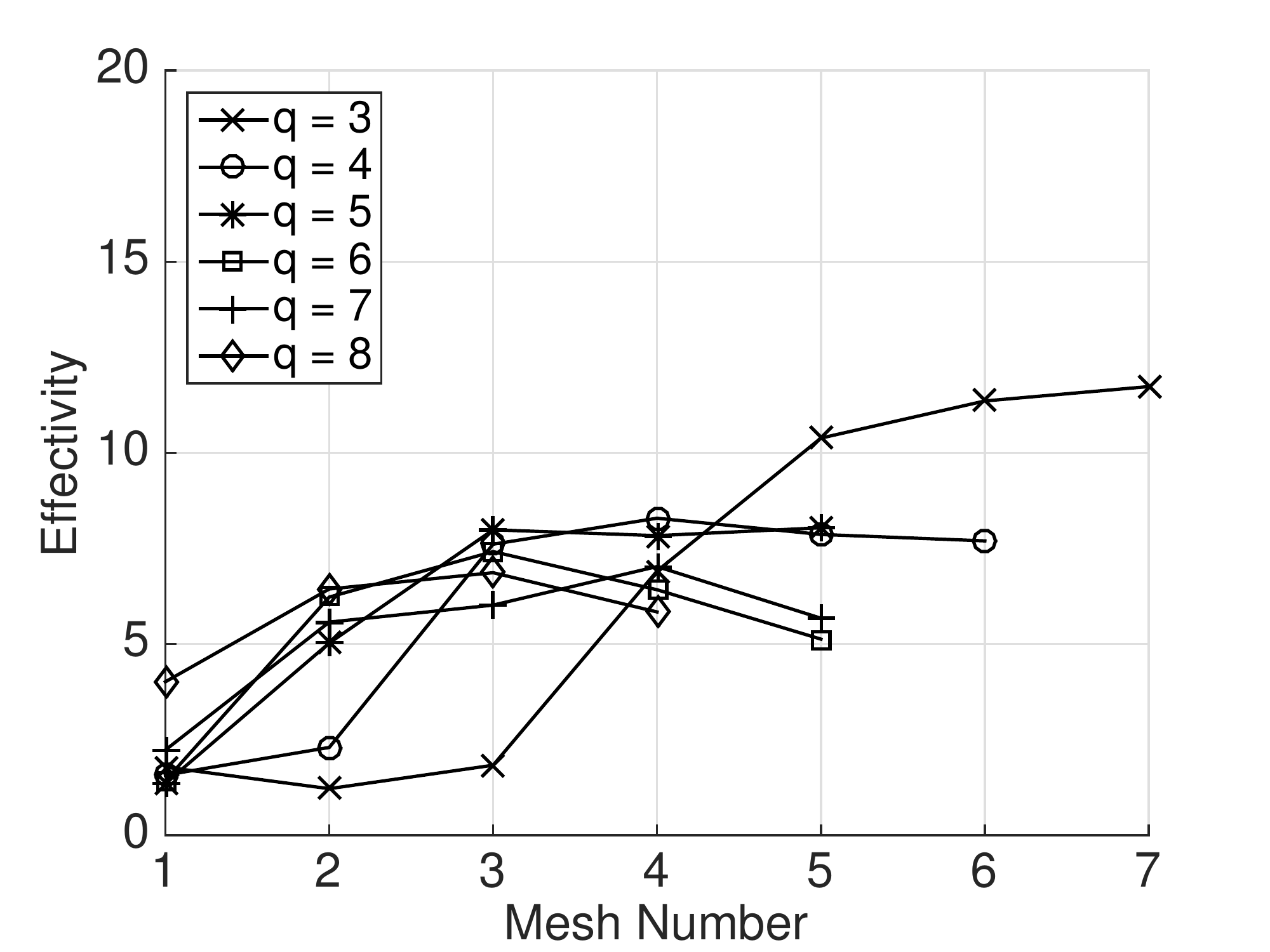}}
    \subfloat[$k=40$, $\mathcal{E}_R$]{\label{fig:eff:robin:40}\includegraphics[width=0.33\textwidth]{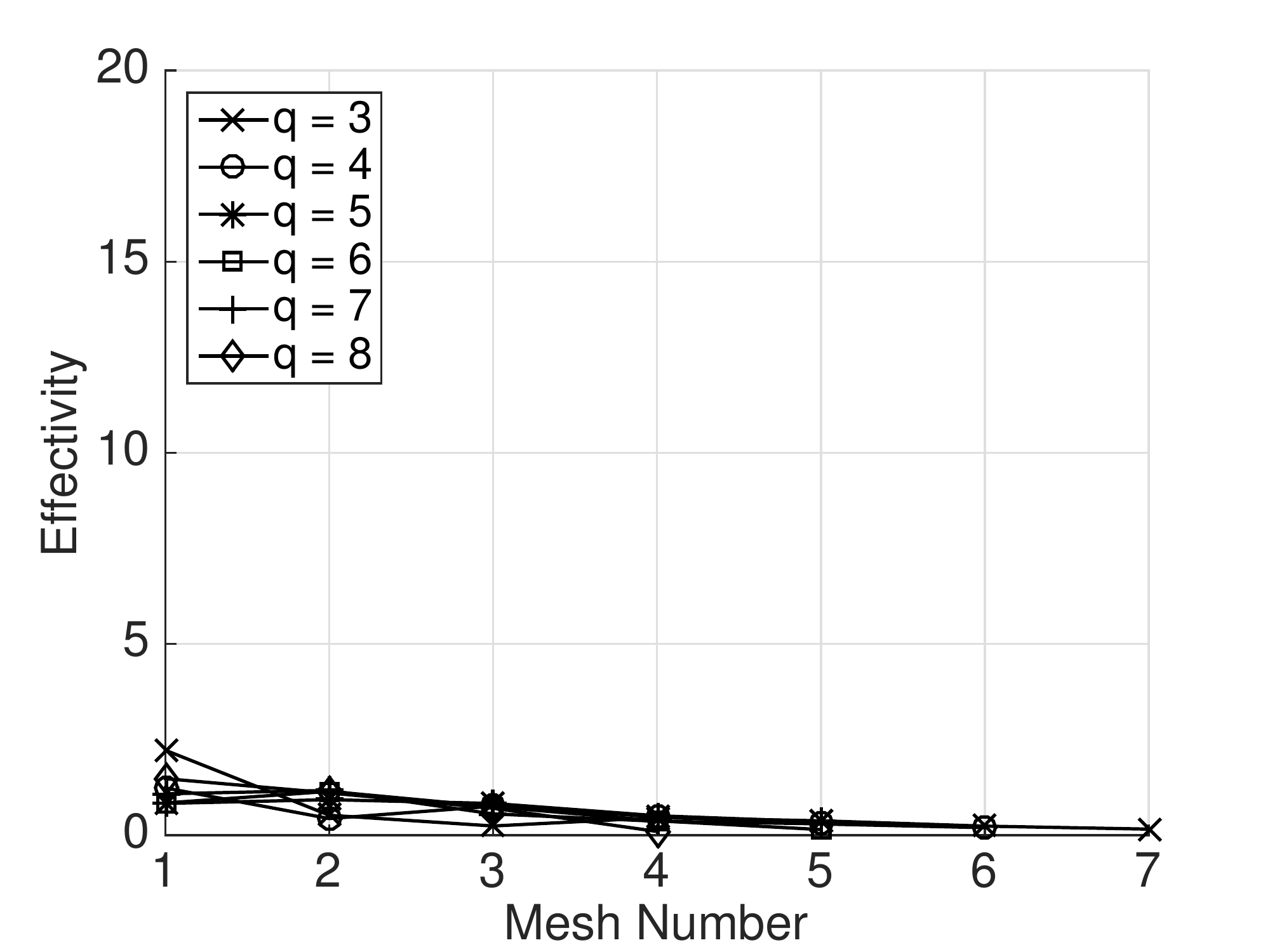}} \\
    \subfloat[$k=50$, $\mathcal{E}_{\jmp{u_{hp}}}$]{\label{fig:eff:u:50}\includegraphics[width=0.33\textwidth]{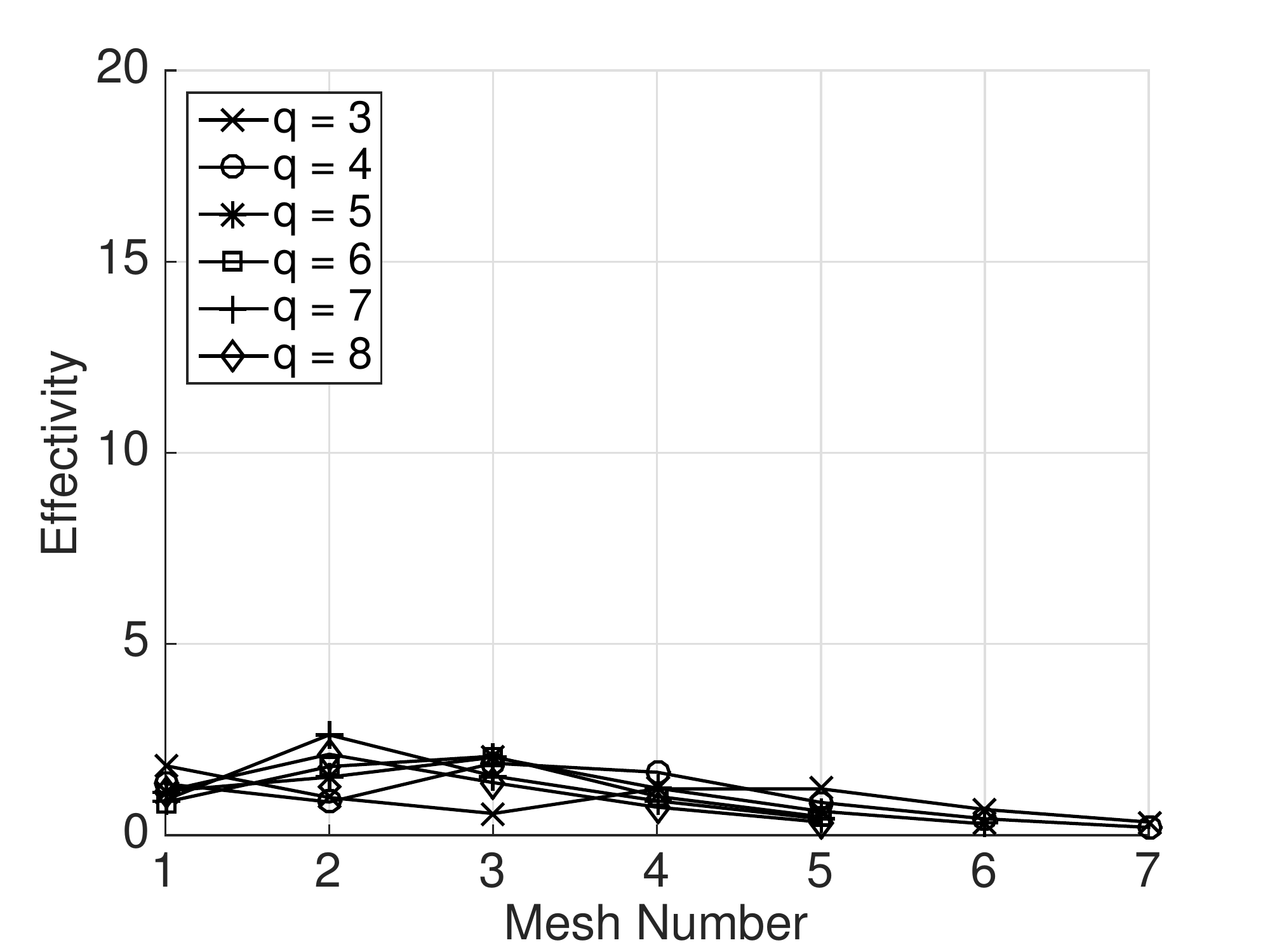}}
    \subfloat[$k=50$, $\mathcal{E}_{\jmp{\nabla u_{hp}}}$]{\label{fig:eff:gradu:50}\includegraphics[width=0.33\textwidth]{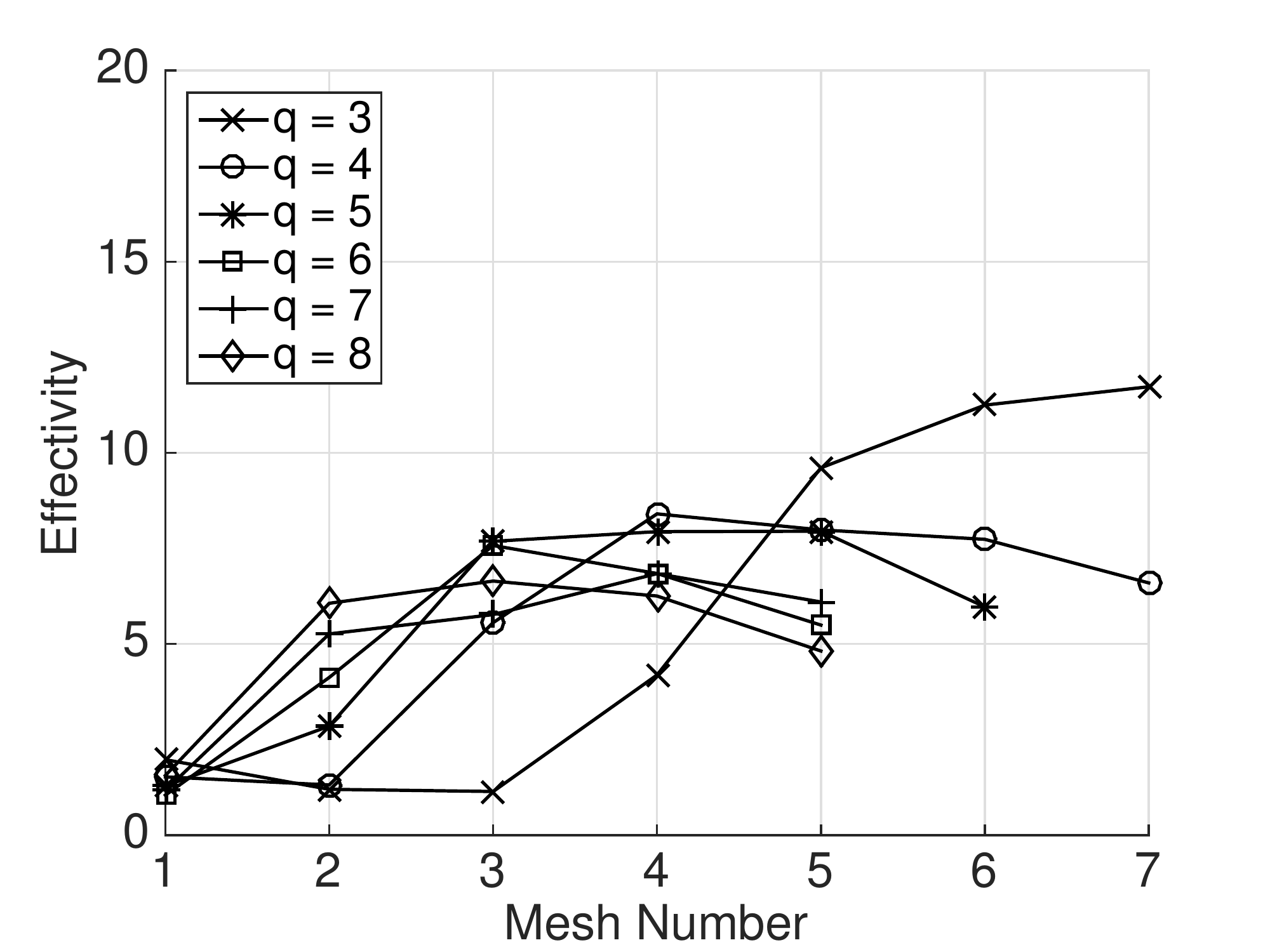}}
    \subfloat[$k=50$, $\mathcal{E}_R$]{\label{fig:eff:robin:50}\includegraphics[width=0.33\textwidth]{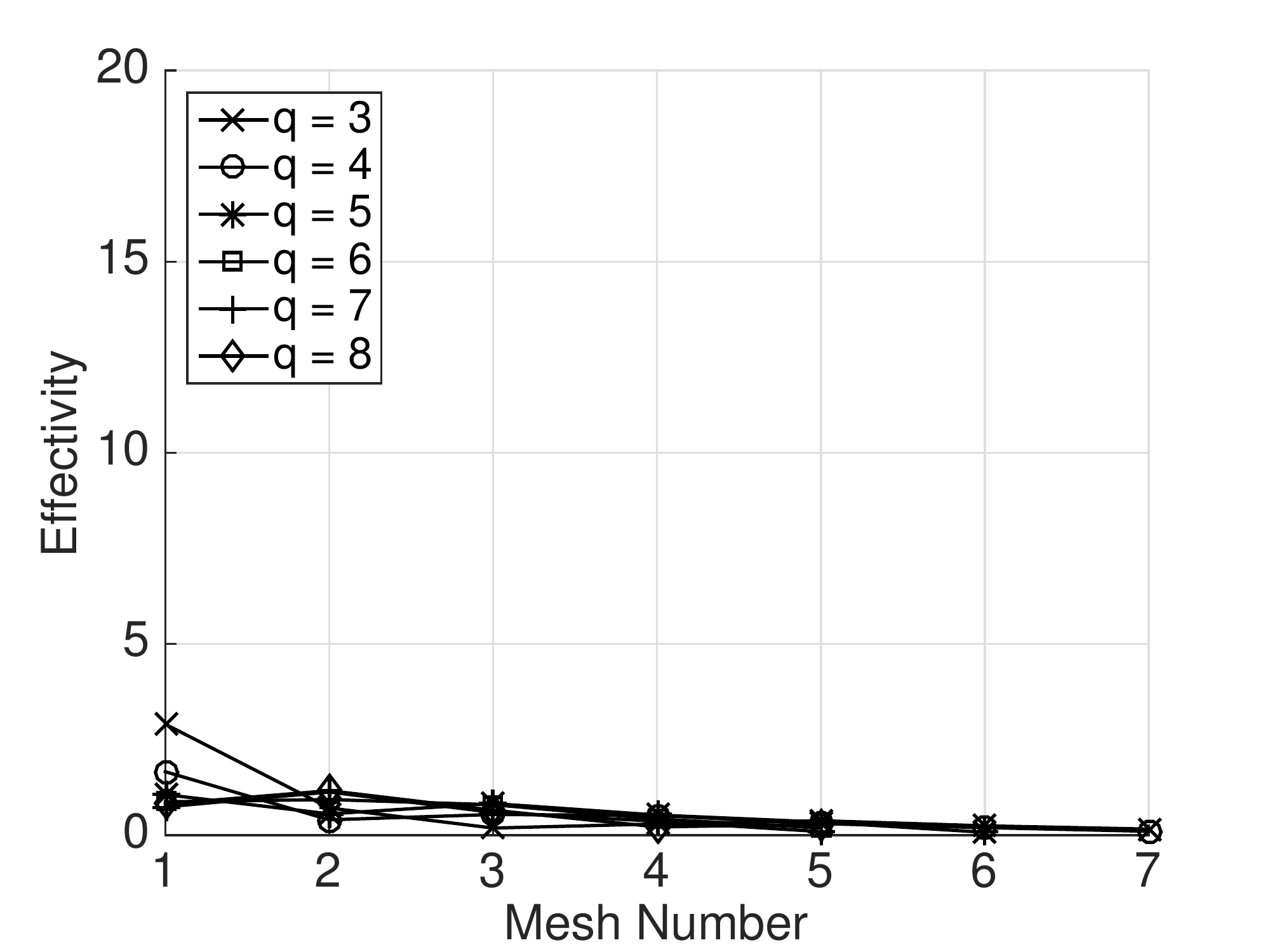}} \\
    \caption{Effectivities of individual components of the error indicators for $h$--refinement with fixed effective polynomial degree of the smooth analytical Hankel solution with different wavenumbers.}
    \label{fig:eff:parts}
\end{figure}
\subsection{Efficiency of the \emph{a posteriori} error indicator}\label{section:aposteriori_effectivity}

The selection of the exponents of $h_K$ and $q_K$ in the weights present in~\eqref{eqn:error_indicator},
together with the independence on the wavenumber $k$, have been determined by numerical experimentation. 
To this end, we considered the example presented in the previous section, cf.~\eqref{eqn:hankal_anal}, whereby the
numerical approximation is computed
on a series of uniform computational meshes, with uniform effective polynomial degrees $q$, 
for a range of wave numbers $k$. In each case, we computed the effectivity index of each constituent
term arising in ${\mathfrak E}(u_h,h,\vect{p})$, whereby the dependency on $h_K$, $q_K$, and $k$ 
was eliminated; note that with the removal of $h_K$, $q_K$, and $k$, the effectivity index
of each term is computed by dividing by $\norm{u-u_{hp}}_{L^2(\Omega)}$. More precisely,
effectivity indices were computed for $q=3,\dots,8$ and $k=20,30,40,50$, based on starting from 
a uniform $4\times 4$ mesh consisting of square elements. On the basis of these results,
the dependence of each term on $h_K$, $q_K$, and $k$ was established. The final
effectivity indices for the correctly scaled empirical {\em a posteriori} error indicator 
${\mathfrak E}(u_h,h,\vect{p})$, i.e., ${\mathfrak E}(u_h,h,\vect{p})/\norm{u-u_{hp}}_{L^2(\Omega)}$ 
are presented in Figure~\ref{fig:eff}. Here, we observe that that the effectivity indices have 
roughly the same values for all the selected values of $h$, $q$, and $k$; however,
at higher wave numbers, pre-asymptotic behaviour leads to some increase in the effectivity indices
as the mesh is refined, due to the fact that the mesh size is too large for the wavelength. 
We note that this behaviour is more noticeable in the case when $q=3$.

Finally, we compute the effectivity index for each individual term arising in the
definition of the error indicator ${\mathfrak E}(u_h,h,\vect{p})$, cf.~\eqref{eqn:error_indicator}. 
More precisely, we define
\begin{align*}
\mathcal{E}_{\jmp{u_{hp}}} &\coloneqq \frac{\left(\sum_{K\in\mesh} \norm*{\alpha^{\nicefrac{1}{2}} h_K^{\nicefrac12}q_F^{-\nicefrac12} \jmp{u_{hp}}}_{L^2(\partial K \setminus \partial\Omega)}^2\right)^{\nicefrac12}}{\norm*{u-u_{hp}}_{L^2(\Omega)}}, \\
\mathcal{E}_{\jmp{\nabla u_{hp}}} &\coloneqq \frac{\left(\sum_{K\in\mesh} \norm*{\beta^{\nicefrac12} h_K^{\nicefrac32}q_K^{-\nicefrac32} \jmp{\nabla u_{hp}}}_{L^2(\partial K \setminus \partial\Omega)}^2\right)^{\nicefrac12}}{\norm*{u-u_{hp}}_{L^2(\Omega)}}, \\
\mathcal{E}_{R} &\coloneqq \frac{\left(\sum_{K\in\mesh} \norm*{\delta^{\nicefrac{1}{2}} h_K^{\nicefrac32}q_K^{-\nicefrac32} \left( g_R - \nabla u_{hp}\cdot\vect{n}_F + ik u_{hp}\right)}^2_{L^2(\partial K \cap \Gamma_R)}\right)^{\nicefrac12}}{\norm*{u-u_{hp}}_{L^2(\Omega)}};
\end{align*}
the results for the case when $k=20,30,40,50$ are depicted in Figure~\ref{fig:eff:parts}. 
Here, we observe that each individual effectivity index is roughly constant 
for all the selected values of $h$, $q$, and $k$, except within the pre-asymptotic region. 
For this smooth problem, we clearly observe that the dominant part of the error indicator 
involves the jump in the gradient of the numerical solution.

\begin{remark}
We note there that we have only computed the weightings for the interior and Robin faces. In the case of Dirichlet boundary conditions we assume that the weighting scales the same as the term involving $\jmp{u_{hp}}$.
\end{remark}

\subsection{$hp$--Adaptive refinement}
In this section we consider computational performance of the proposed $hp$--adaptive refinement algorithm,
with directional adaptivity, 
for a range of test problems in both two- and three-dimensions. 
To this end, employ the fixed fraction refinement strategy to mark elements for refinement; 
throughout this section, we set the refinement fraction equal to 25\% 
of the elements with the largest contribution to the error bound.
Furthermore, we allow the meshes $\mesh$ to be `1-irregular', i.e., each face of any 
element $K\in\mesh$ contains at most one hanging node (which, for simplicity, we assume
to be at the barycenter of the corresponding face) and each edge of each face contains at most one
hanging node (yet again assumed to be at the barycenter of the edge). We also only allow the 
effective polynomial degree $q_K$ to vary by one between neighbouring elements.

For each test problem, we compare the performance of employing $hp$--adaptive refinement with 
$h$--adaptivity. In the latter case, we consider a standard $h$--adaptive algorithm, i.e., adaptive
mesh refinement without directional adaptivity, as well as an $h$--adaptive strategy which incorporates
directional adaptivity; here, we shall consider the two cases when directional adaptivity is either 
undertaken only on the elements marked for refinement, as well as the case when it is performed on all
elements in the computational mesh. In the $hp$--setting, similar comparisons will be made, in addition
to studying the case when directional adaptivity is only performed on elements marked for $p$--refinement.

We note that when $hp$--refinement is exploited we often reach a point where the 
$L^2$--norm of the error and {\em a posteriori} error bound stagnates, in the sense that
both quantities no longer tend to zero, and indeed may start to oscillate, as further refinement
is undertaken. This is caused by the fact that as the relative magnitude of $q_K$,
with respect to $h_Kk$, becomes large,
the local plane wave bases are very ill-conditioned. 
In this situation, we simply stop the numerical experiments and discard further results; however, 
possible improvements based on ensuring $q_K$ is well behaved with respect to $h_Kk$ could be 
implemented; cf. \cite{Cessenat1998, Huttunen2002, Luostari2013} for details.

\subsubsection{Example 1 --- Smooth solution (Hankel function)} \label{sec-hprefine-hankel}
\begin{figure}[pt]
    \configfigure
    \subfloat[$k=20$; $h$--refinement]{\label{fig:hankel:error:20h}\includegraphics[width=0.4\textwidth]{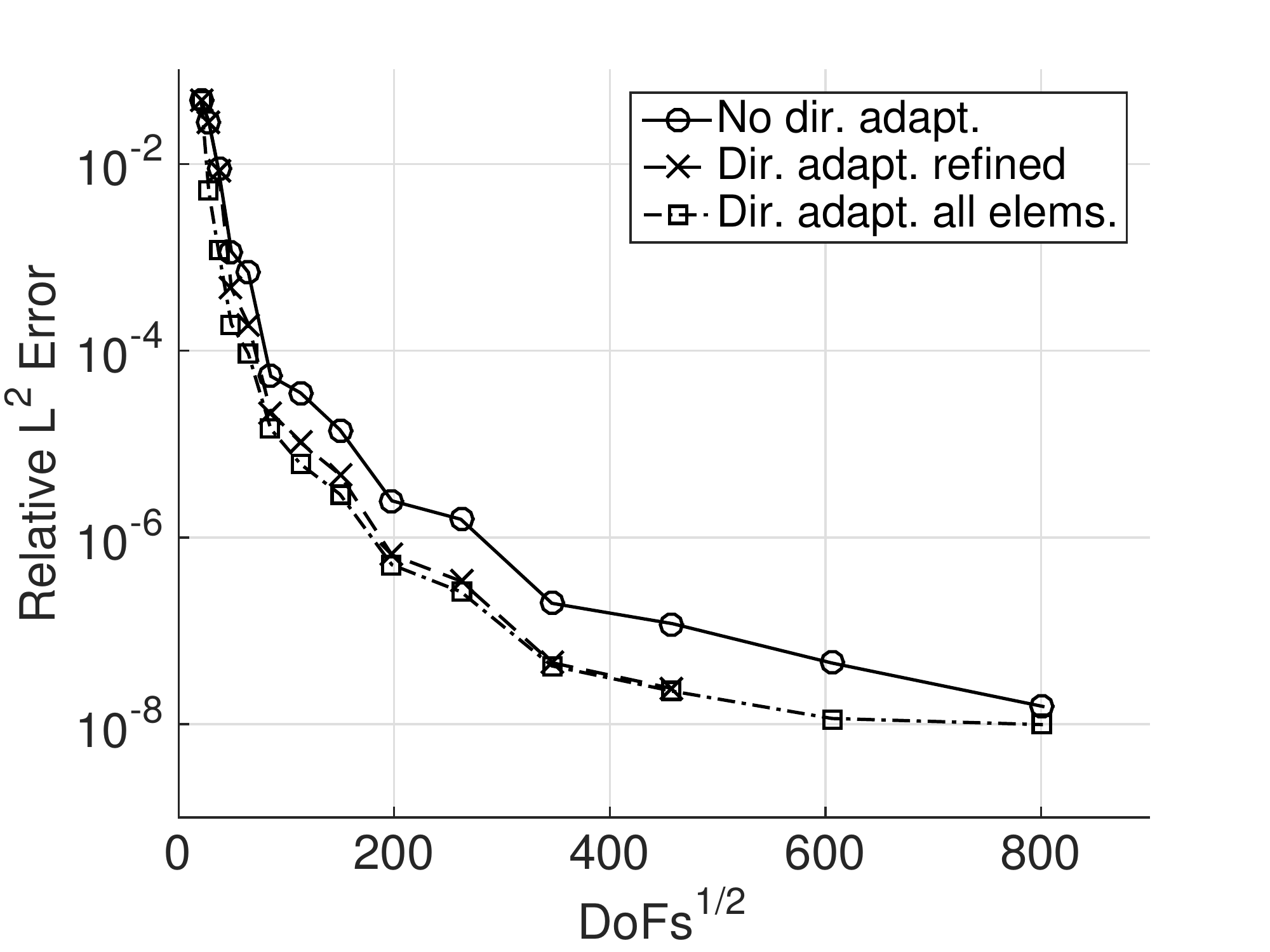}}
    \subfloat[$k=20$; $h$--refinement]{\label{fig:hankel:eff:20h}\includegraphics[width=0.4\textwidth]{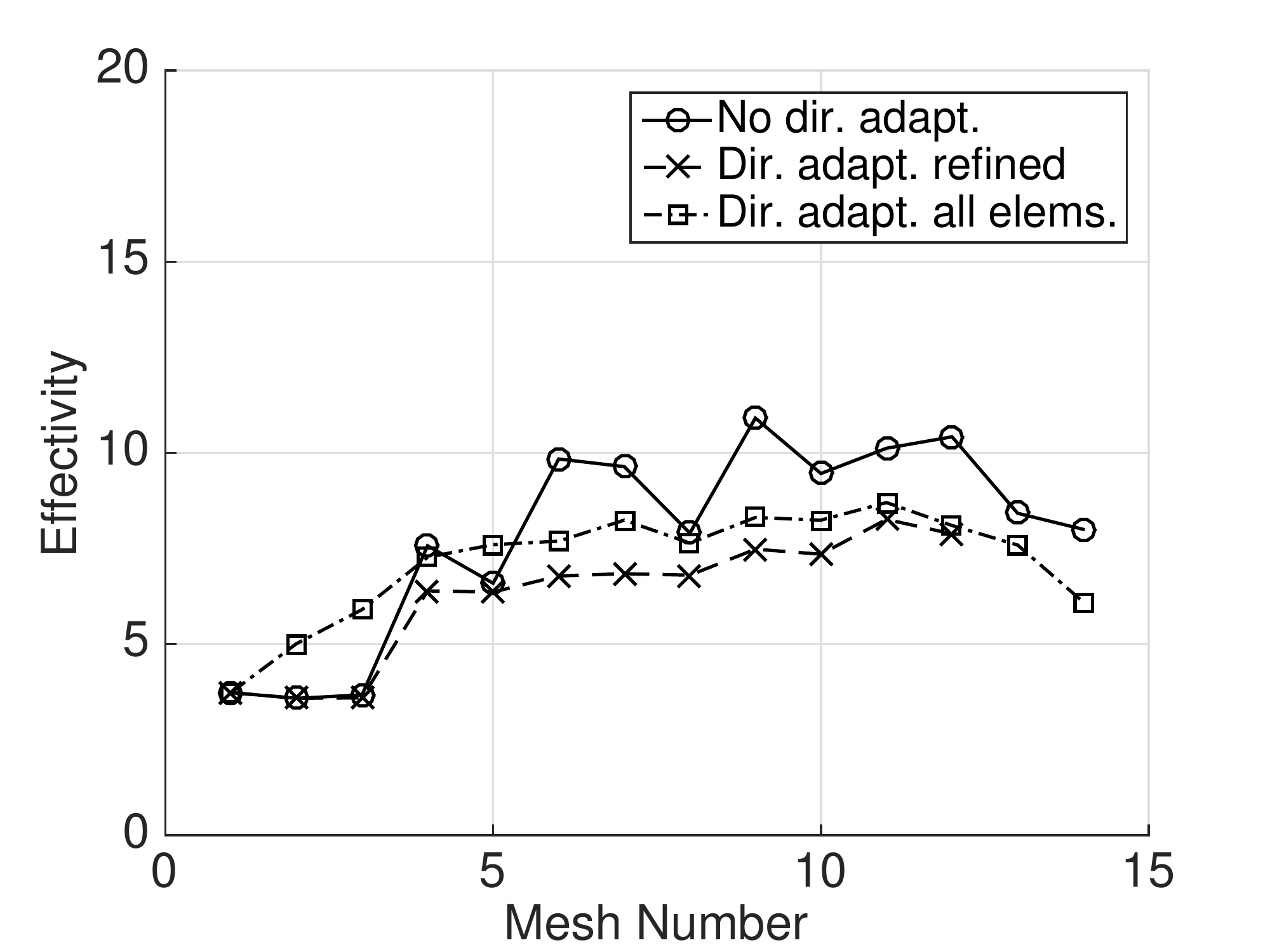}} \\
    \subfloat[$k=20$; $hp$--refinement]{\label{fig:hankel:error:20hp}\includegraphics[width=0.4\textwidth]{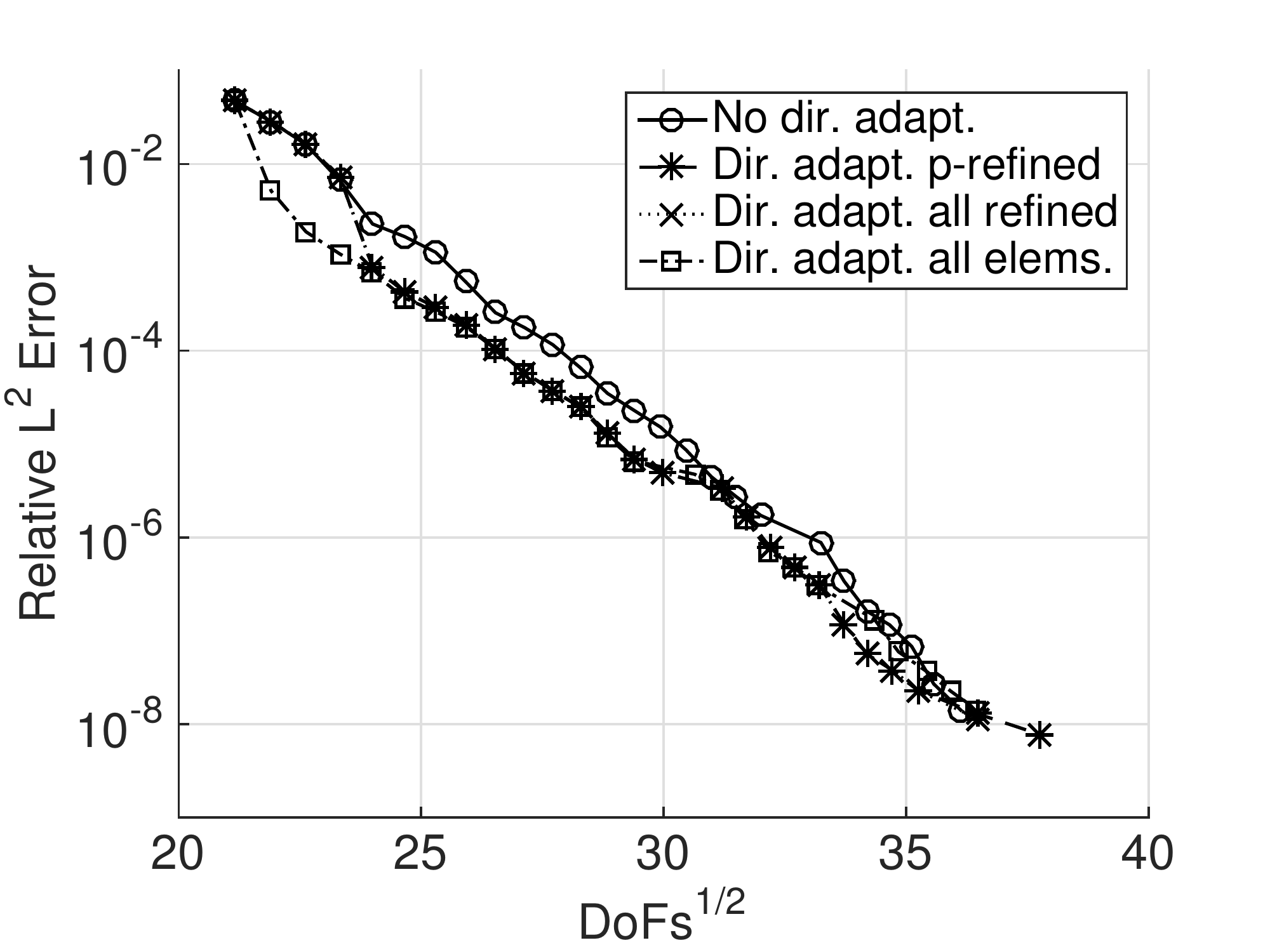}}
    \subfloat[$k=20$; $hp$--refinement]{\label{fig:hankel:eff:20hp}\includegraphics[width=0.4\textwidth]{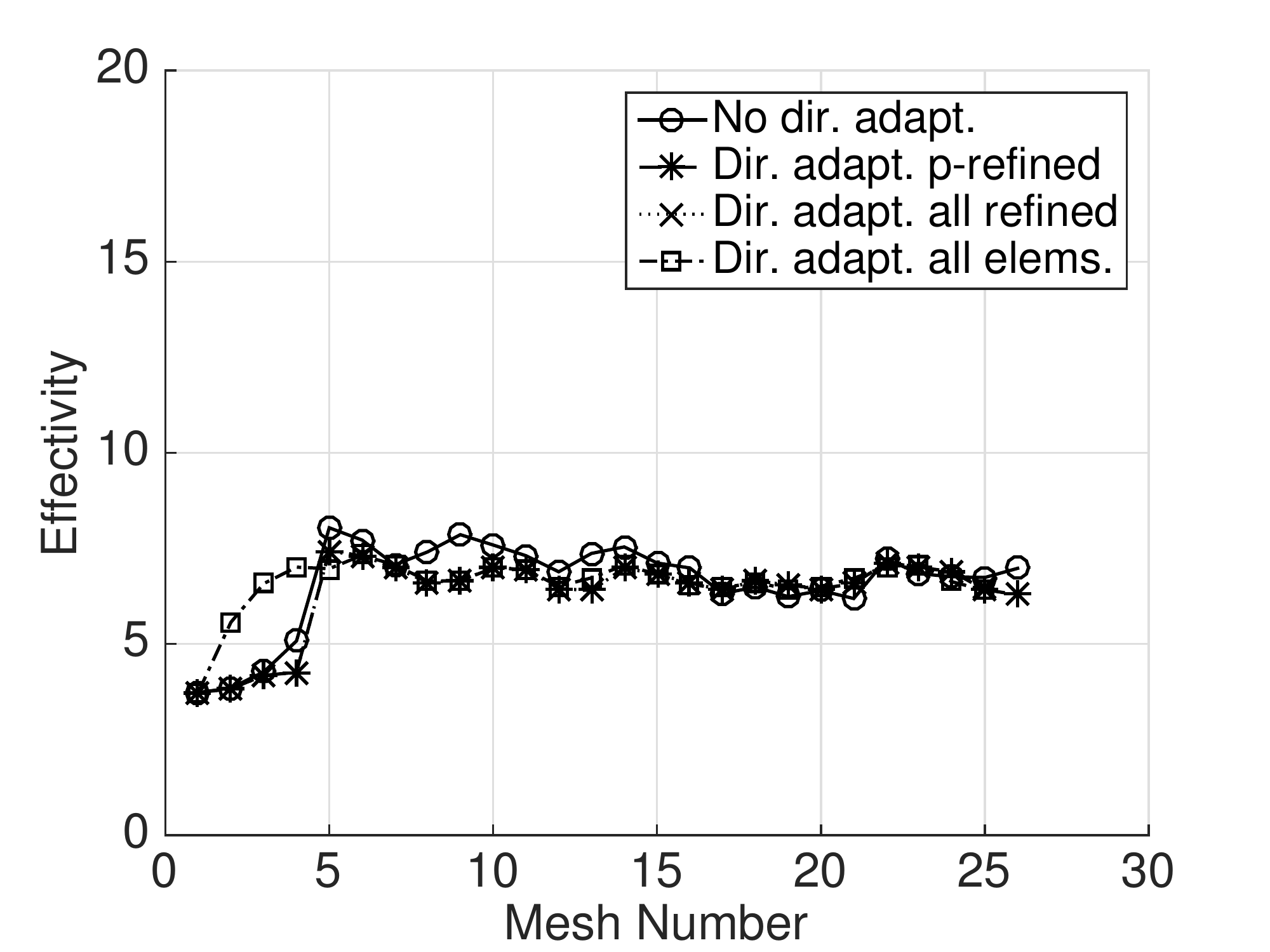}} \\
    \subfloat[$k=50$; $h$--refinement]{\label{fig:hankel:error:50h}\includegraphics[width=0.4\textwidth]{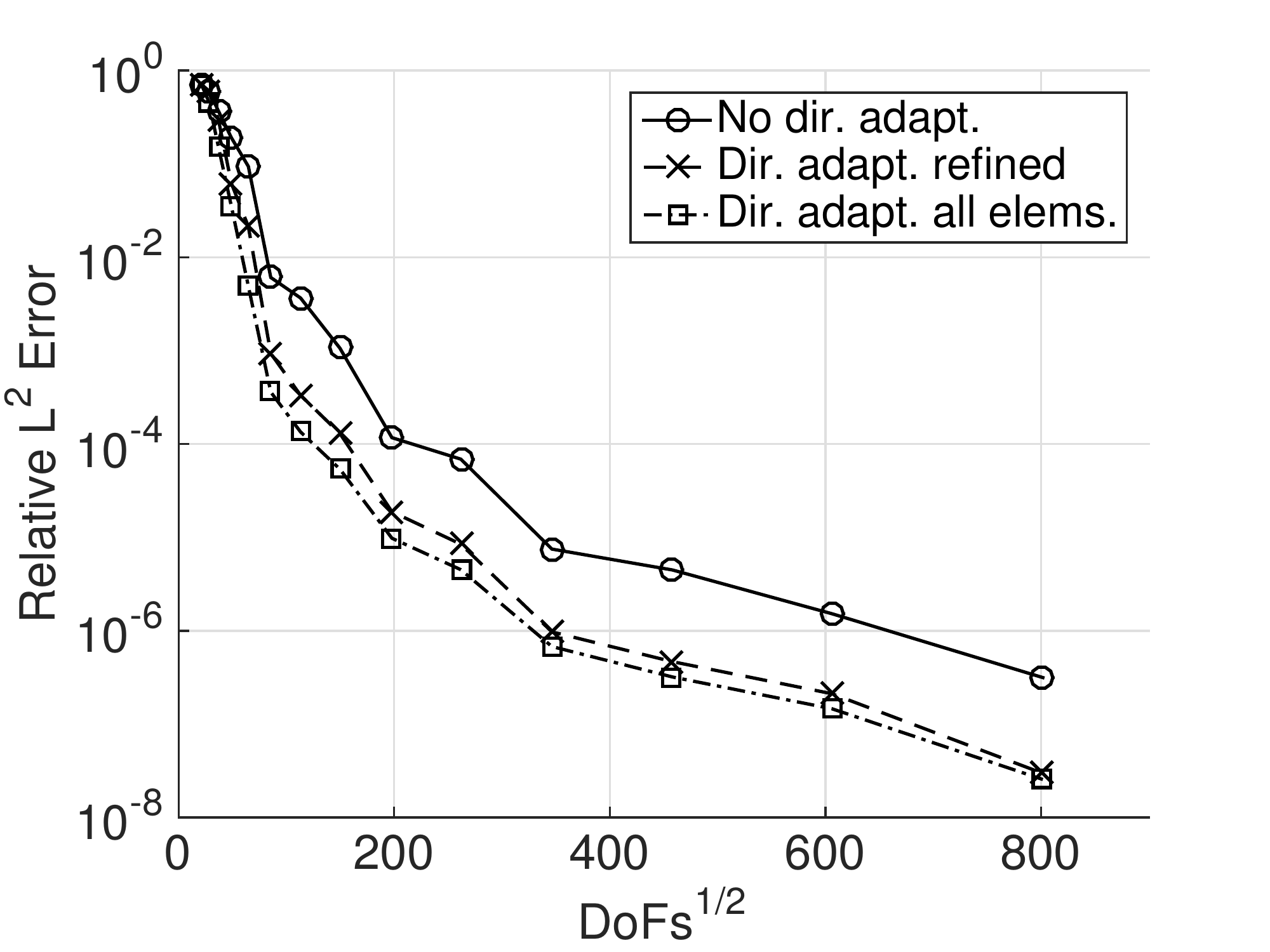}}
    \subfloat[$k=50$; $h$--refinement]{\label{fig:hankel:eff:50h}\includegraphics[width=0.4\textwidth]{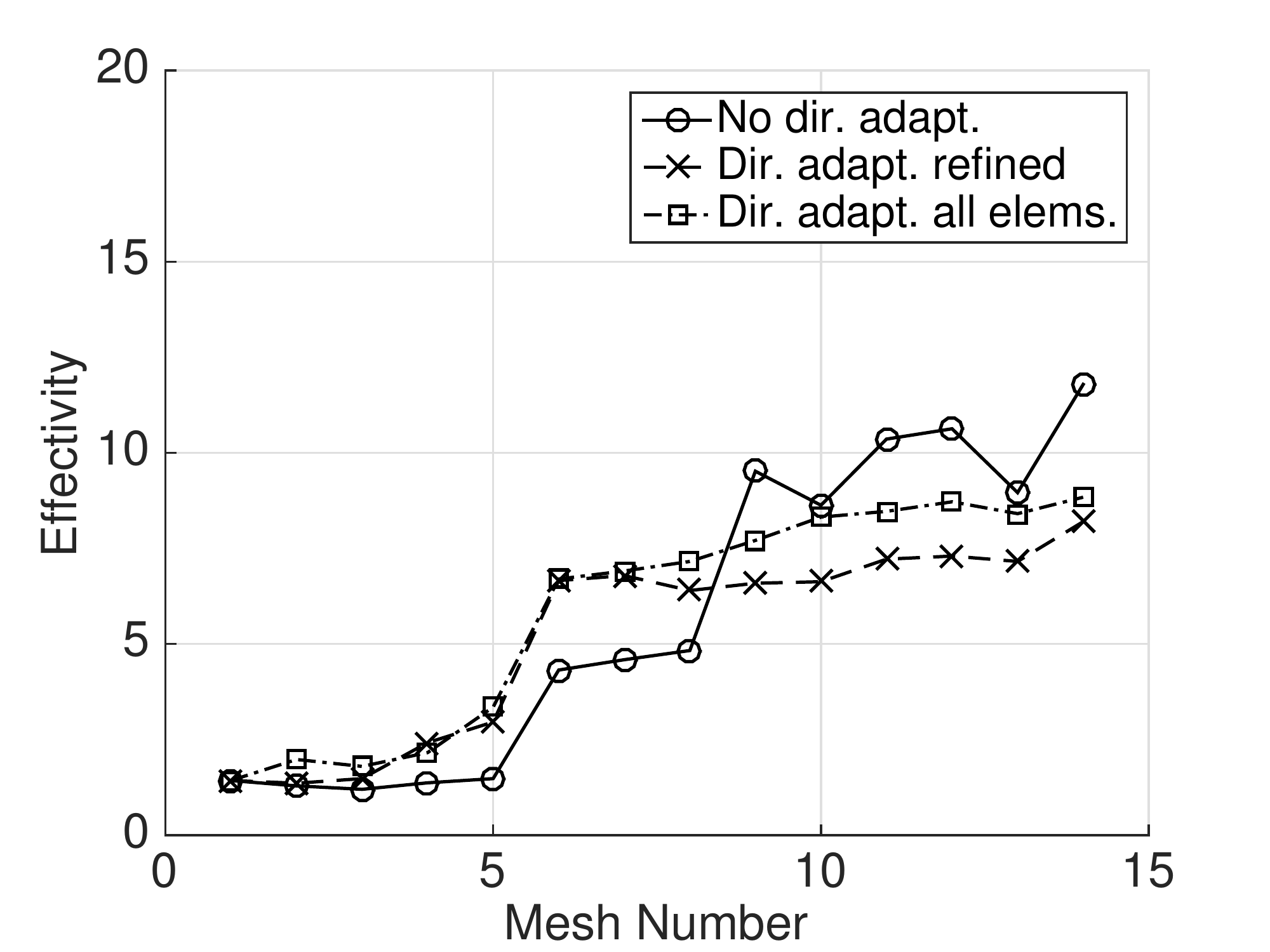}} \\
    \subfloat[$k=50$; $hp$--refinement]{\label{fig:hankel:error:50hp}\includegraphics[width=0.4\textwidth]{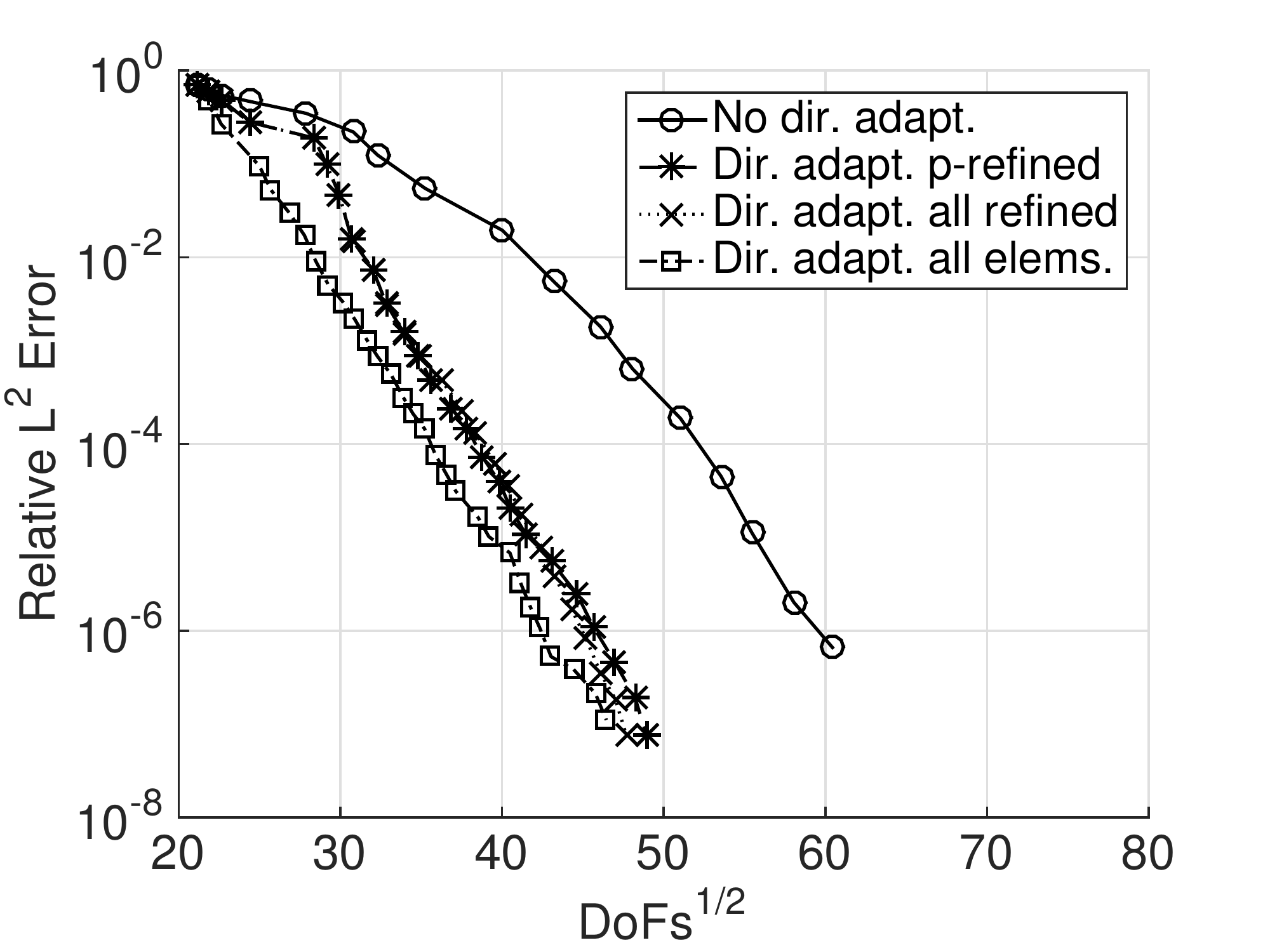}}
    \subfloat[$k=50$; $hp$--refinement]{\label{fig:hankel:eff:50hp}\includegraphics[width=0.4\textwidth]{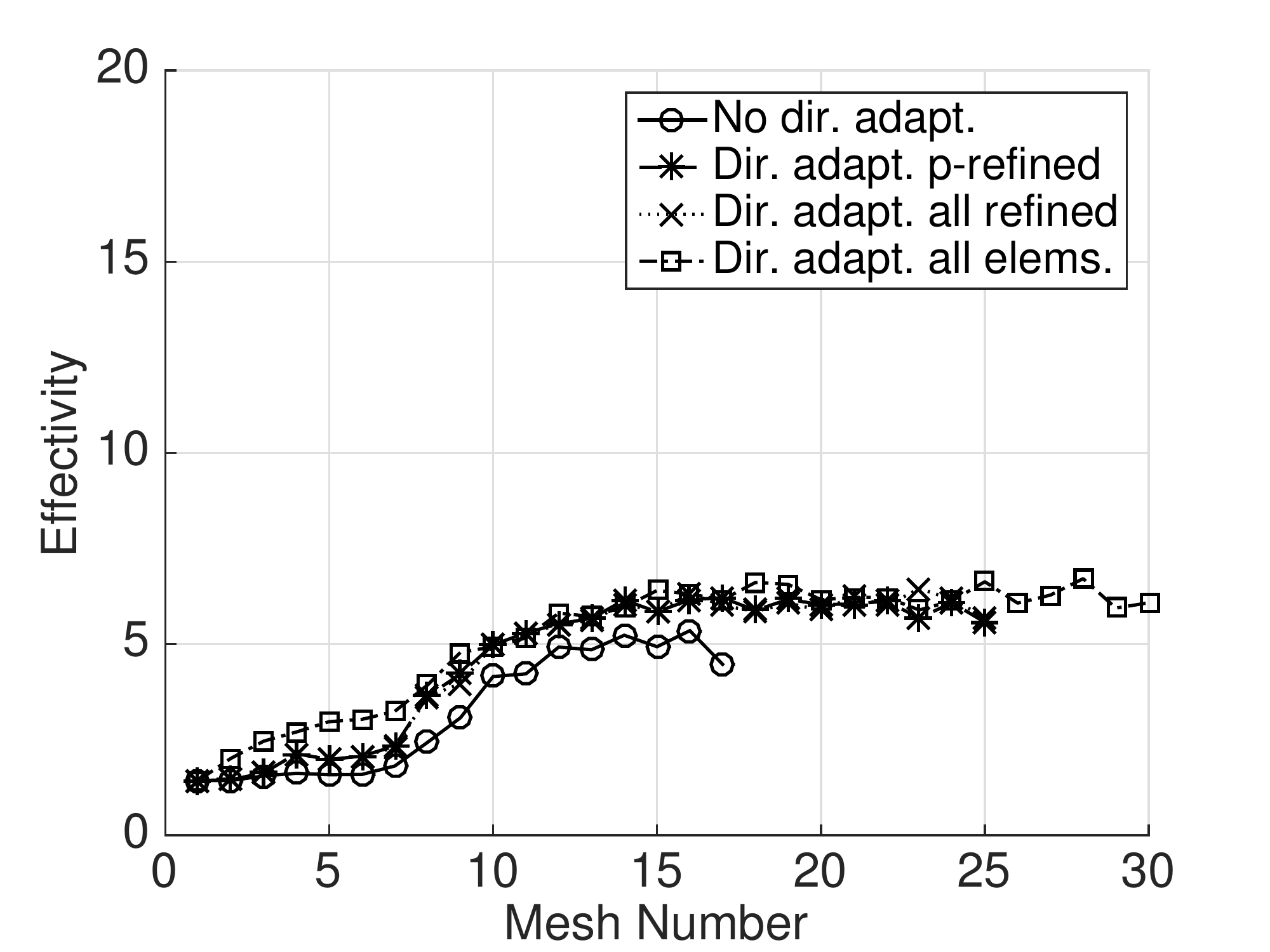}}
    \caption{Example 1: \protect\subref{fig:hankel:error:20h} $L^2$-error and \protect\subref{fig:hankel:eff:20h} Effectivity index for $h$--refinement with wavenumber $k=20$; \protect\subref{fig:hankel:error:20hp} $L^2$-error and \protect\subref{fig:hankel:eff:20hp} Effectivity index for $hp$--refinement with $k=20$; \protect\subref{fig:hankel:error:50h} $L^2$-error and \protect\subref{fig:hankel:eff:50h} Effectivity index for $h$--refinement with $k=50$; \protect\subref{fig:hankel:error:50hp} $L^2$-error and \protect\subref{fig:hankel:eff:50hp} Effectivity index for $hp$--refinement with $k=50$.}
\end{figure}
\begin{figure}[pt]
    \configfigure
    \subfloat[$k=20$]{\label{fig:hankel:errorcompare:20}\includegraphics[width=0.4\textwidth]{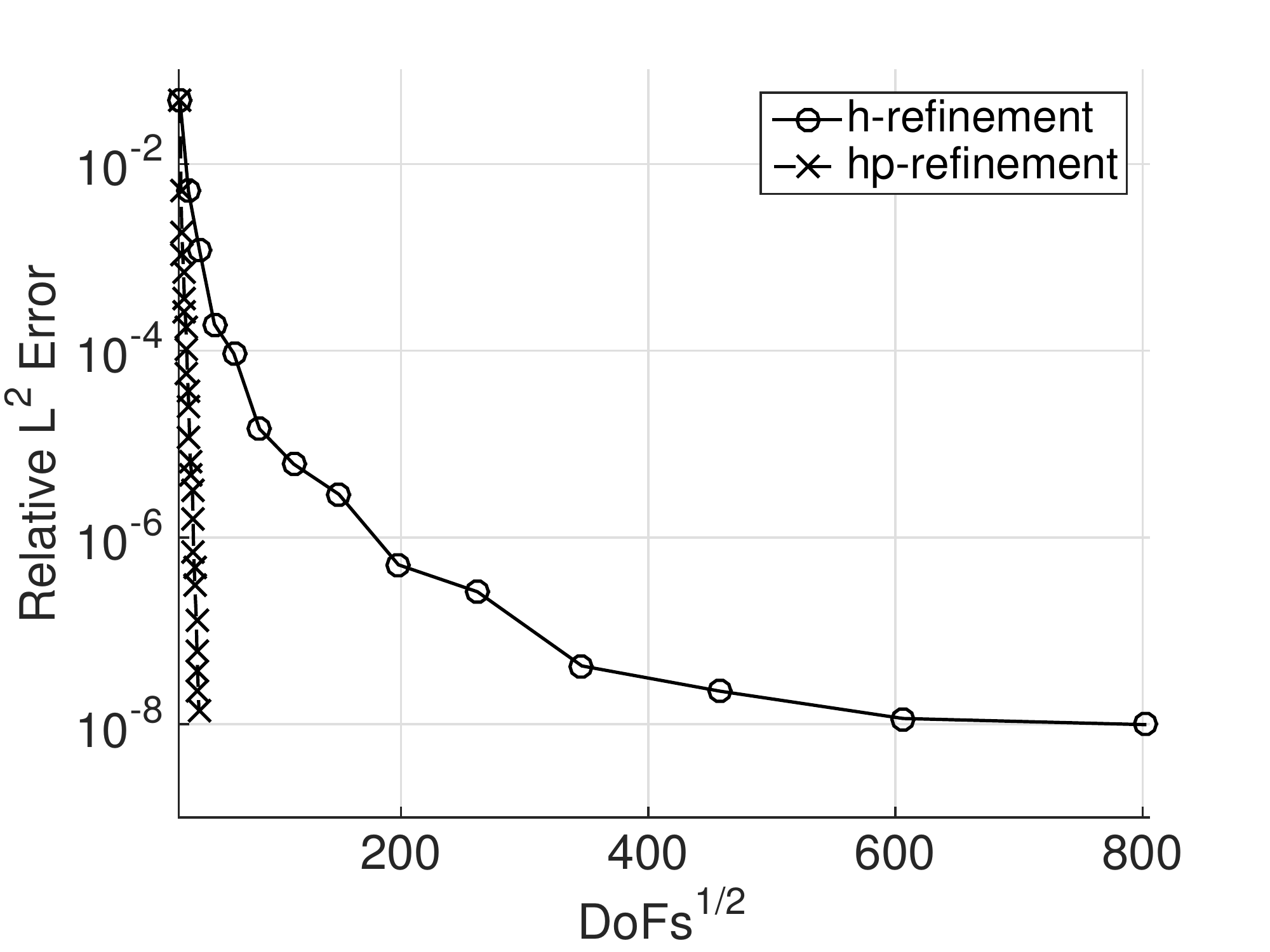}}
    \subfloat[$k=50$]{\label{fig:hankel:errorcompare:50}\includegraphics[width=0.4\textwidth]{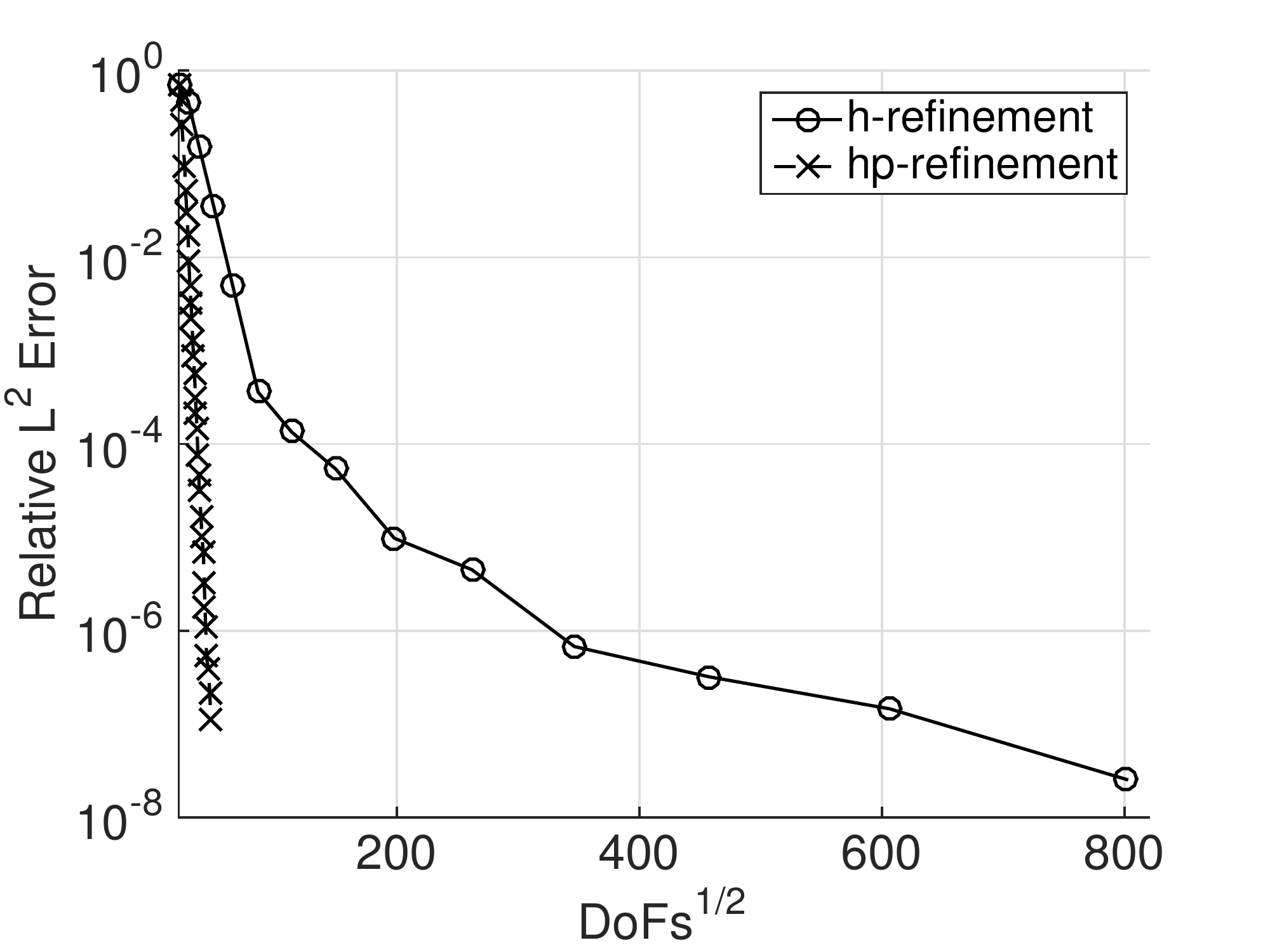}}
    \caption{Example 1: Comparison of relative $L^2$-error for $h$-- and $hp$--refinement, with direction adaptivity on all elements, for wavenumbers \protect\subref{fig:hankel:errorcompare:20} $k=20$ and \protect\subref{fig:hankel:errorcompare:50} $k=50$.}
    \label{fig:hankel:errorcompare}
\end{figure}
\begin{figure}[pt]
    \configfigure
    \subfloat[$k=20$]{\label{fig:hankel:mesh:h:20}\includegraphics[width=0.4\textwidth]{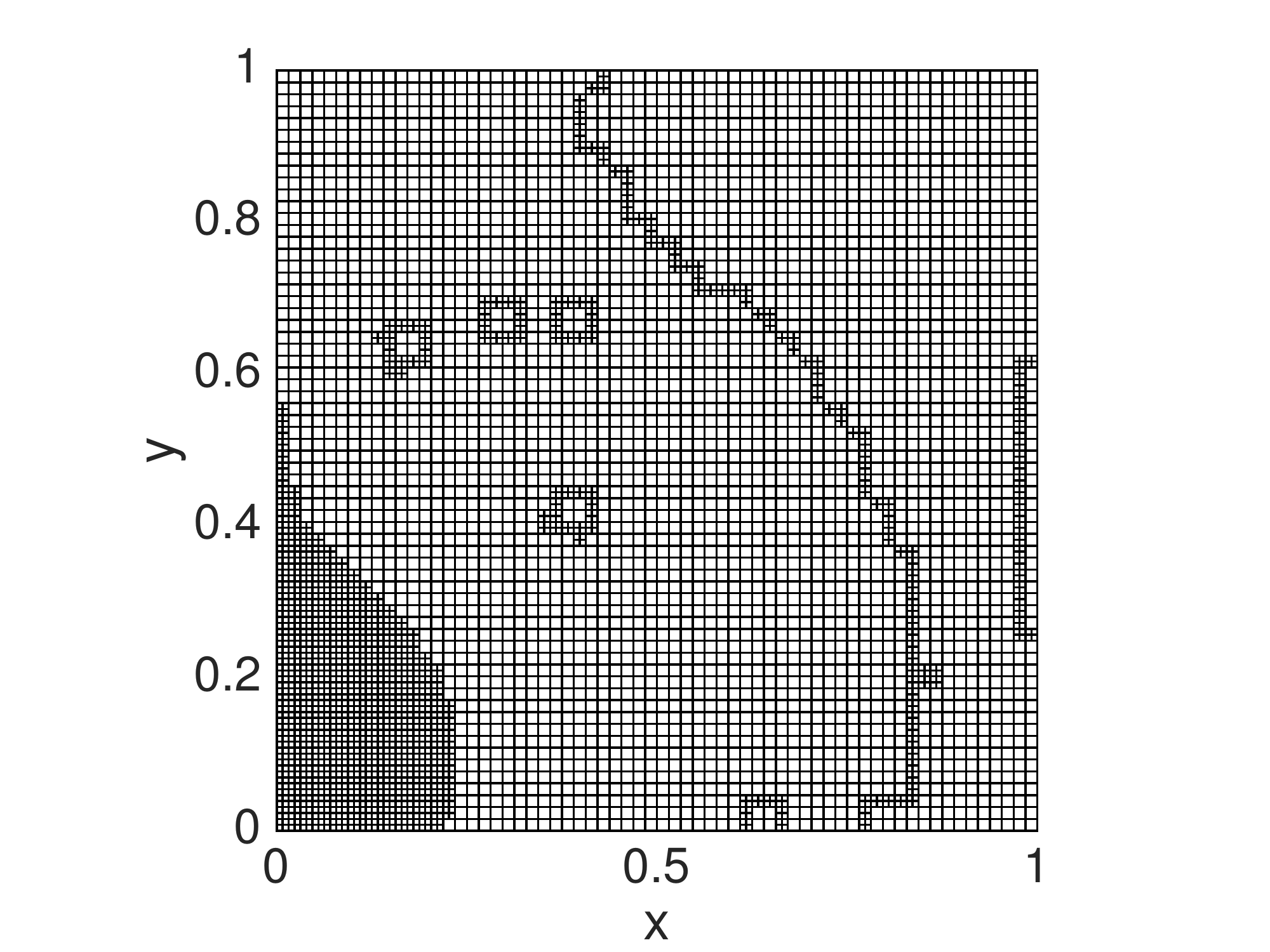}}
    \subfloat[$k=20$]{\label{fig:hankel:mesh:hp:20}\includegraphics[width=0.4\textwidth]{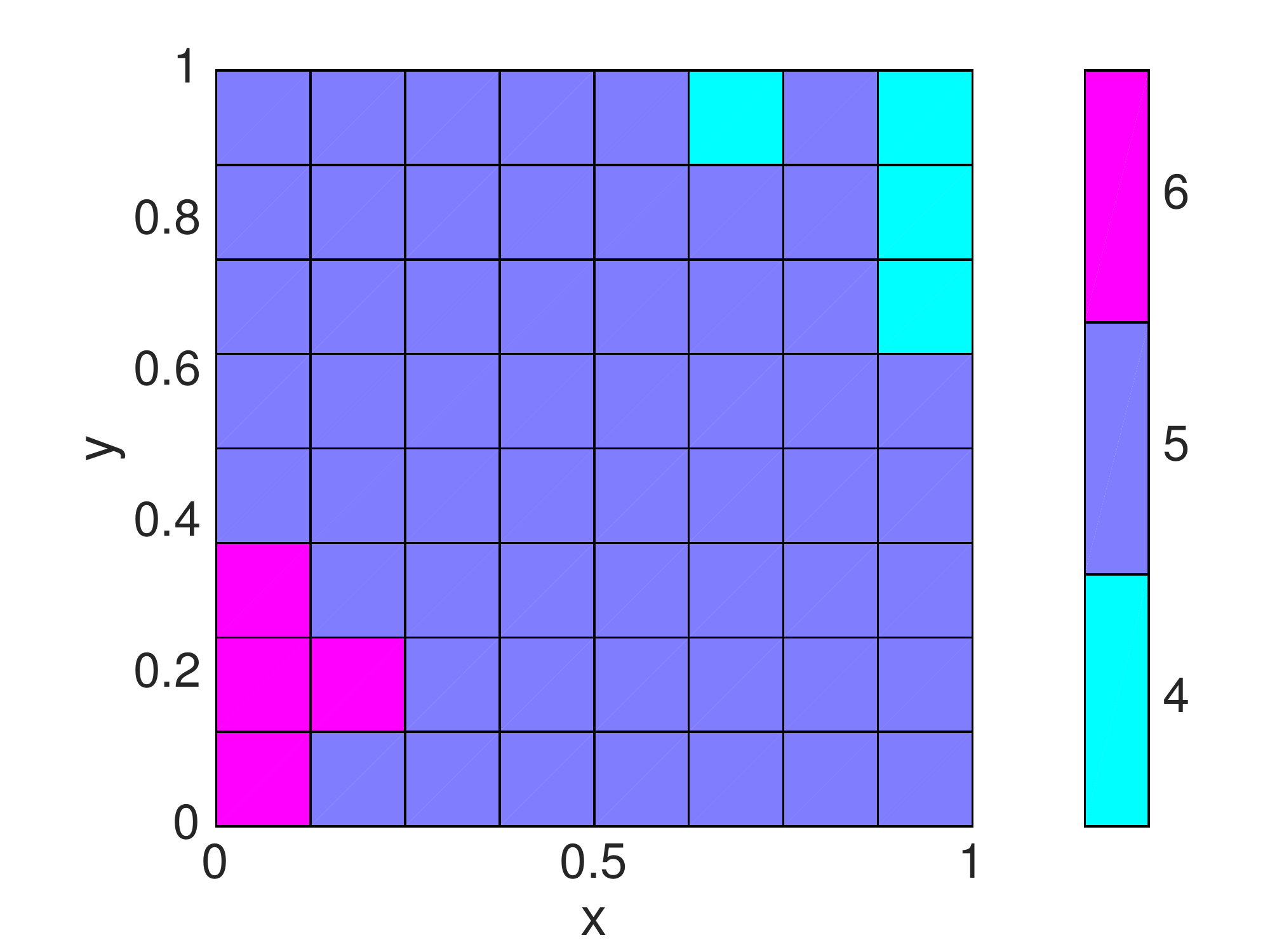}} \\
    \subfloat[$k=50$]{\label{fig:hankel:mesh:h:50}\includegraphics[width=0.4\textwidth]{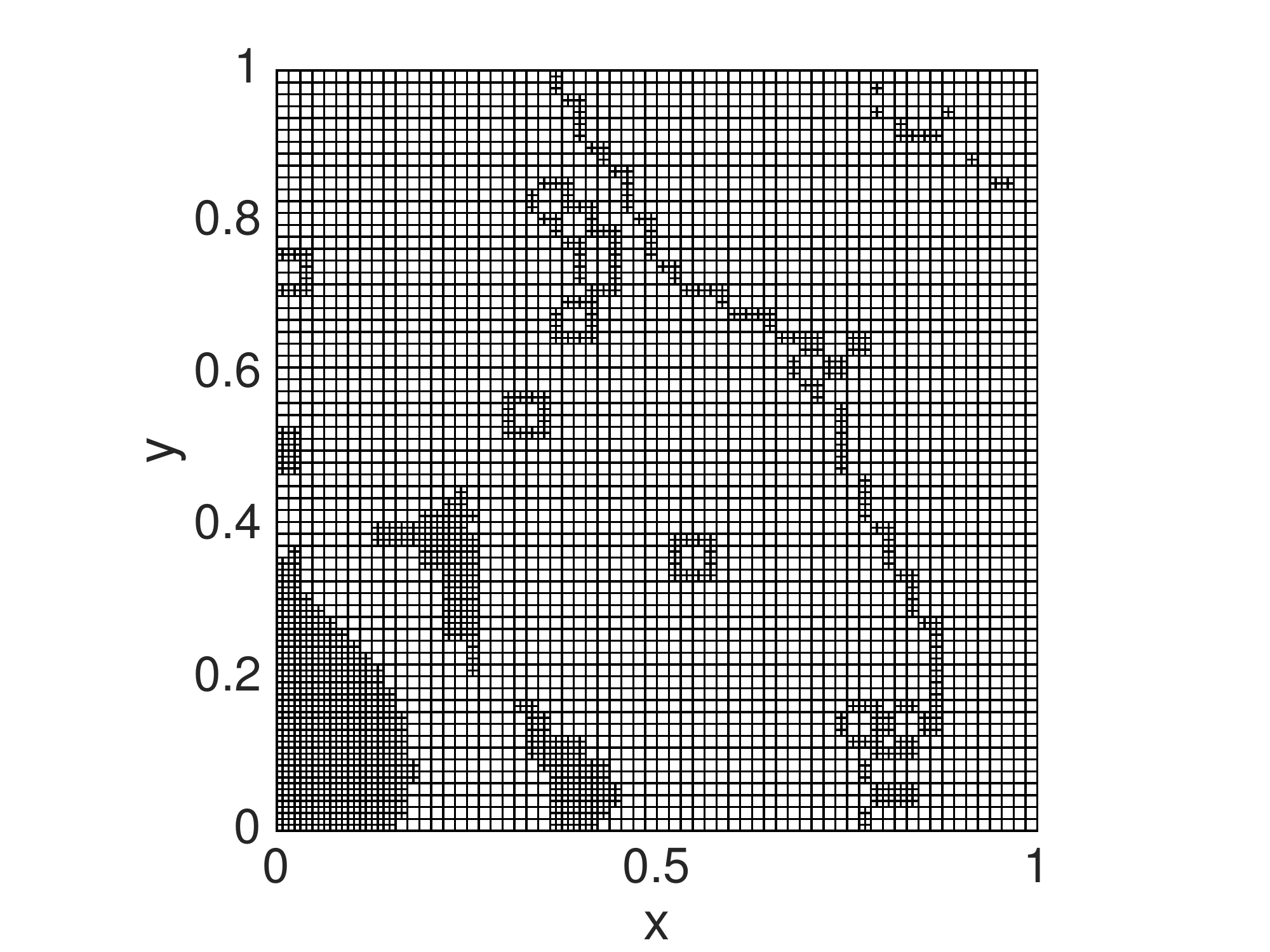}}
    \subfloat[$k=50$]{\label{fig:hankel:mesh:hp:50}\includegraphics[width=0.4\textwidth]{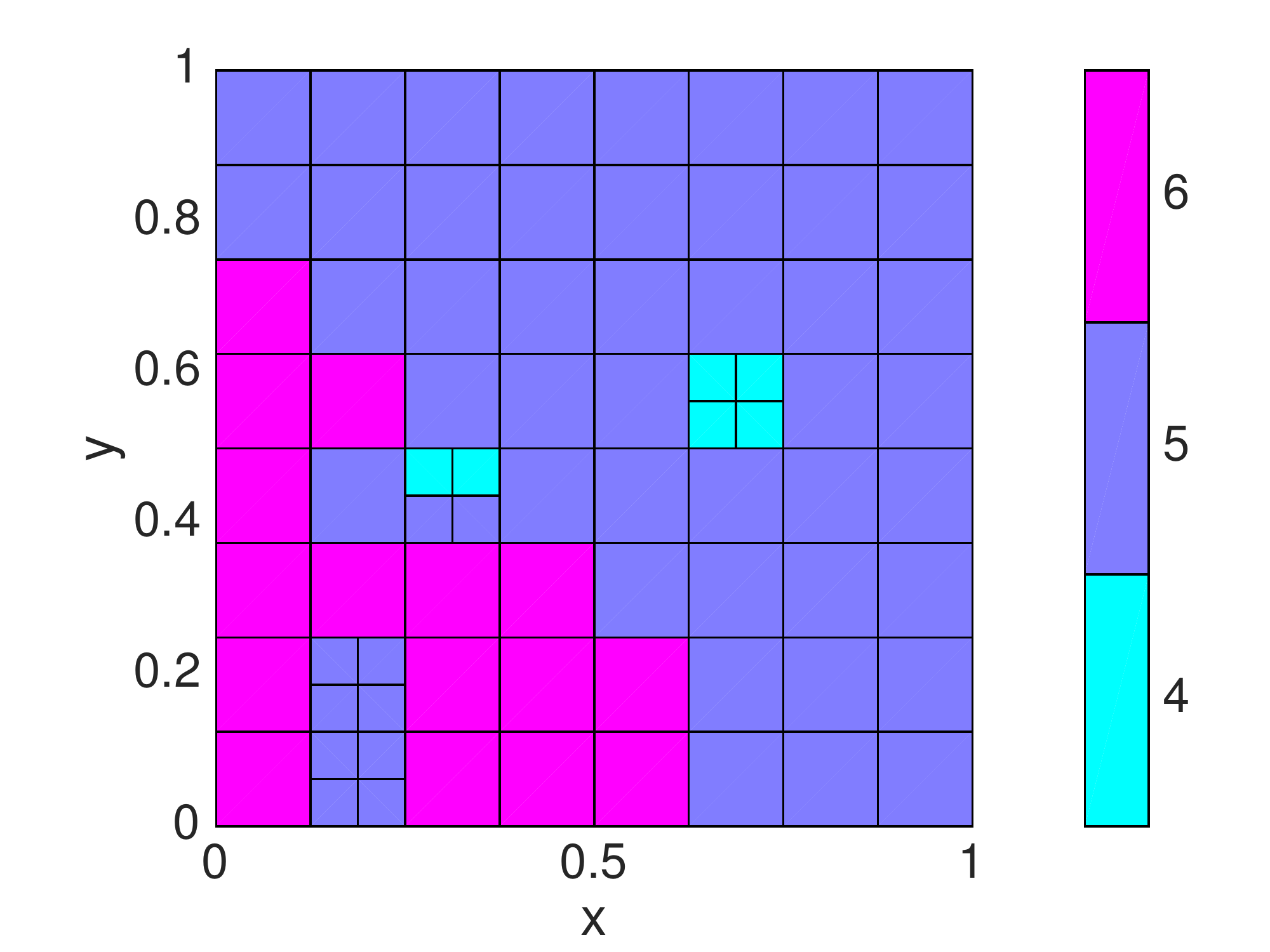}}
    \caption{Example 1: Meshes after 8 \protect\subref{fig:hankel:mesh:h:20} $h$-- and \protect\subref{fig:hankel:mesh:hp:20} $hp$--refinements for wavenumber $k=20$; meshes after 8 \protect\subref{fig:hankel:mesh:h:50} $h$-- and \protect\subref{fig:hankel:mesh:hp:50} $hp$--refinements for wavenumber $k=50$.}
    \label{fig:hankel}
\end{figure}

In this section, we again consider the problem outlined in Section~\ref{section:plane_wave_refine}.
Furthermore, we select the initial mesh to consist of $8\times 8$ uniform square elements and set
$q_K=3$ on each $K\in\mesh$. Firstly, in Figures~\ref{fig:hankel:error:20h} and \ref{fig:hankel:error:50h} 
we compare the relative error in the $L^2$-norm to the number of degrees of freedom in the
TDG space $V_{\vect{p}}(\mathcal{T}_h)$, when $h$--refinement is employed, with the wavenumbers $k=20$ 
and $k=50$, respectively. In each case, we consider the performance of the underlying
adaptive algorithm when both the standard TDG scheme (without direction adaptivity) is employed,
as well as the corresponding method with directional adaptivity; in this latter setting, we 
consider the cases when either directional adaptivity is undertaken on only the elements marked 
for refinement, as well as when it is exploited on every element in the computational
mesh $\mesh$. 
Analogous results are presented in Figures~\ref{fig:hankel:error:20hp} and \ref{fig:hankel:error:50hp} 
in the $hp$--setting, respectively; here, we compare standard $hp$--refinement, with $hp$--adaptivity incorporating
directional adaptivity. In the latter case, different directional adaptivity strategies are considered: 
firstly, directional adaptivity is performed only on elements marked for $p$--refinement;
secondly, directional adaptivity is undertaken on all elements marked for refinement;
finally, directional adaptivity is applied to every element in $\mesh$. In the $hp$--setting
we observe exponential convergence of the error as the finite element space is adaptively
enriched: indeed, on a linear-log scale, the convergence lines are roughly straight. Thereby,
it is clear that the exploitation of the proposed $hp$--refinement algorithm, with directional
adaptivity, leads to a significant reduction in the $L^2$-norm of error, for a given number of degrees of
freedom, when compared to the same quantity computed with $h$--refinement alone; cf. Figure~\ref{fig:hankel:errorcompare}.

In both the $h$-- and $hp$--refinement cases, we generally observe that the error is decreased when
directional refinement is employed. Moreover, it is evident in the $hp$--setting that 
selecting more elements for directional refinement generally leads to a smaller error, for a given
number of degrees of freedom; this is particularly noticeable in the case when $k=50$. 
In Figures~\ref{fig:hankel:eff:20h}, \ref{fig:hankel:eff:20hp}, \ref{fig:hankel:eff:50h}, 
and \ref{fig:hankel:eff:50hp} we plot the effectivity indices for each of the above refinement 
strategies for the case when $k=20,50$; here, we observe that they remain roughly constant 
during adaptive $h$--/$hp$--mesh refinement, and are roughly the same for the two different wavenumbers, 
with the notable exception of the pre-asymptotic region for $k=50$. 

Finally, in Figures~\ref{fig:hankel:mesh:h:20}--\ref{fig:hankel:mesh:hp:50}, we show the meshes 
after 8 $h$-- and $hp$--refinements, with directional adaptivity employed on all elements, for 
both $k=20$ and $k=50$; here, the $hp$--meshes show the effective polynomial degree $q_K$ for each 
element. Given the smoothness of the analytical solution on $\Omega$, we observe that the resulting 
computational meshes are almost uniform; indeed, in the $hp$--setting almost uniform $p$--refinement
has been undertaken.

\subsubsection{Example 2 --- Singular solution}
\begin{figure}[pt]
    \configfigure
    \subfloat[$k=20$; $h$--refinement]{\label{fig:besselsin:error:20h}\includegraphics[width=0.4\textwidth]{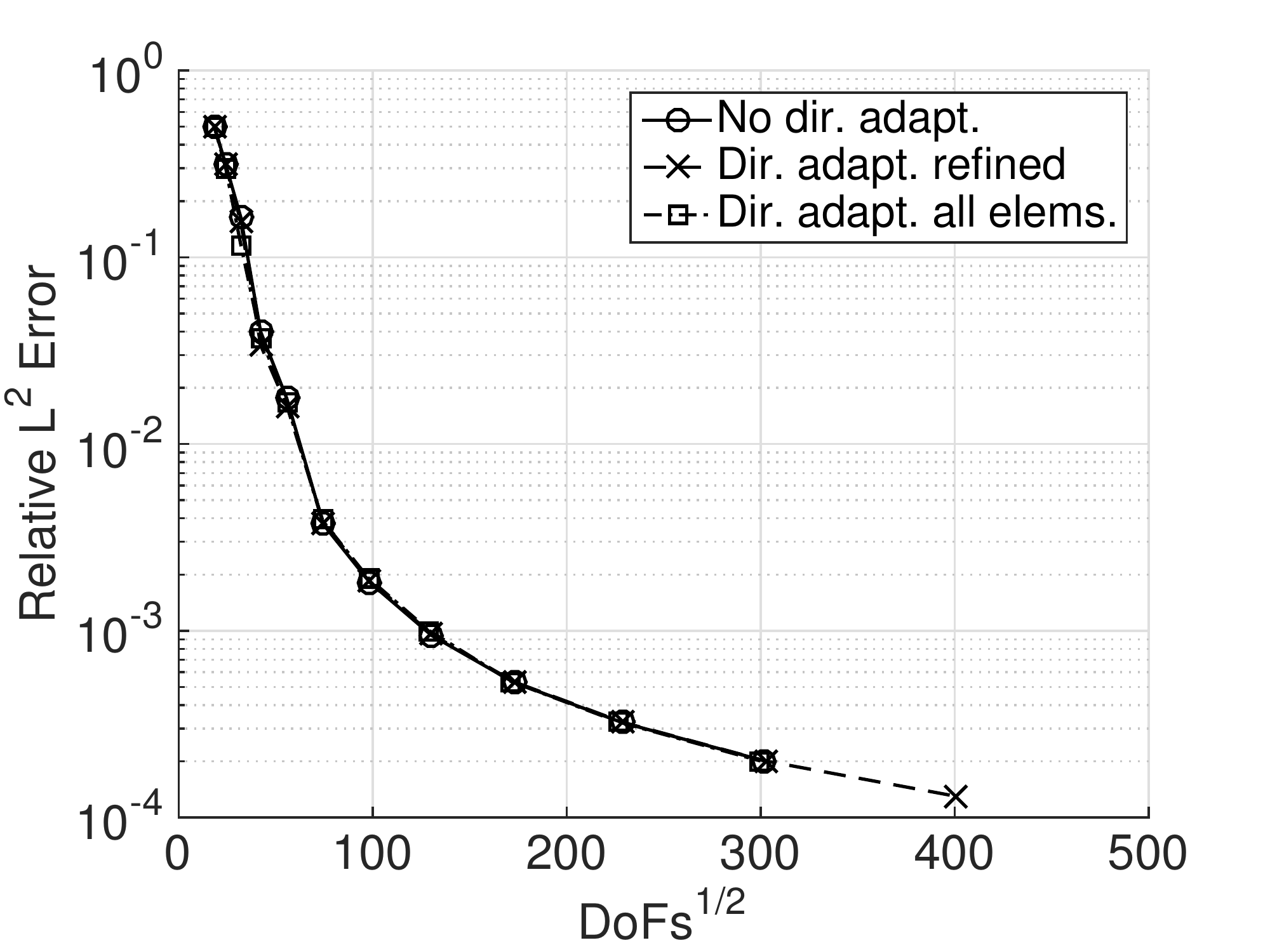}}
    \subfloat[$k=20$; $h$--refinement]{\label{fig:besselsin:eff:20h}\includegraphics[width=0.4\textwidth]{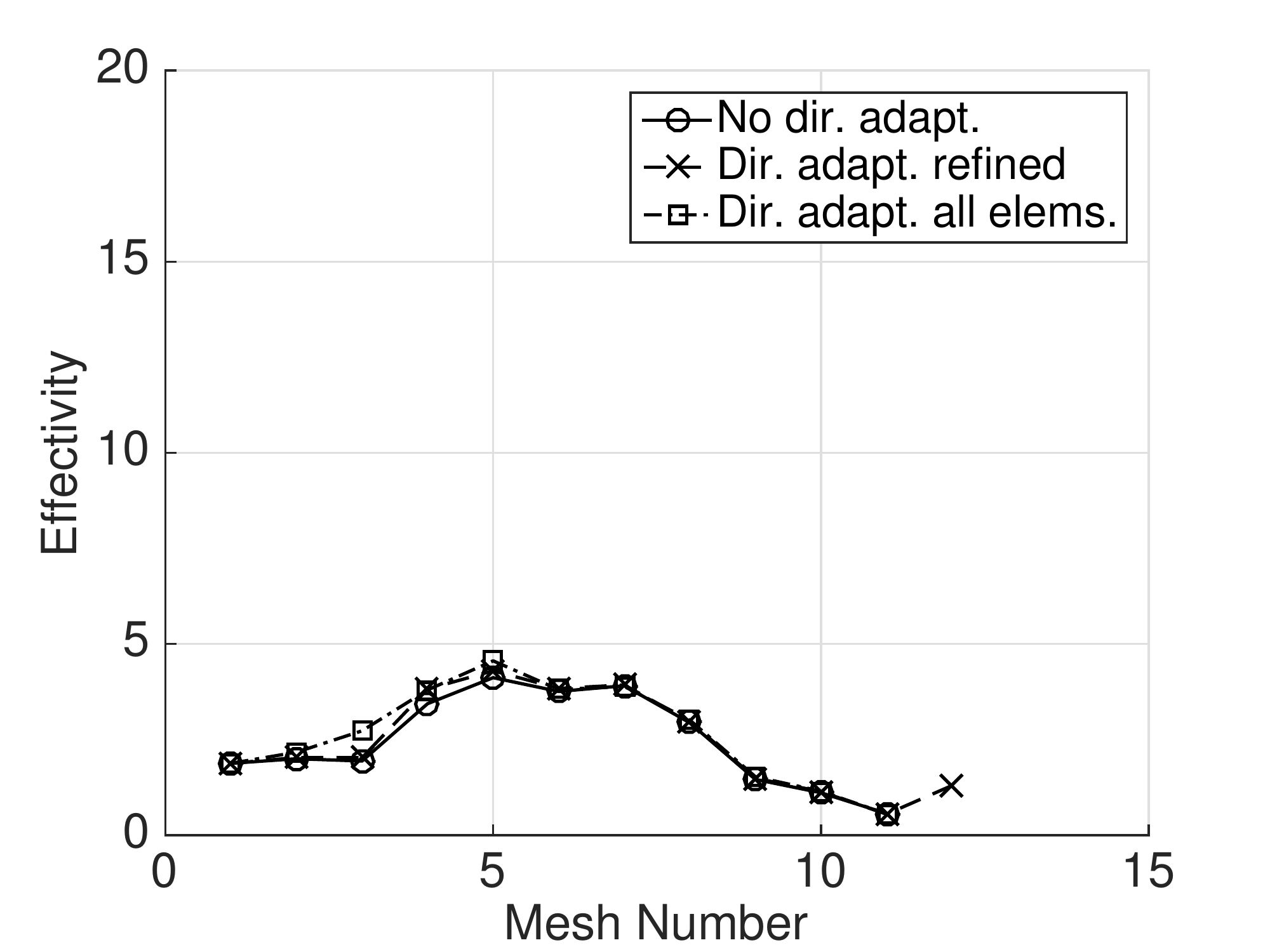}} \\
    \subfloat[$k=20$; $hp$--refinement]{\label{fig:besselsin:error:20hp}\includegraphics[width=0.4\textwidth]{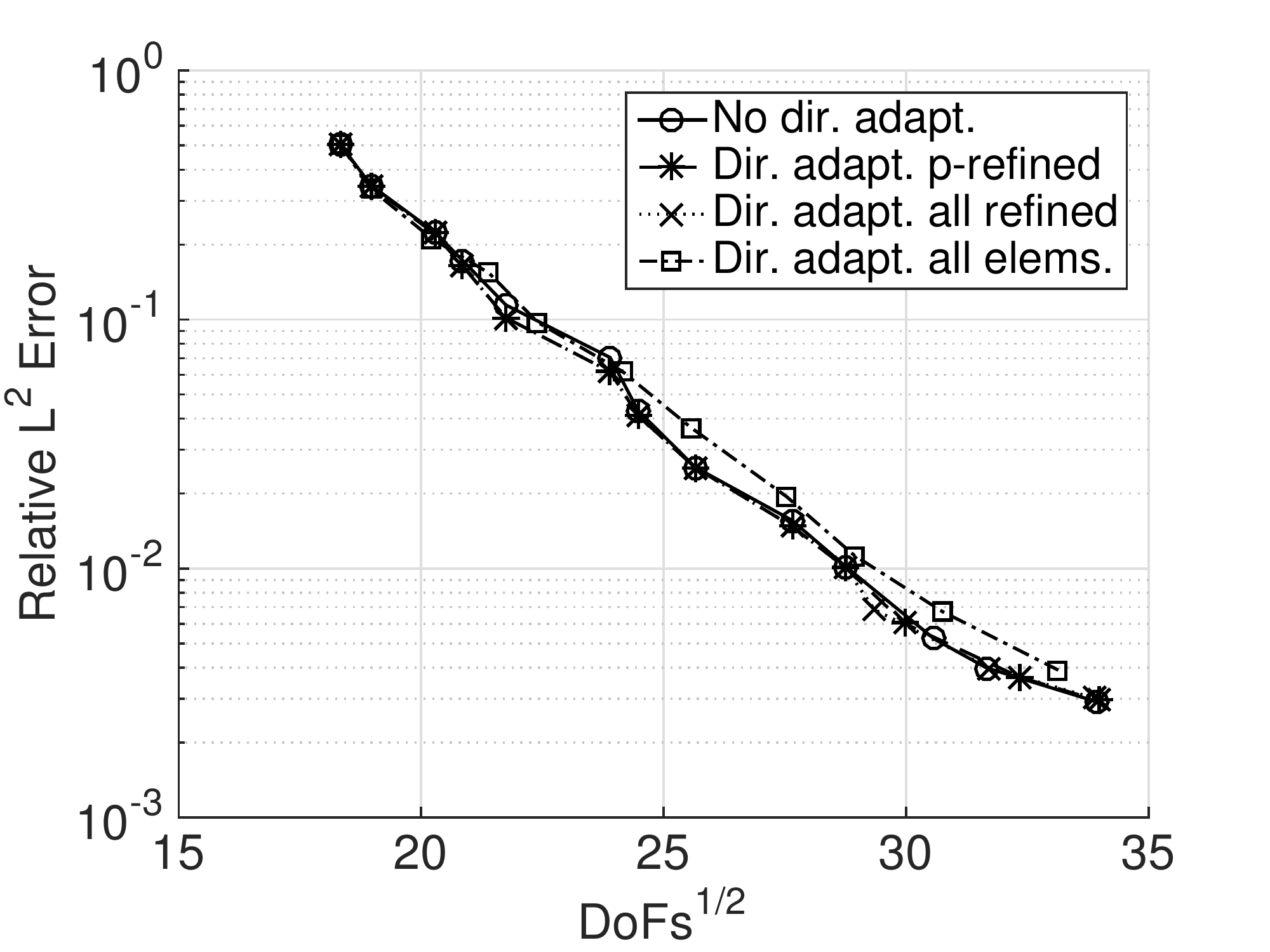}}
    \subfloat[$k=20$; $hp$--refinement]{\label{fig:besselsin:eff:20hp}\includegraphics[width=0.4\textwidth]{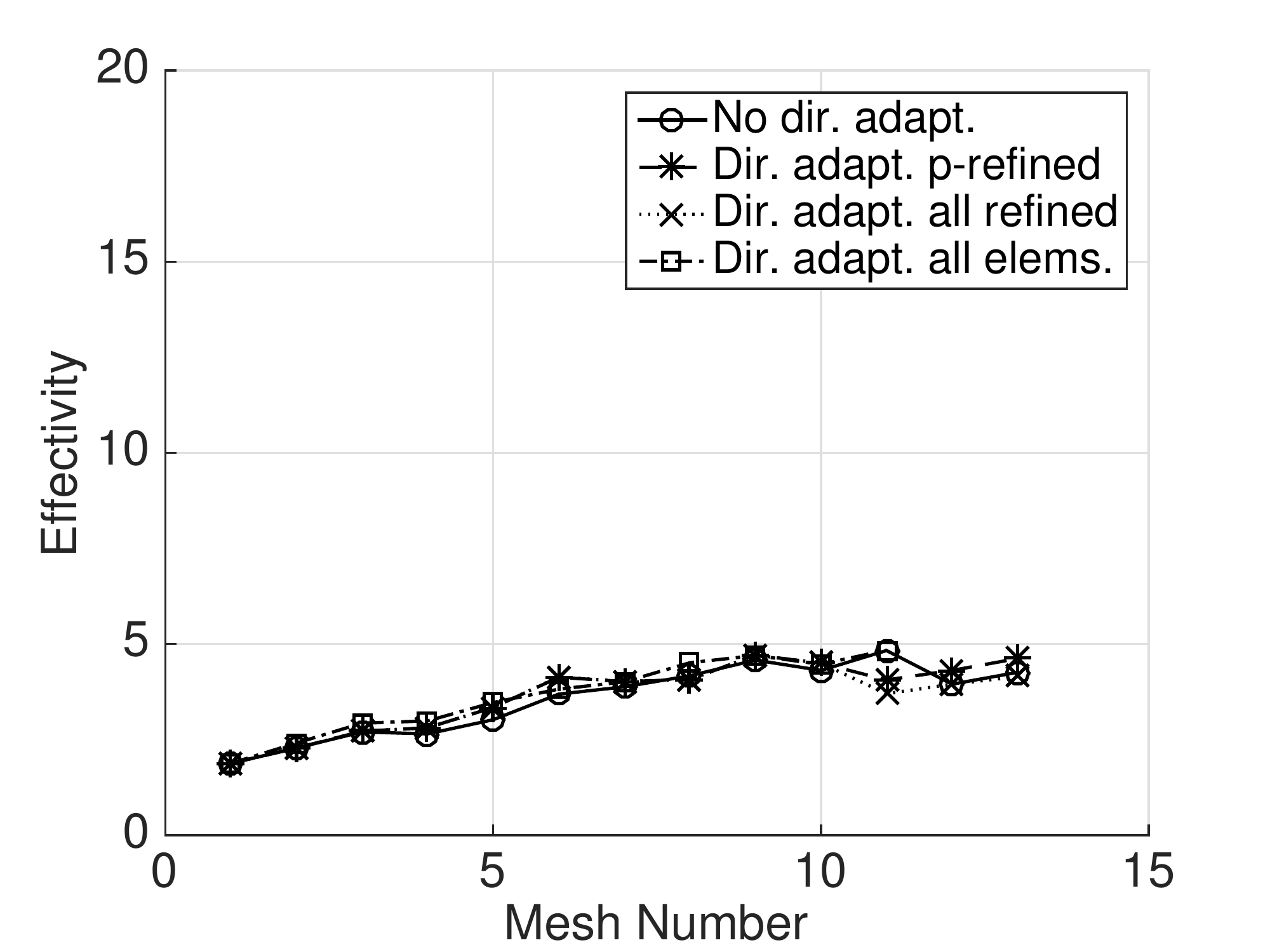}} \\
    \subfloat[$k=50$; $h$--refinement]{\label{fig:besselsin:error:50h}\includegraphics[width=0.4\textwidth]{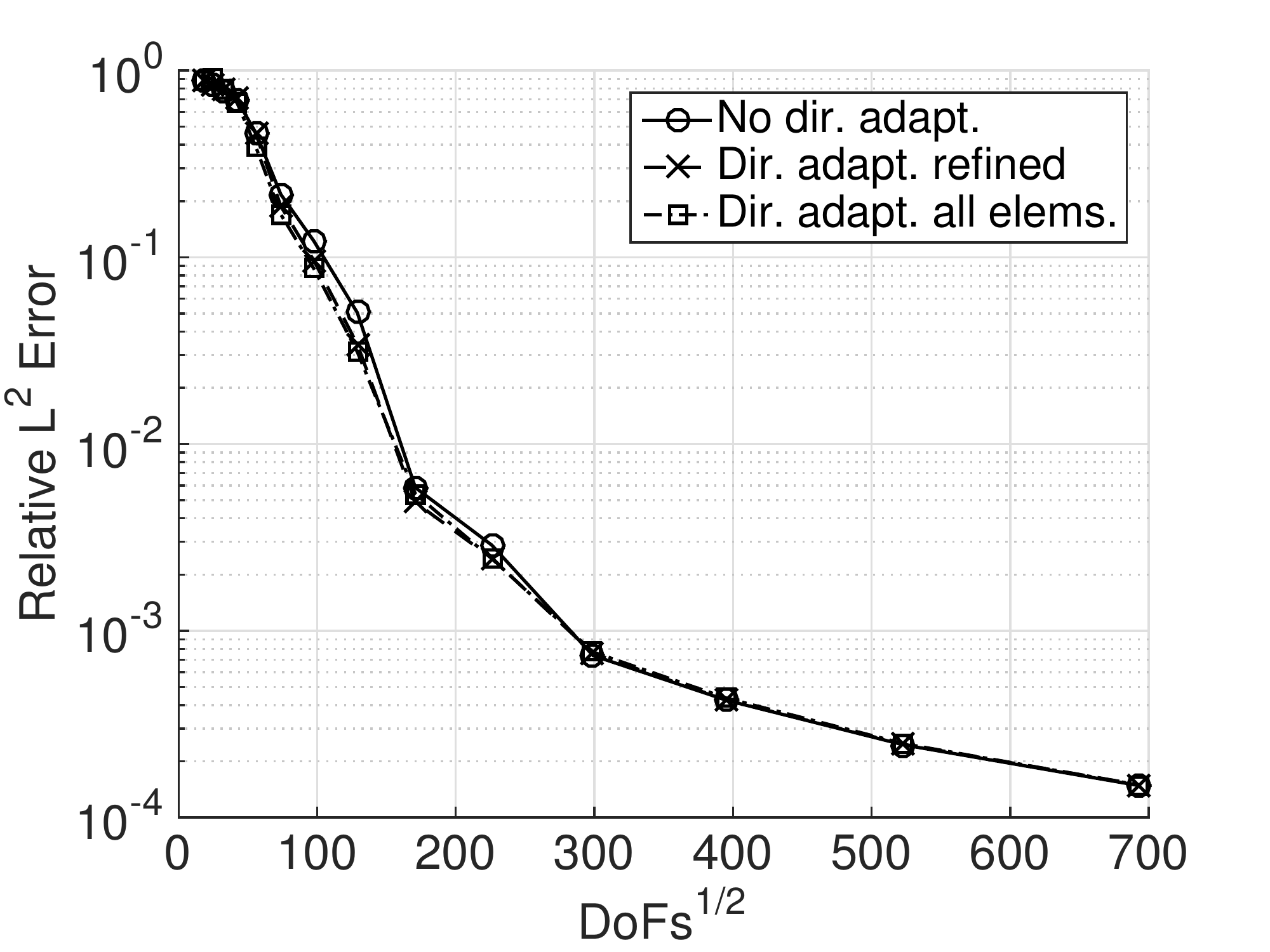}}
    \subfloat[$k=50$; $h$--refinement]{\label{fig:besselsin:eff:50h}\includegraphics[width=0.4\textwidth]{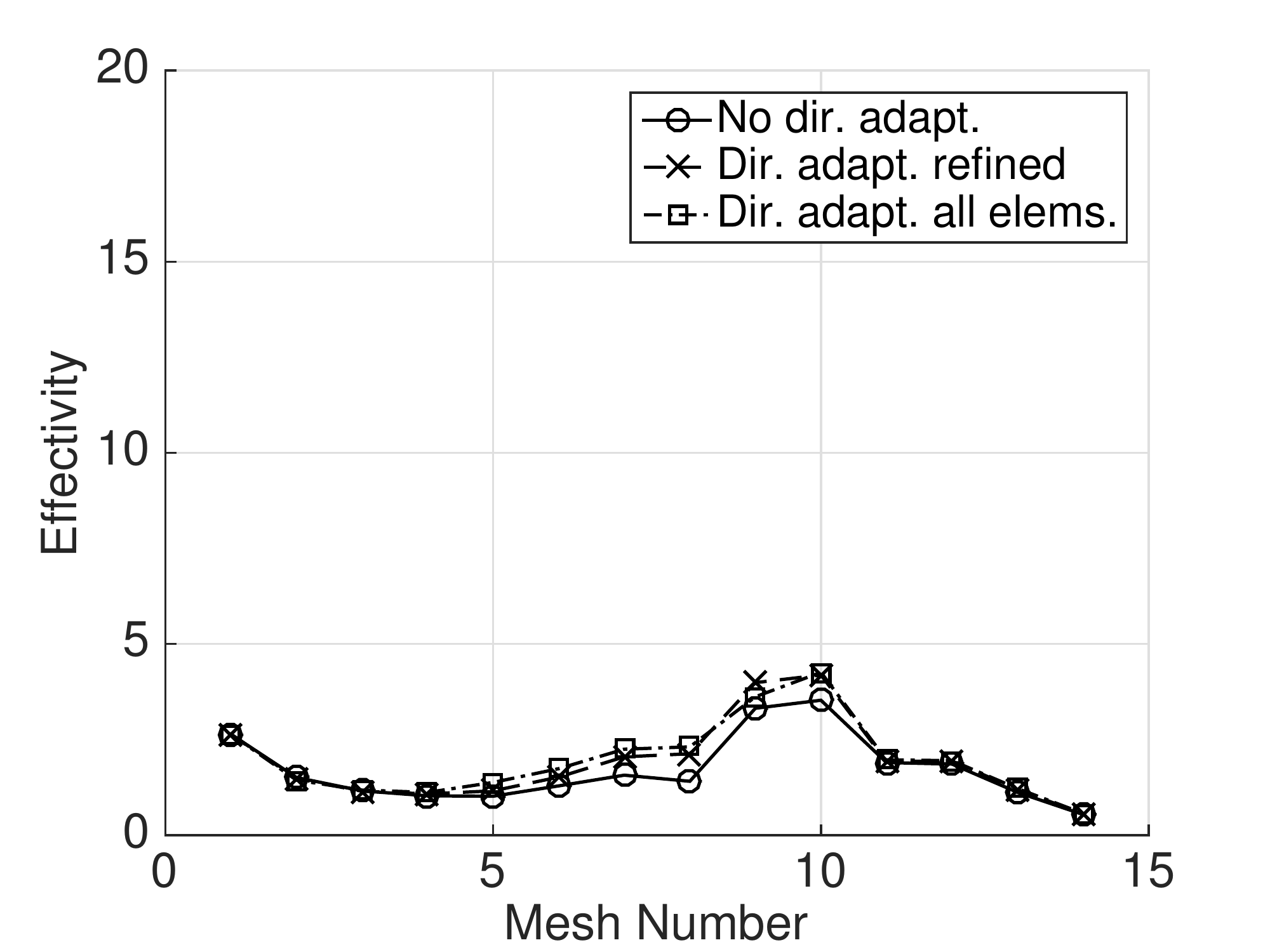}} \\
    \subfloat[$k=50$; $hp$--refinement]{\label{fig:besselsin:error:50hp}\includegraphics[width=0.4\textwidth]{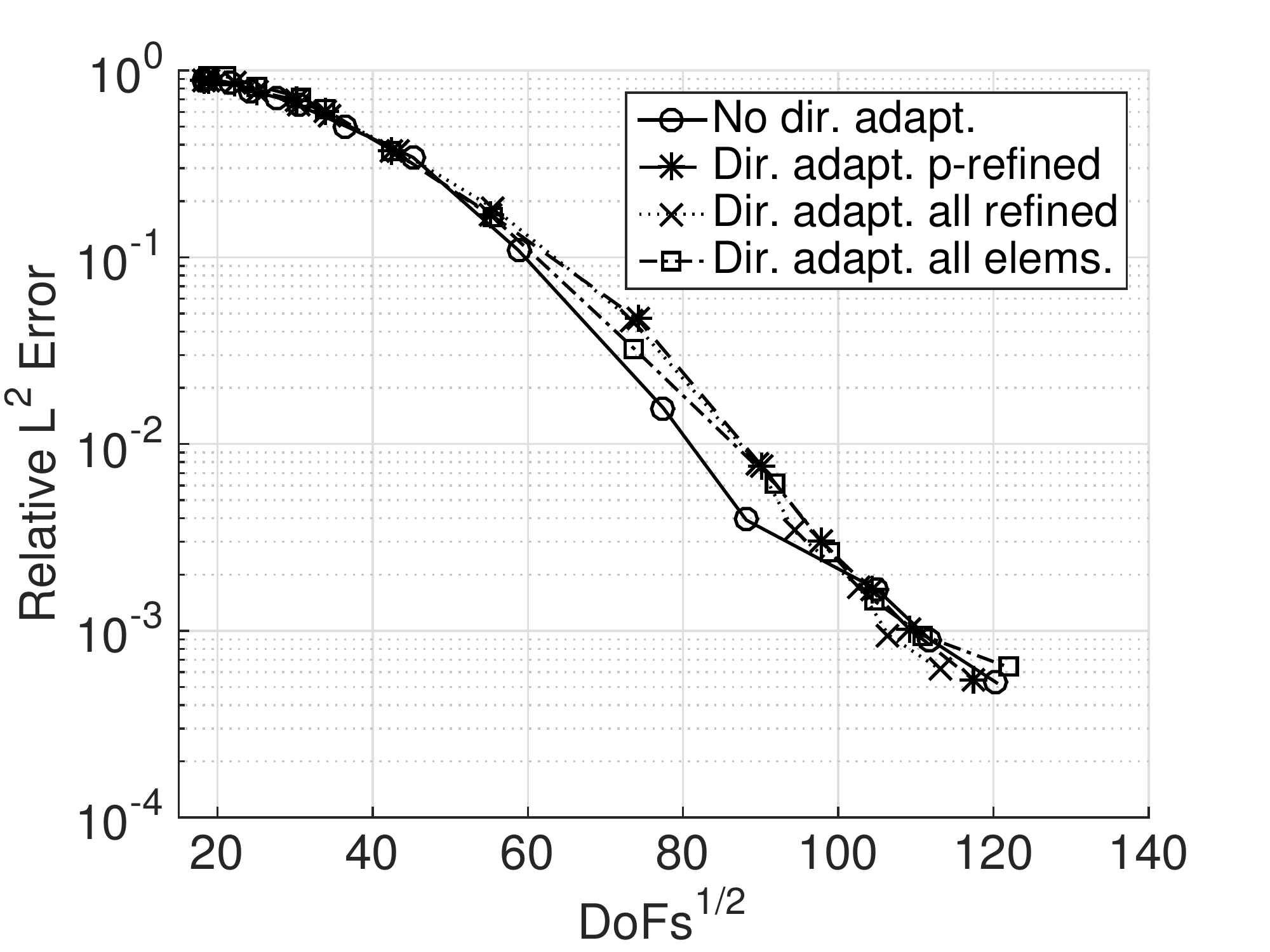}}
    \subfloat[$k=50$; $hp$--refinement]{\label{fig:besselsin:eff:50hp}\includegraphics[width=0.4\textwidth]{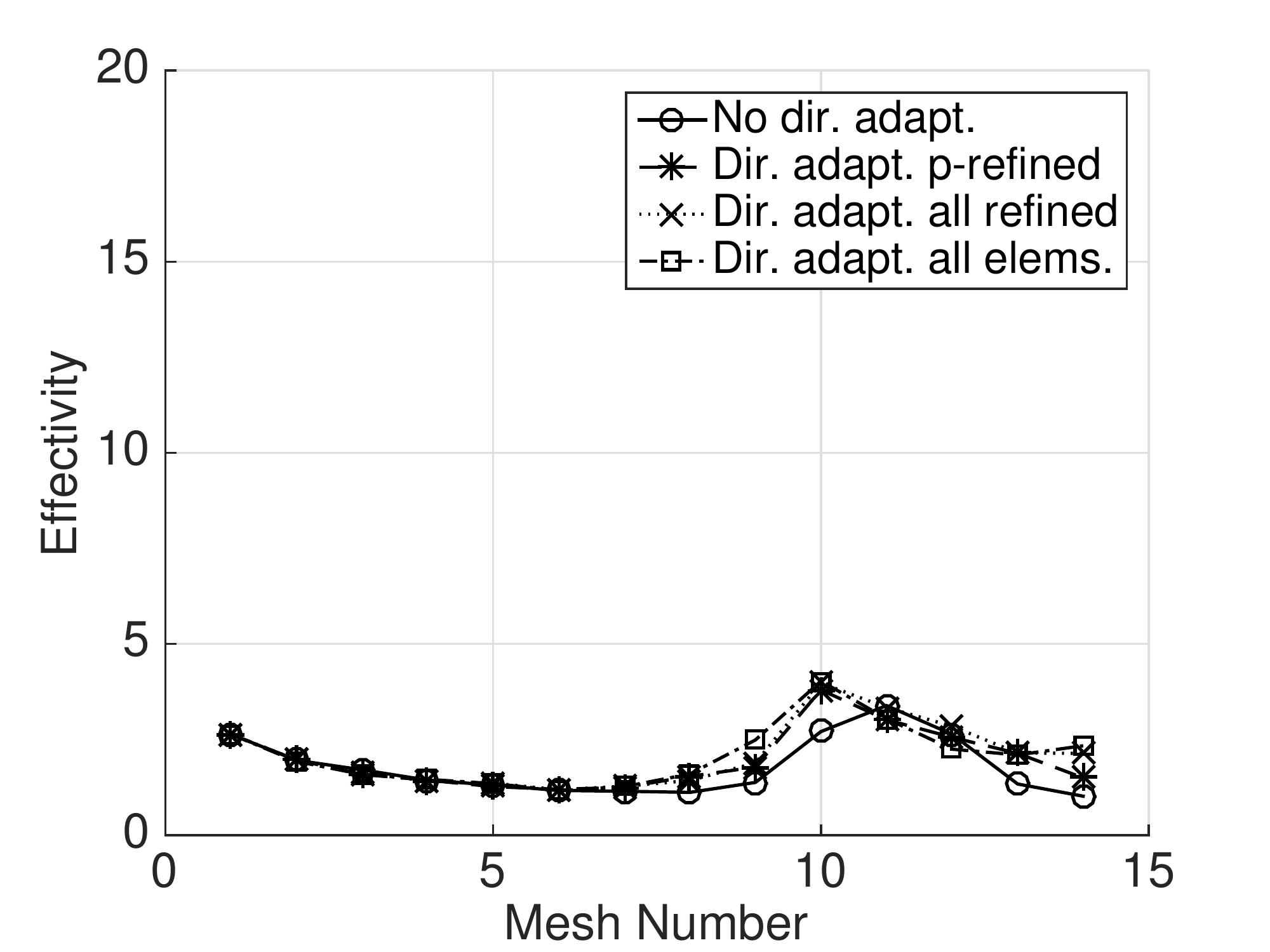}}
    \caption{Example 2: \protect\subref{fig:besselsin:error:20h} $L^2$-error and \protect\subref{fig:besselsin:eff:20h} Effectivity index for $h$--refinement with wavenumber $k=20$; \protect\subref{fig:besselsin:error:20hp} $L^2$-error and \protect\subref{fig:besselsin:eff:20hp} Effectivity index for $hp$--refinement with $k=20$; \protect\subref{fig:besselsin:error:50h} $L^2$-error and \protect\subref{fig:besselsin:eff:50h} Effectivity index for $h$--refinement with $k=50$; \protect\subref{fig:besselsin:error:50hp} $L^2$-error and \protect\subref{fig:besselsin:eff:50hp} Effectivity index for $hp$--refinement with $k=50$.}
    \label{fig:besselsin}
\end{figure}
\begin{figure}[pt]
    \configfigure
    \subfloat[$k=20$]{\label{fig:besselsin:mesh:h:20}\includegraphics[width=0.4\textwidth]{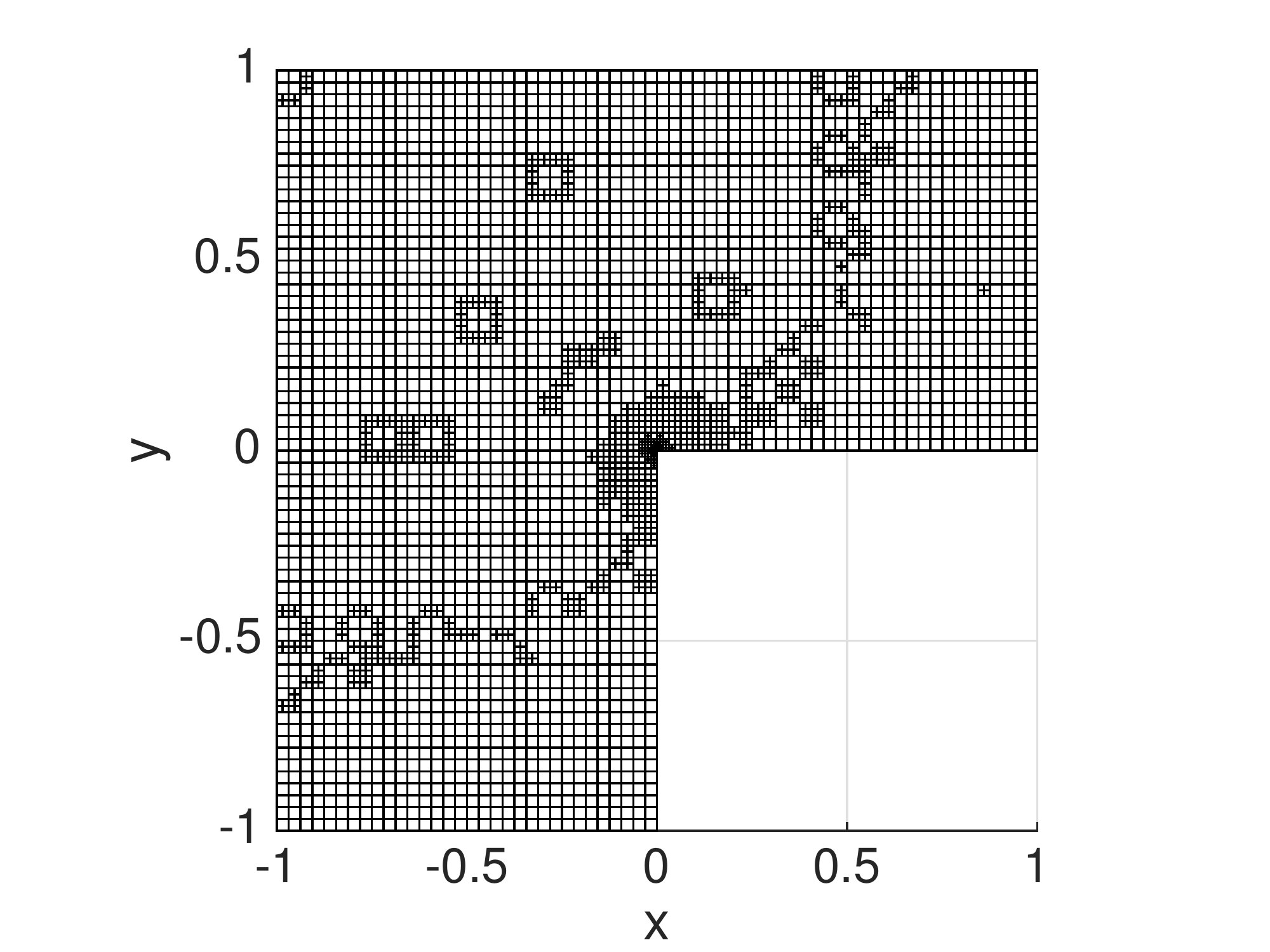}}
    \subfloat[$k=20$]{\label{fig:besselsin:mesh:hp:20}\includegraphics[width=0.4\textwidth]{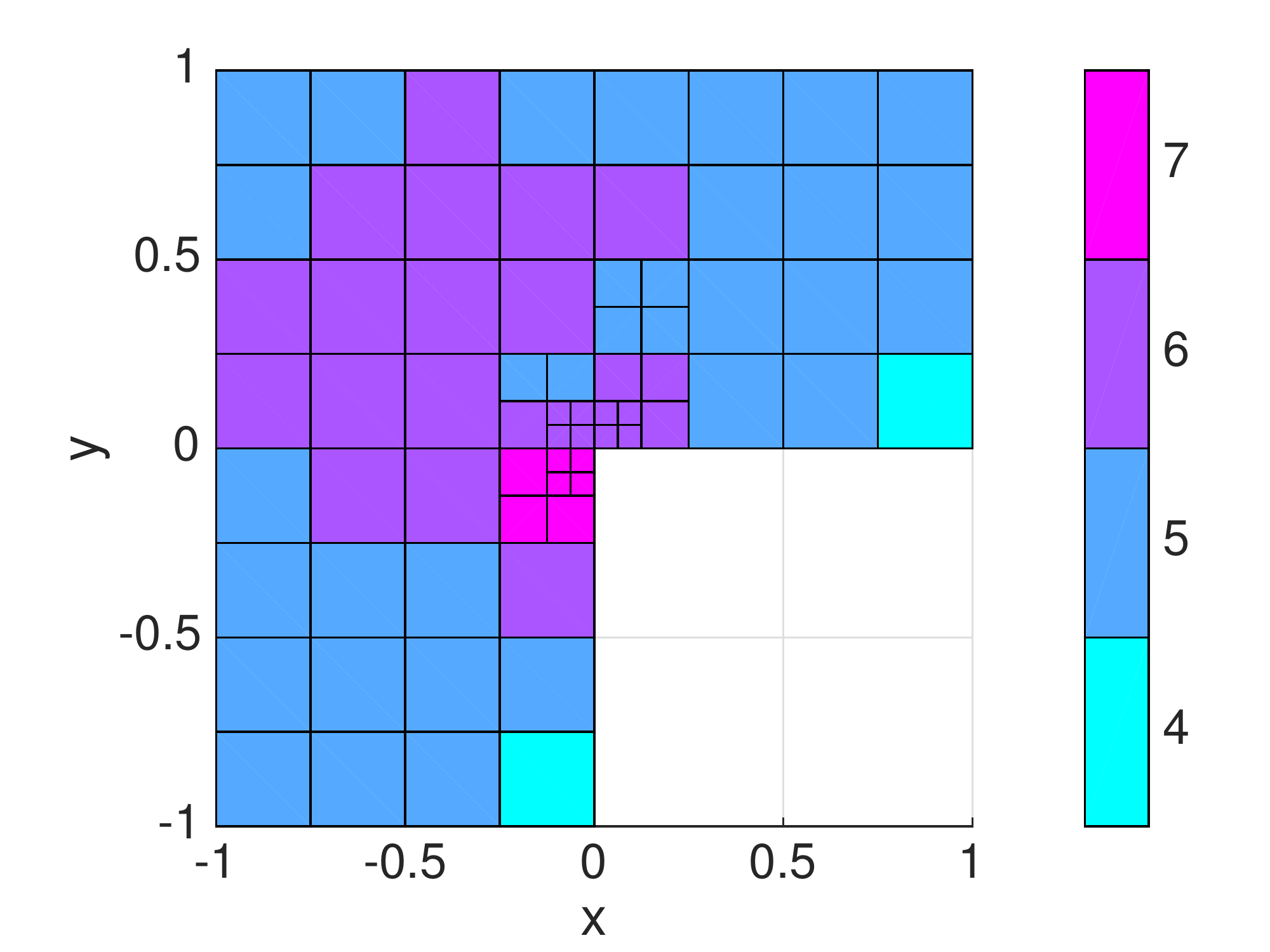}}\\
    \subfloat[$k=50$]{\label{fig:besselsin:mesh:h:50}\includegraphics[width=0.4\textwidth]{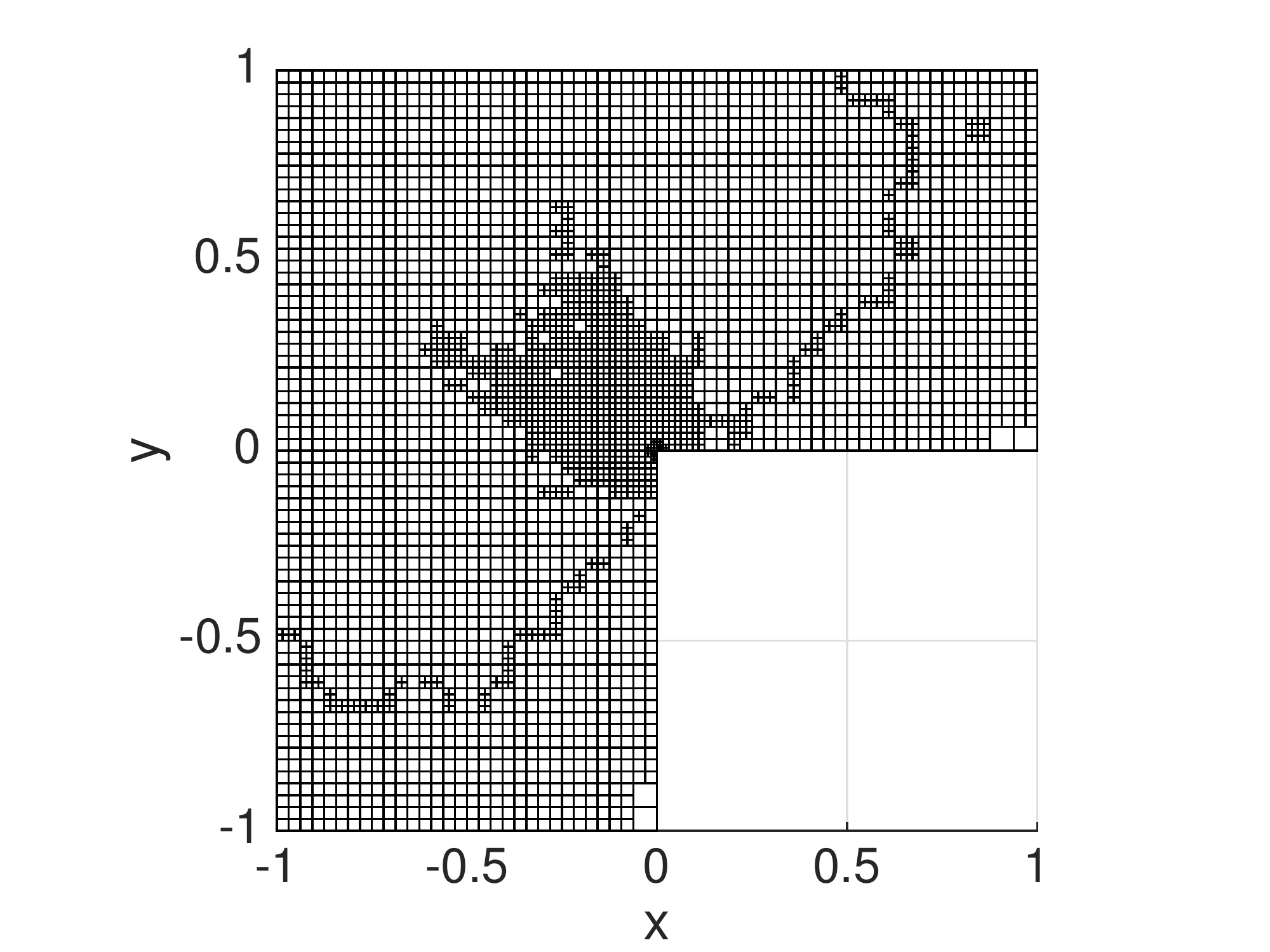}}
    \subfloat[$k=50$]{\label{fig:besselsin:mesh:hp:50}\includegraphics[width=0.4\textwidth]{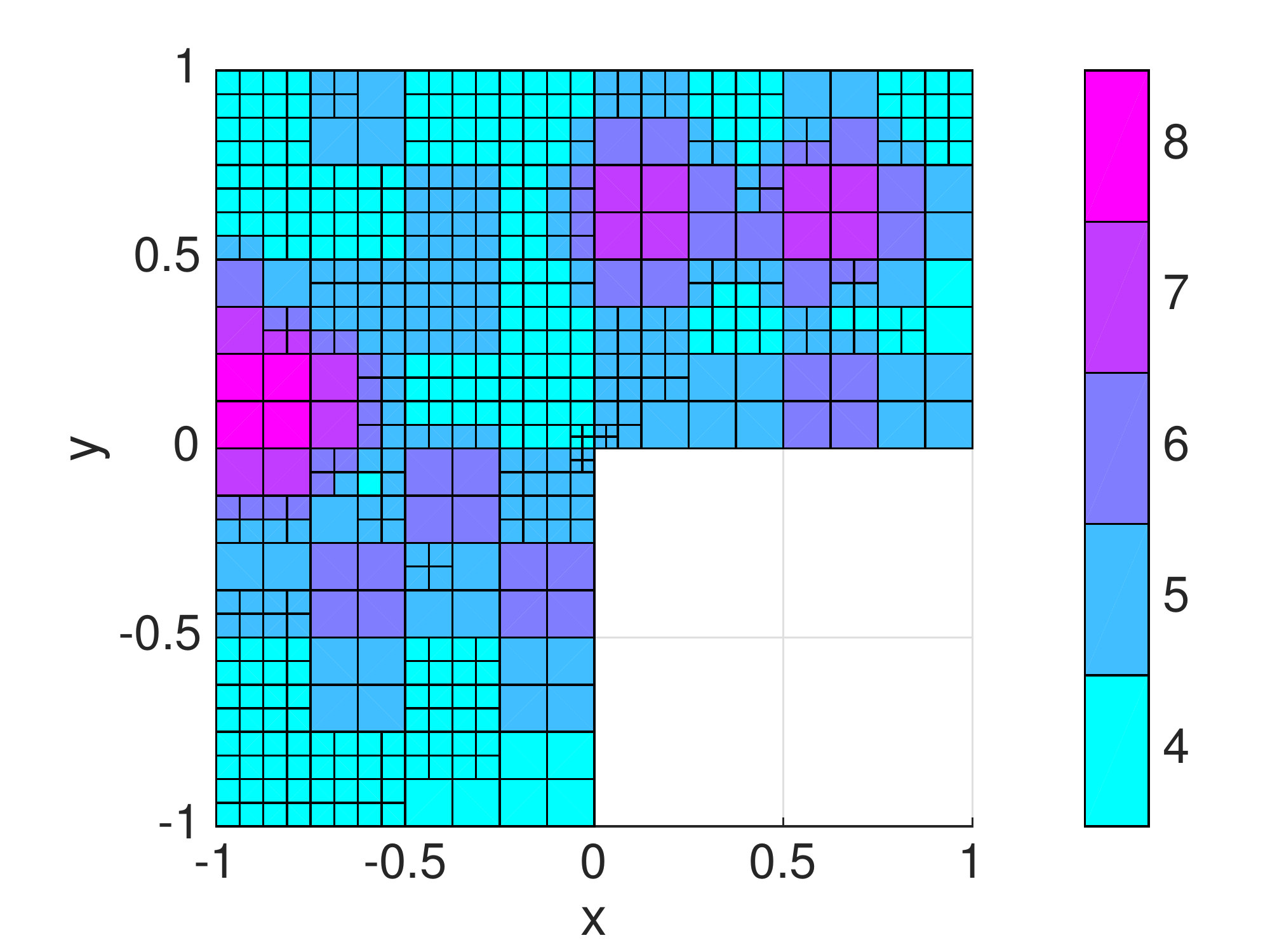}}
    \caption{Example 2: Meshes after 8 \protect\subref{fig:besselsin:mesh:h:20} $h$-- and \protect\subref{fig:besselsin:mesh:hp:20} $hp$--refinements for wavenumber $k=20$; meshes after 8 \protect\subref{fig:besselsin:mesh:h:50} $h$-- and \protect\subref{fig:besselsin:mesh:hp:50} $hp$--refinements for wavenumber $k=50$.}
    \label{fig:besselsin:mesh}
\end{figure}

In this second example, we consider problem~\eqref{eqn:helmholtz} posed  on the L-shaped domain $\Omega=(-1,1)^2\setminus(0,1)\times(-1,1)$, $\Gamma_R=\partial\Omega$, and $\Gamma_D\equiv\emptyset$, with Robin boundary condition $g_R$ selected so that the analytical solution is given, in polar coordinates $(r,\varphi)$, by
\[
    u(r,\theta) = \mathcal{J}_{\nicefrac23}(kr)\sin(\nicefrac{2\theta}{3});
\]
we note that the gradient of $u$ has a singularity at the origin. 

As in the previous example, we again compare the performance of the $h$-- and $hp$--adaptive refinement algorithms, both in the standard setting, as well as when directional adaptivity is employed; here, we again consider the analogous directional refinement strategies employed in Section~\ref{sec-hprefine-hankel}. To this end, in Figures~\ref{fig:besselsin:error:20h} and \ref{fig:besselsin:error:50h} we compare the relative error in the $L^2$-norm with the number of degrees of freedom in the TDG space $V_{\vect{p}}(\mathcal{T}_h)$ when $h$--refinement is employed for $k=20$ and $k=50$, respectively; the respective convergence plots in the $hp$--setting are given in Figures~\ref{fig:besselsin:error:20h} and \ref{fig:besselsin:error:50h}. Here, we observe that although exploiting $hp$--refinement leads to exponential convergence of the relative $L^2$-norm of the error as $V_{\vect{p}}(\mathcal{T}_h)$ is enriched, in both the $h$-- and $hp$--settings, we observe that the magnitude of the error, computed both with and without directional refinement, is roughly identical; i.e., directional refinement does not lead to any reduction in the computed TDG solution when either $h$--/$hp$--refinement is employed. We note that, for this particular problem, this behaviour is expected, since the error in the computed TDG solution is dominated by the influence of the singularity at the origin, rather than local wave propagation. 

In Figures~\ref{fig:besselsin:eff:20h}, \ref{fig:besselsin:eff:20hp}, \ref{fig:besselsin:eff:50h}, and \ref{fig:besselsin:eff:50hp} we plot the effectivity indices when both $h$-- and $hp$--refinement is employed for the case when $k=20,50$. In all cases, we observe that the effectivity indices are roughly constant for this singular problem, though when $h$--refinement is employed, on highly refined meshes, we see a slight drop in the computed effectivity indices. Finally, in Figures~\ref{fig:besselsin:mesh:h:20}--\ref{fig:besselsin:mesh:hp:50}, we show the meshes after 8 $h$-- and $hp$--refinements, with direction adaptivity employed on all elements, for both $k=20$ and $k=50$. As we would expect, in both the $h$-- and $hp$--settings, mesh subdivision is concentrated in the vicinity of the singularity located at the origin; away from this region, the $h$--refinement algorithm employs almost uniform mesh subdivision, while the $hp$--refinement strategy employs the necessary combination of local mesh refinement and local polynomial enrichment, as required, to reduce the error in the computed TDG solution.

\subsubsection{Example 3 --- Transmission/internal reflection}

\begin{figure}[t]
    \configfigure
    \subfloat[$\theta_i=29\degrees$]{\label{fig:reflection:anal}\includegraphics[width=0.4\textwidth]{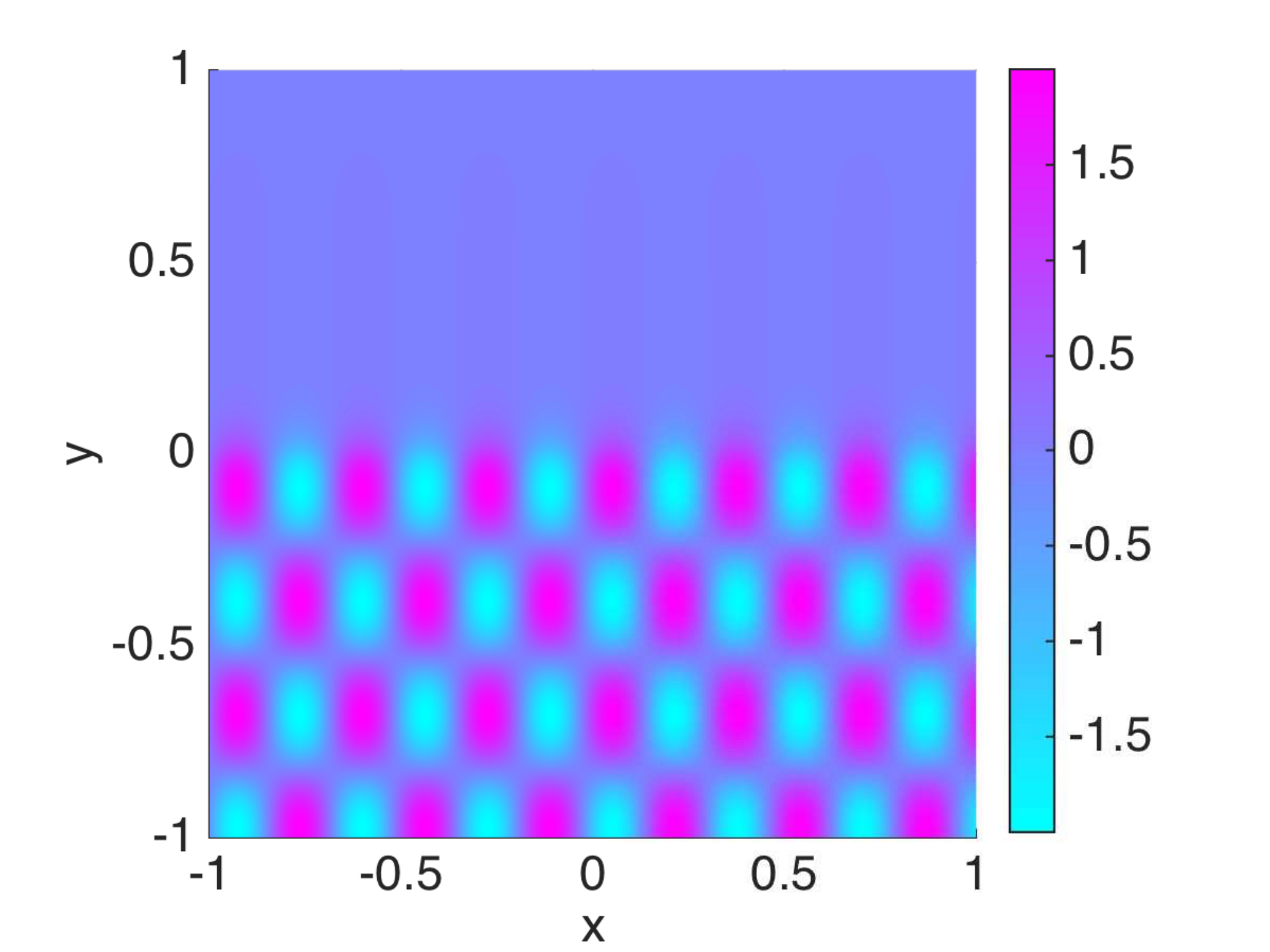}}
    \subfloat[$\theta_i=69\degrees$]{\label{fig:refraction:anal}\includegraphics[width=0.4\textwidth]{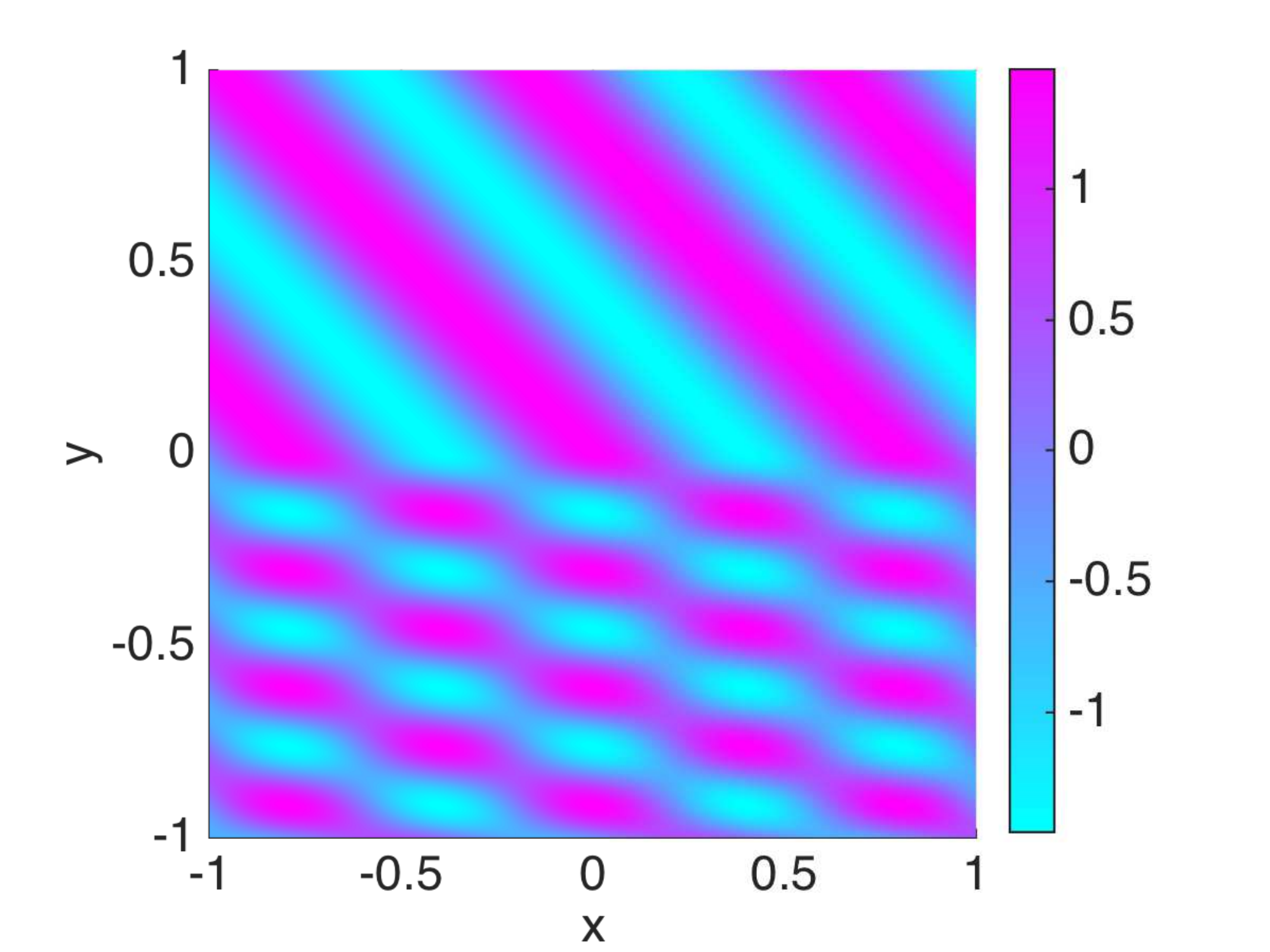}}
    \caption{Example 3: Analytical solutions (real part) when \protect\subref{fig:reflection:anal} $\theta_i=29\degrees$ resulting in internal reflection, and \protect\subref{fig:refraction:anal} $\theta_i=69\degrees$ resulting in refraction.}
    \label{fig:reflection_anal}
\end{figure}

We now consider the case of transmission and internal reflection of a plane wave across a fluid-fluid interface in the domain $\Omega=(-1,1)^2$, $\Gamma_R\equiv \emptyset$, and $\Gamma_D=\partial\Omega$, with two different refractive indices, cf.~\cite[Section 6.3]{Kapita2015}. The interface between the two regions is located at $y=0$; in this setting the wavenumber $k$ is given by the piecewise constant function
\[
    k(x,y) = \begin{cases}
        k_1 \coloneqq \omega n_1 & \text{if } y \leq 0,\\
        k_2 \coloneqq \omega n_2 & \text{if } y > 0,
    \end{cases}
\]
where, we select $\omega=11$, $n_1=2$, and $n_2=1$. Throughout this section we impose an appropriate inhomogeneous Dirichlet boundary condition, so that the analytical solution $u$ to \eqref{eqn:helmholtz} is given, for a constant $0\leq\theta_i\leq\nicefrac\pi2$, by
\[
    u(x,y) = \begin{cases}
        T \e^{i(K_1x+K_2y)} & \text{if } y>0 ,\\
        \e^{ik_1(x\cos(\theta_i)+y\sin(\theta_i))} + R \e^{ik_1(x\cos(\theta_i)-y\sin(\theta_i))} & \text{if } y<0,
    \end{cases}
\]
where $K_1=k_1\cos(\theta_i)$, $K_2 = \sqrt{k_2^2-k_1^2\sin^2(\theta_i)}$,
\[
    R=-\frac{K_2-k_1\sin(\theta_i)}{K_2+k_1\sin(\theta_i)},
\]
and $T=1+R$.
We note that there exists a critical angle $\theta_{crit}$, such that when $\theta_i>\theta_{crit}$ the wave is refracted, while $\theta_i<\theta_{crit}$ results in internal reflection, cf.~\cite[Section 6.3]{Kapita2015}. As in \cite{Kapita2015} we perform numerical experiments for both internal reflection ($\theta_i=29\degrees$) and refraction ($\theta_i=69\degrees$). To highlight the reflection and refraction behaviour, in Figures~\ref{fig:reflection:anal} and~\ref{fig:refraction:anal} we show the analytical solution when $\theta_i=29\degrees$ and $\theta_i=69\degrees$, respectively.

\begin{figure}[pt]
    \configfigure
    \subfloat[$\theta_i=29\degrees$; $h$--refinement]{\label{fig:reflection:error:h}\includegraphics[width=0.4\textwidth]{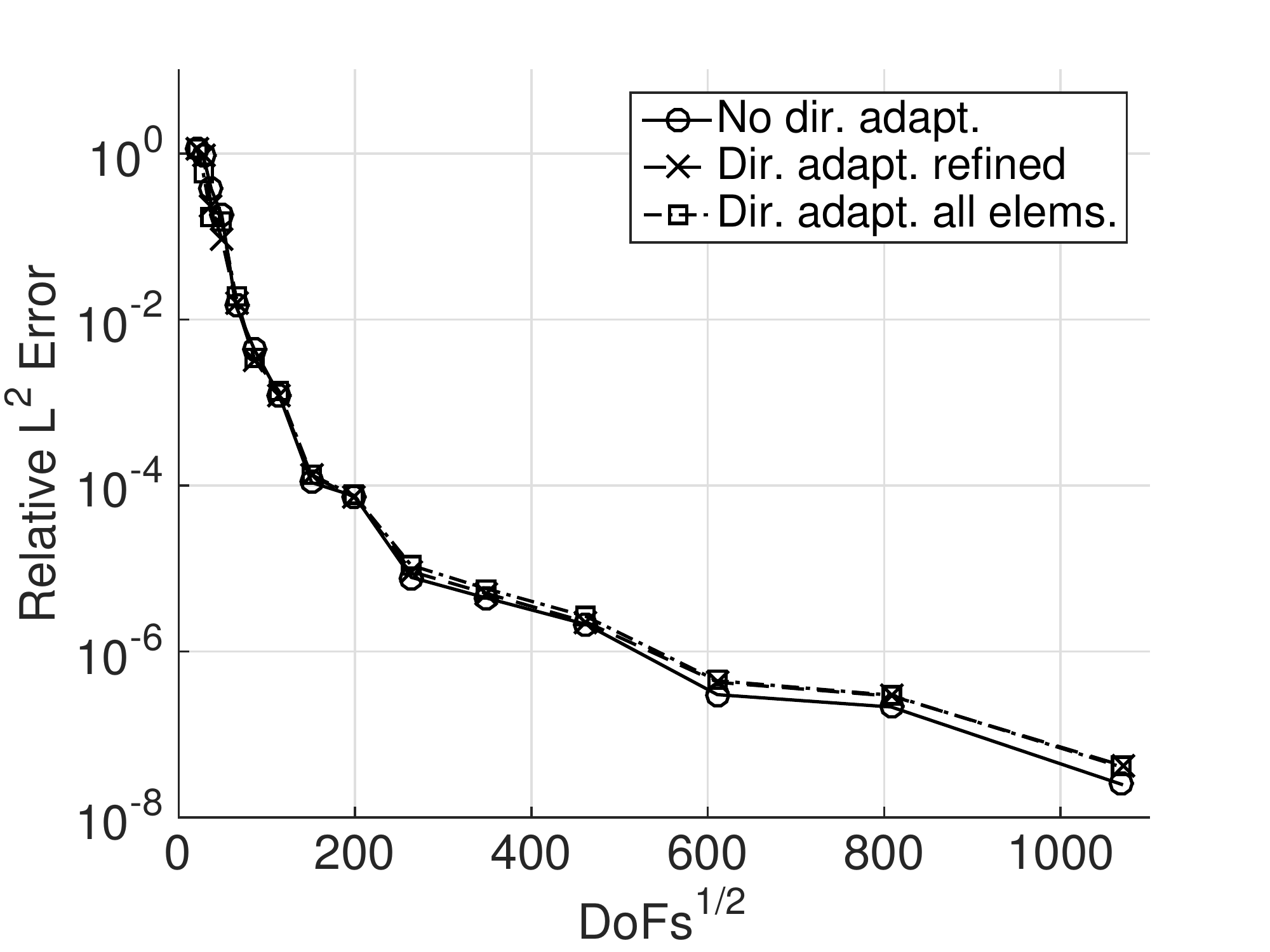}}
    \subfloat[$\theta_i=29\degrees$; $h$--refinement]{\label{fig:reflection:eff:h}\includegraphics[width=0.4\textwidth]{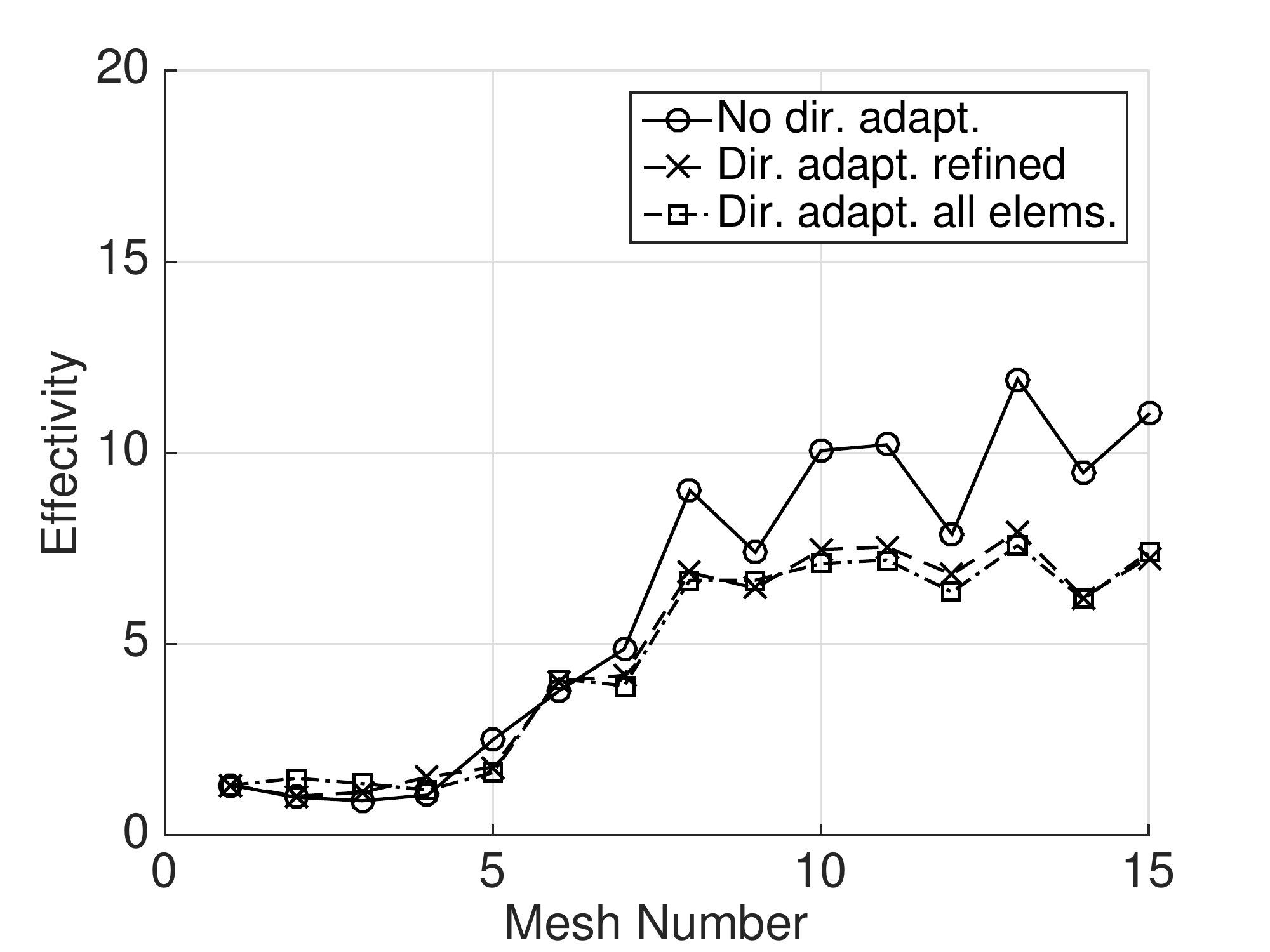}} \\
    \subfloat[$\theta_i=29\degrees$; $hp$--refinement]{\label{fig:reflection:error:hp}\includegraphics[width=0.4\textwidth]{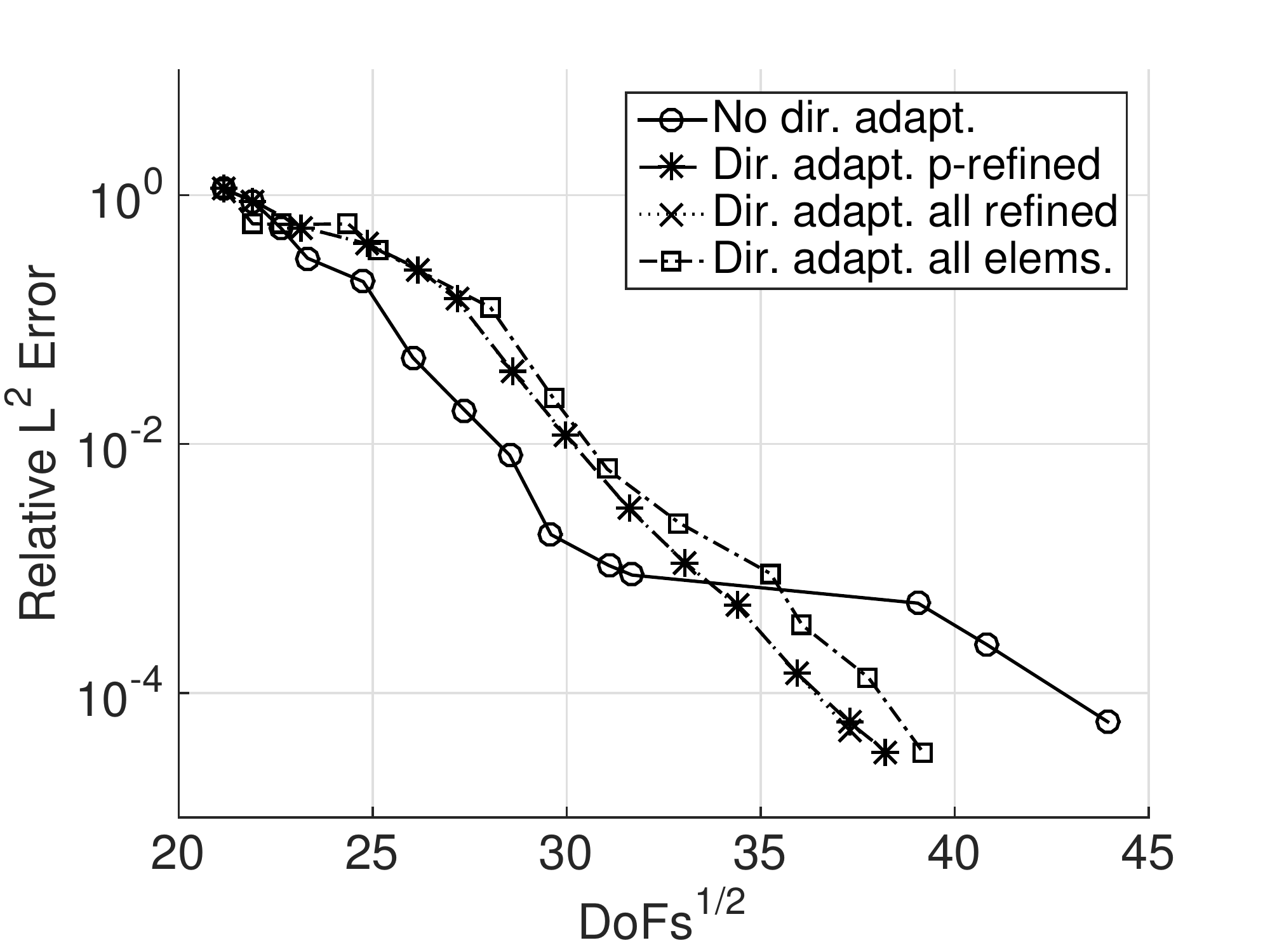}}
    \subfloat[$\theta_i=29\degrees$; $hp$--refinement]{\label{fig:reflection:eff:hp}\includegraphics[width=0.4\textwidth]{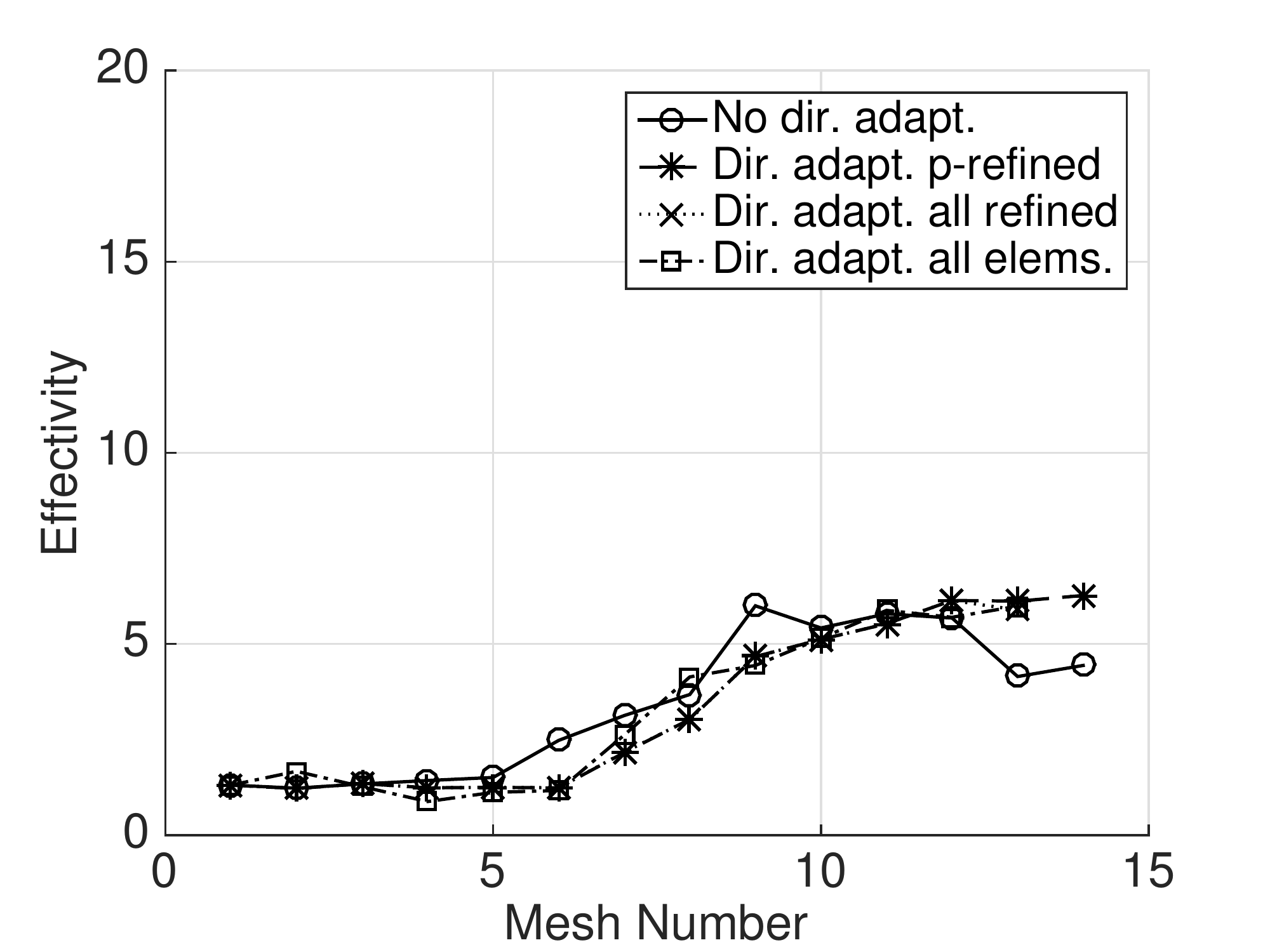}} \\
    \subfloat[$\theta_i=69\degrees$; $h$--refinement]{\label{fig:refraction:error:h}\includegraphics[width=0.4\textwidth]{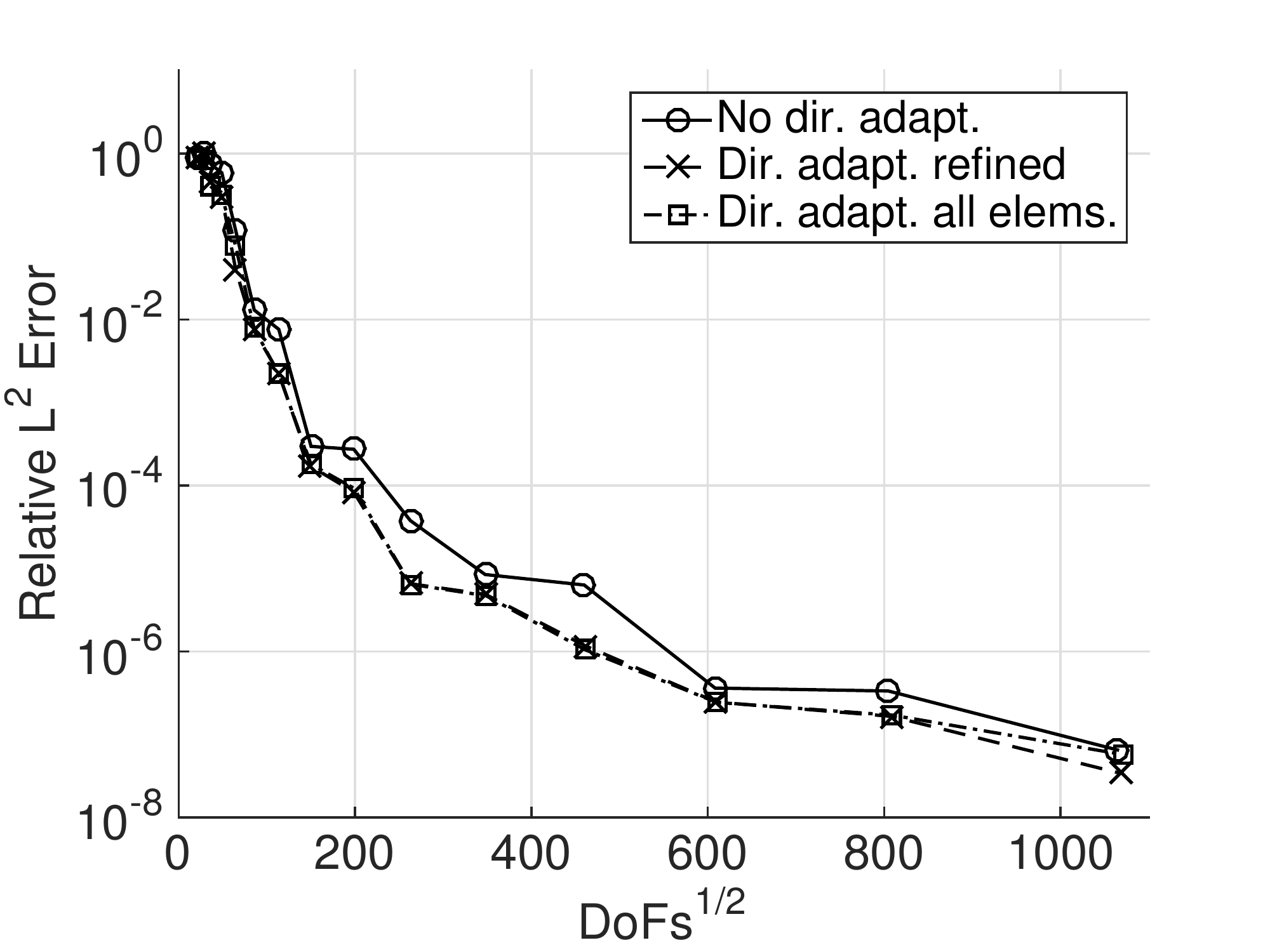}}
    \subfloat[$\theta_i=69\degrees$; $h$--refinement]{\label{fig:refraction:eff:h}\includegraphics[width=0.4\textwidth]{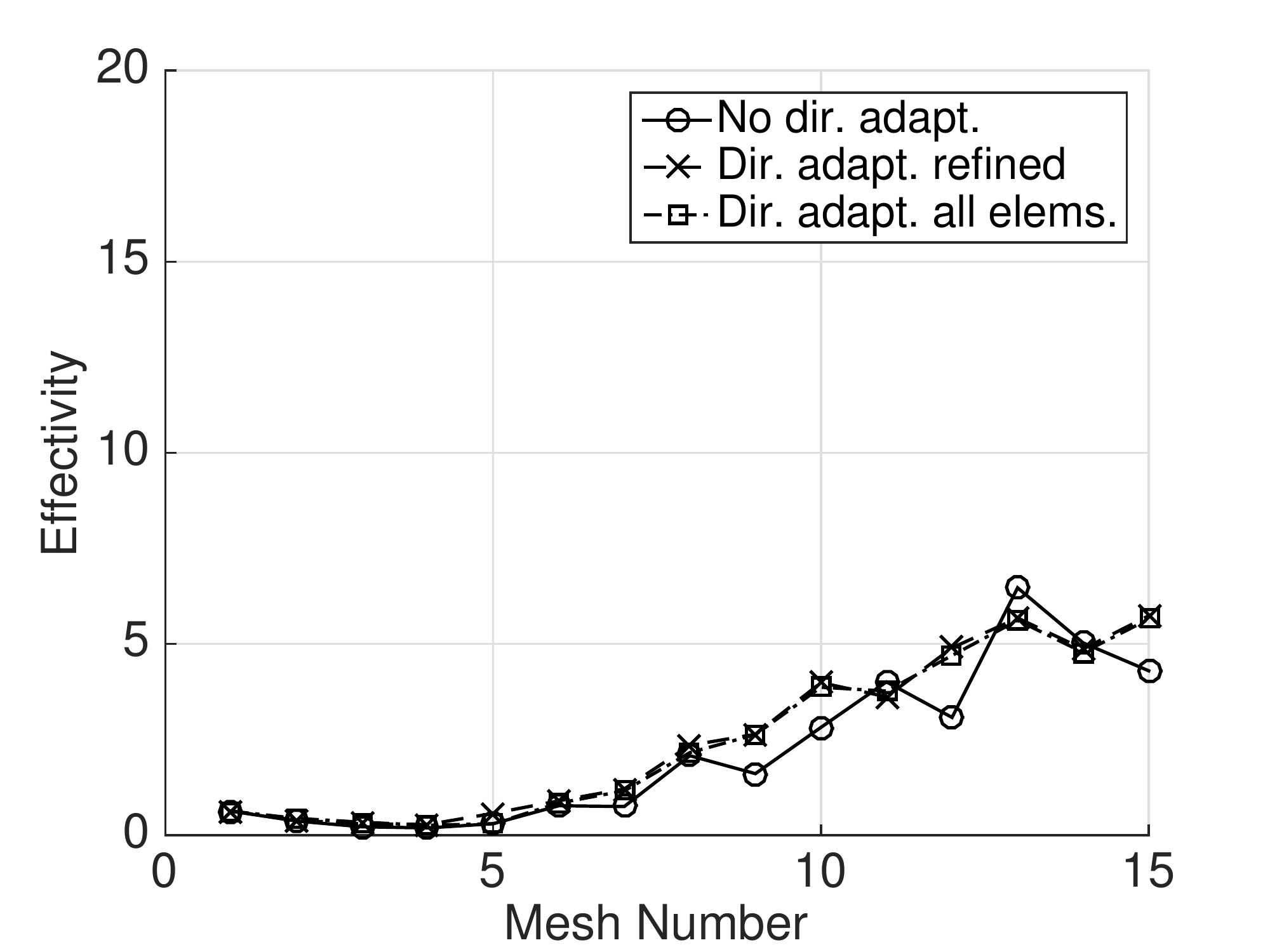}} \\
    \subfloat[$\theta_i=69\degrees$; $hp$--refinement]{\label{fig:refraction:error:hp}\includegraphics[width=0.4\textwidth]{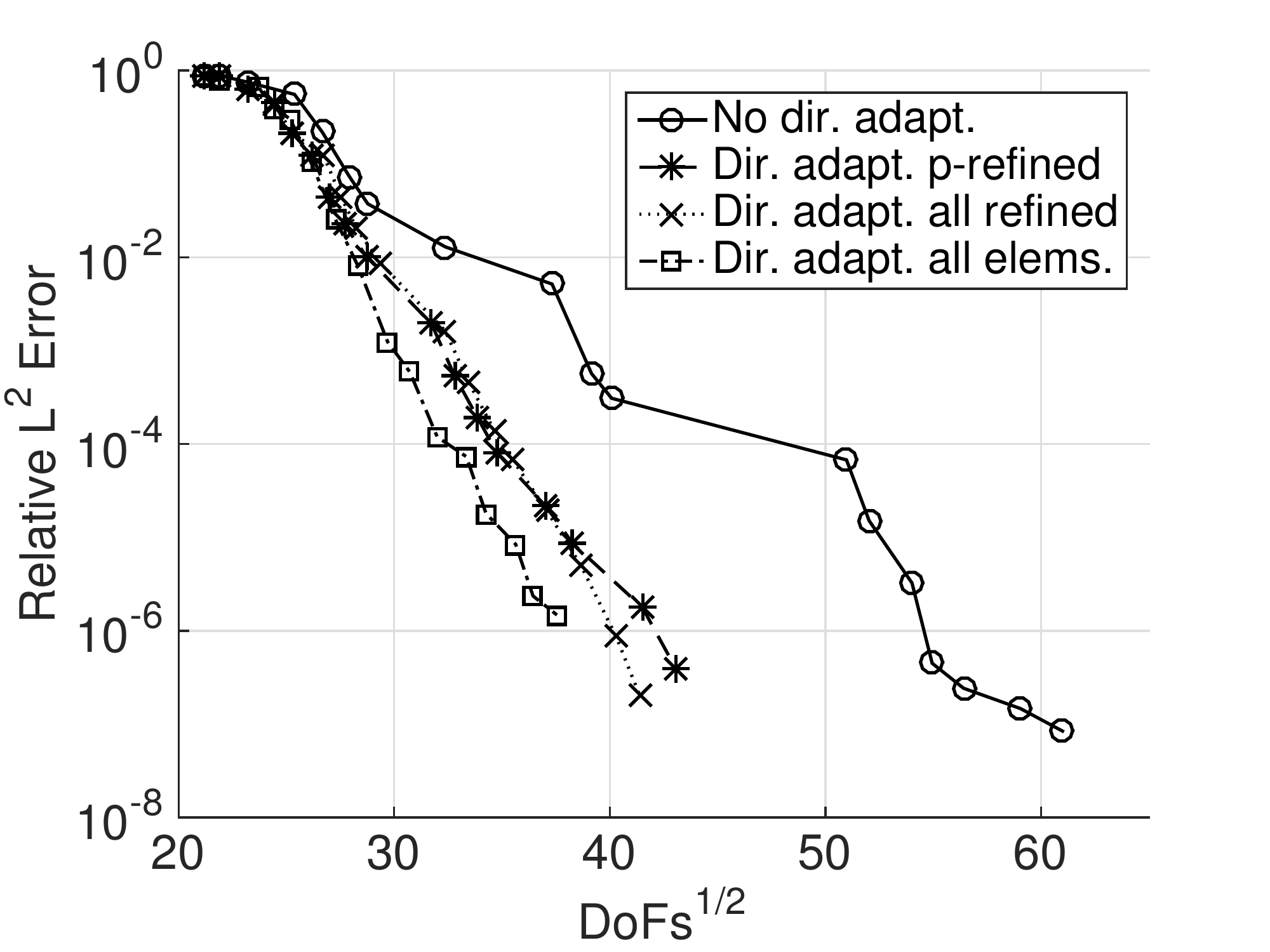}}
    \subfloat[$\theta_i=69\degrees$; $hp$--refinement]{\label{fig:refraction:eff:hp}\includegraphics[width=0.4\textwidth]{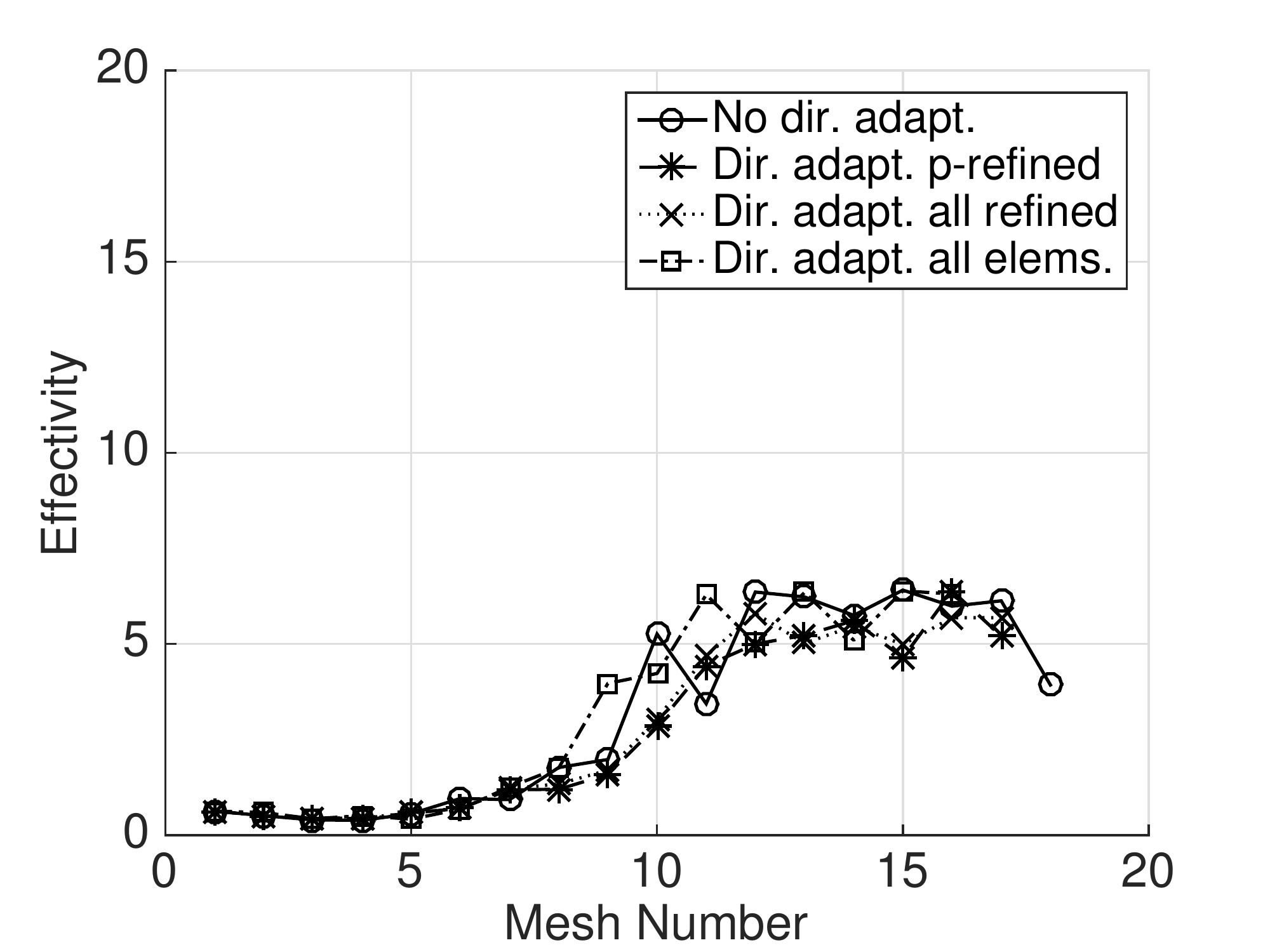}}
    \caption{Example 3: \protect\subref{fig:reflection:error:h} $L^2$-error and \protect\subref{fig:reflection:eff:h} Effectivity index for $h$--refinement with reflection ($\theta_i=29\degrees$); \protect\subref{fig:reflection:error:hp} $L^2$-error and \protect\subref{fig:reflection:eff:hp} Effectivity index for $hp$--refinement with reflection ($\theta_i=29\degrees$); \protect\subref{fig:refraction:error:h} $L^2$-error and \protect\subref{fig:refraction:eff:h} Effectivity index for $h$--refinement with refraction ($\theta_i=69\degrees$); \protect\subref{fig:refraction:error:hp} $L^2$-error and \protect\subref{fig:refraction:eff:hp} Effectivity index for $hp$--refinement with refraction ($\theta_i=69\degrees$).}
    \label{fig:reflection}
\end{figure}
\begin{figure}[pt]
    \configfigure
    \subfloat[$\theta_i=29\degrees$]{\label{fig:reflection:mesh:h}\includegraphics[width=0.4\textwidth]{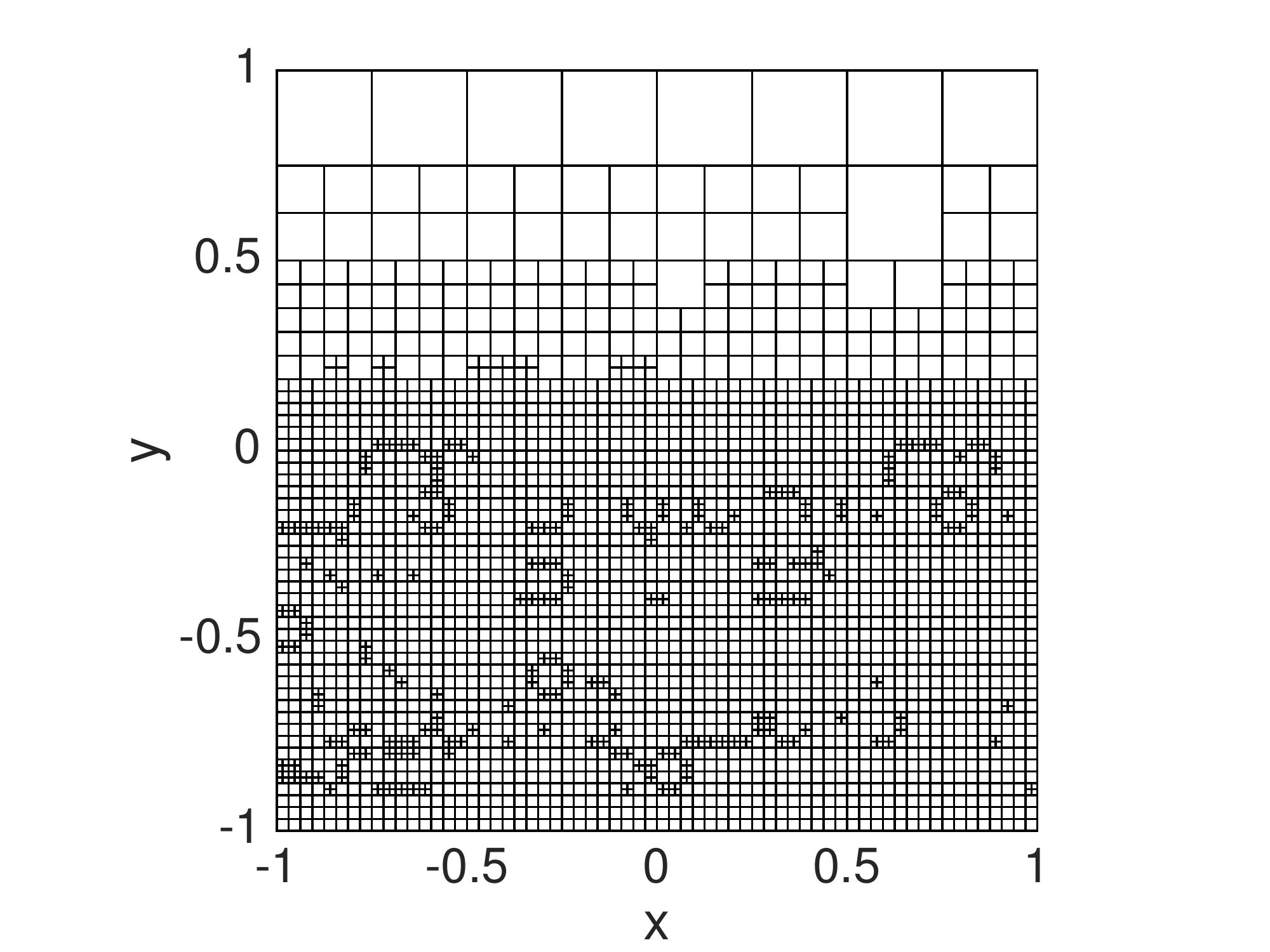}}
    \subfloat[$\theta_i=29\degrees$]{\label{fig:reflection:mesh:hp}\includegraphics[width=0.4\textwidth]{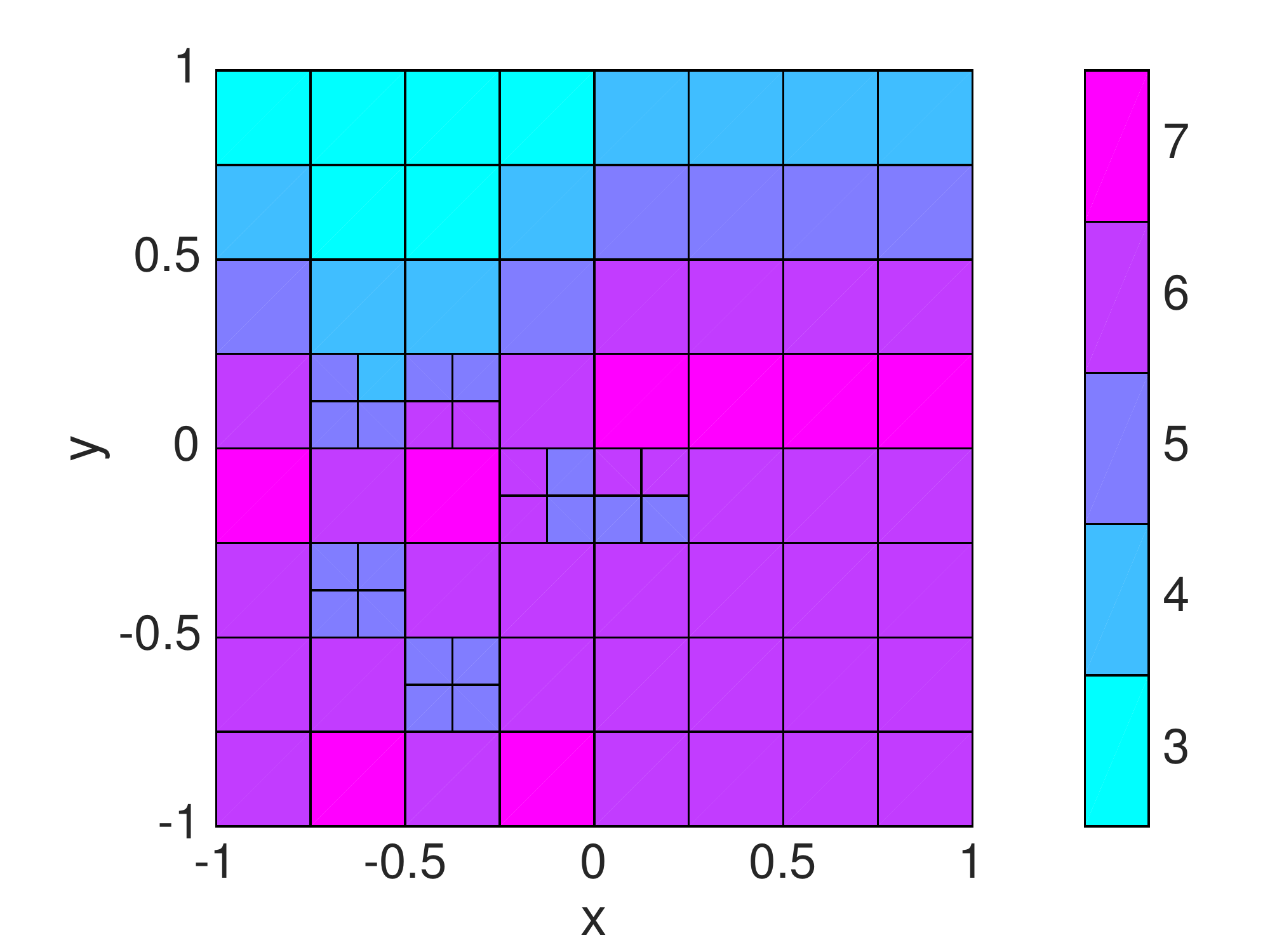}}\\
    \subfloat[$\theta_i=69\degrees$]{\label{fig:refraction:mesh:h}\includegraphics[width=0.4\textwidth]{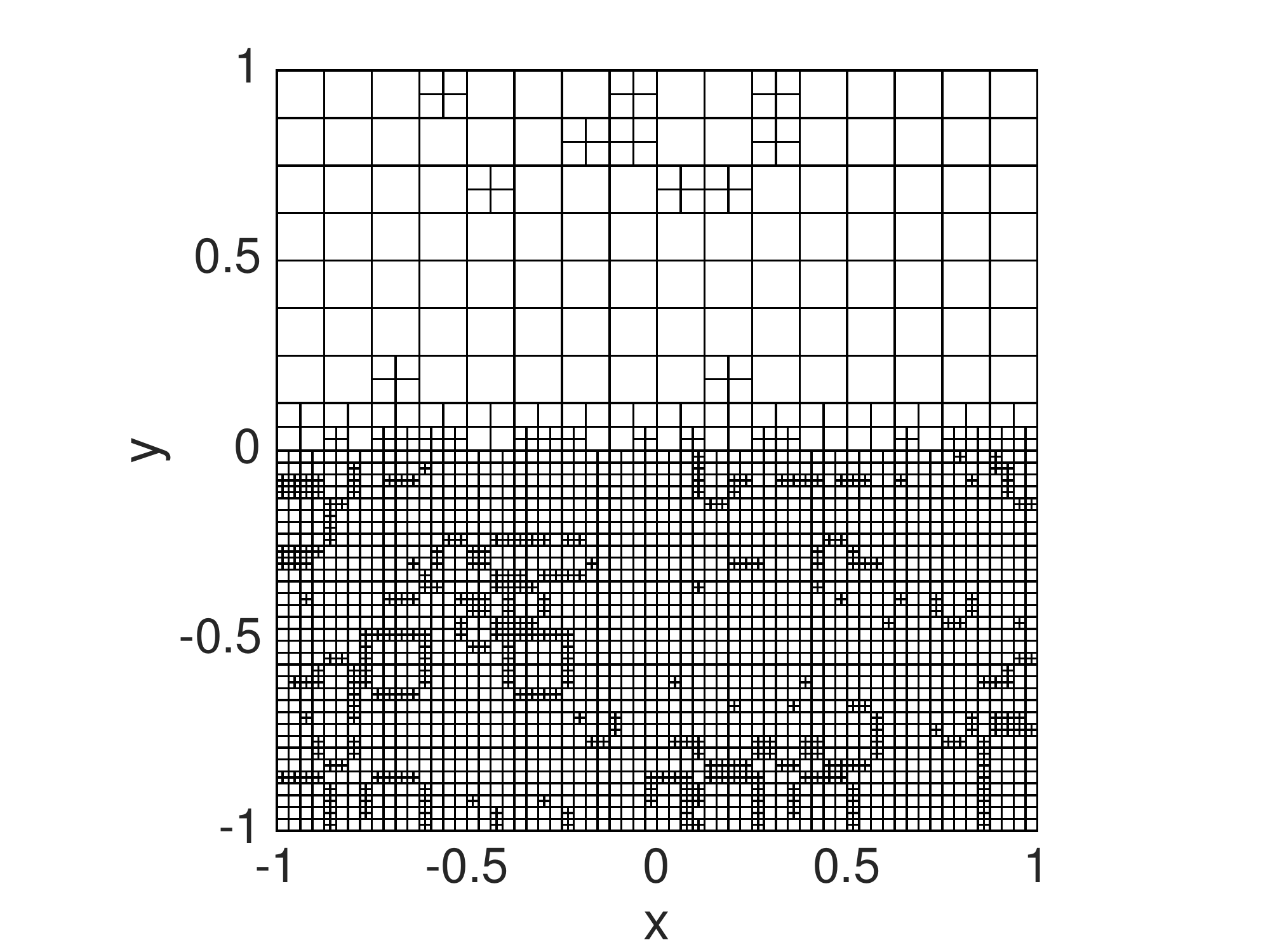}}
    \subfloat[$\theta_i=69\degrees$]{\label{fig:refraction:mesh:hp}\includegraphics[width=0.4\textwidth]{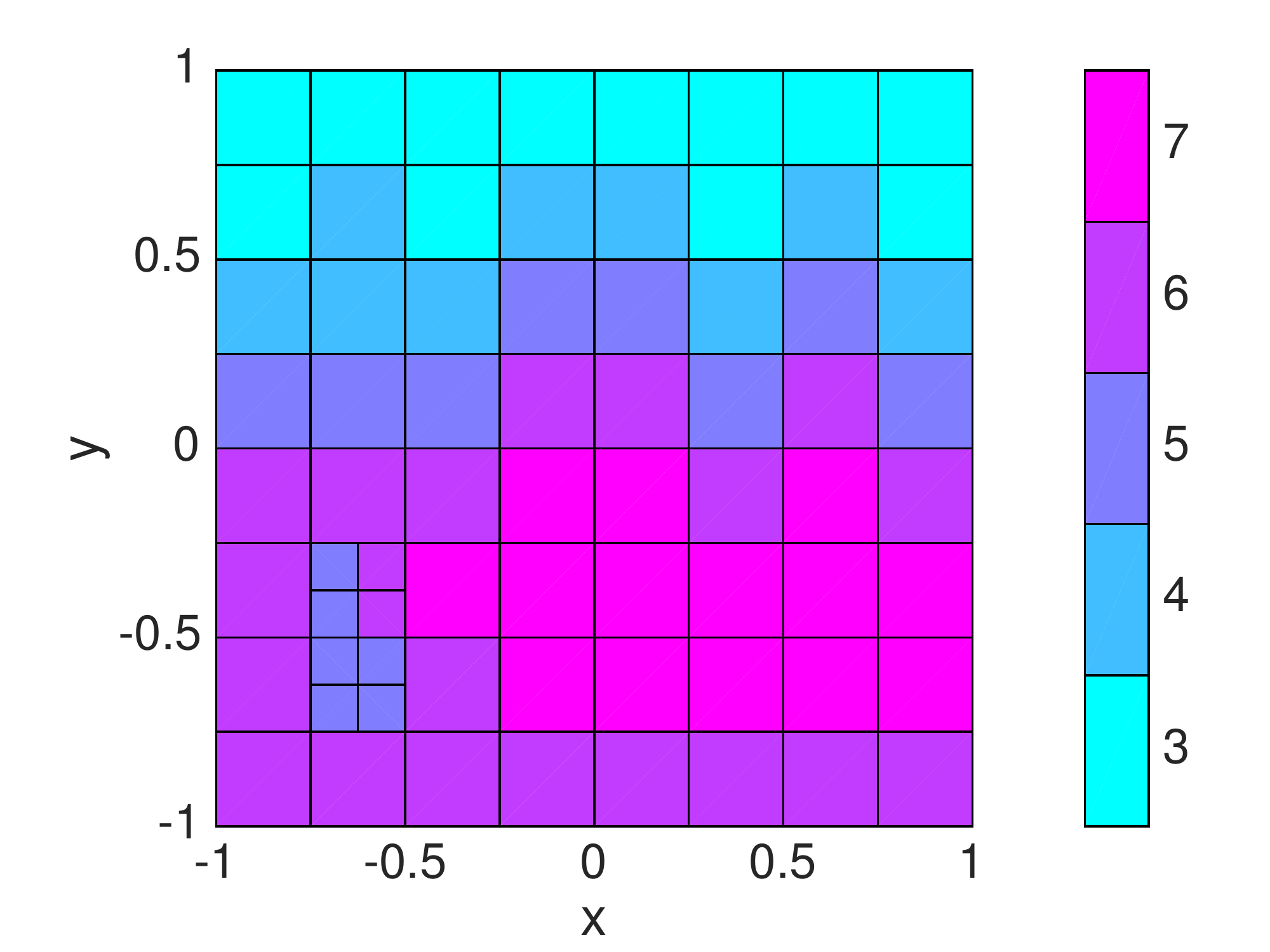}}
    \caption{Example 3: Meshes after 7 \protect\subref{fig:reflection:mesh:h} $h$-- and \protect\subref{fig:reflection:mesh:hp} $hp$--refinements for reflection ($\theta_i=29\degrees$); meshes after 7 \protect\subref{fig:refraction:mesh:h} $h$-- and \protect\subref{fig:refraction:mesh:hp} $hp$--refinements for refraction ($\theta_i=69\degrees$).}
    \label{fig:reflection:mesh}
\end{figure}

To account for the jump in the wavenumber $k$, the value of $k$ present in the integrals along the interface $y=0$ in the TDG scheme \eqref{eqn:bilinear_form} is replaced by $\omega$. We select the initial mesh to consist of $8 \times 8$ uniform square elements, so that the interface between the two materials is captured by the mesh; thereby, the wavenumber is constant in every element, and hence the TDG space~\eqref{eqn:pw_basis} and error indicators~\eqref{eqn:error_indicator} can be easily modified to treat this example by setting the wavenumber for each element equal to the wavenumber of the material within which the element is contained. Firstly, we consider the case when there is an internal reflection, i.e., when $\theta_i=29\degrees$; to this end, in Figures~\ref{fig:reflection:error:h} and \ref{fig:reflection:error:hp} we plot the relative error in the $L^2$-norm against the number of degrees of freedom in $V_{\vect{p}}(\mathcal{T}_h)$ using both $h$-- and $hp$--refinement, respectively. As for the previous numerical experiments, here we again observe exponential convergence of the error when $hp$--refinement is employed. Furthermore, in both the $h$-- and $hp$--version settings, we observe that employing directional adaptivity does not improve the magnitude of the error; indeed, in the $hp$--version setting, initially the standard refinement approach is superior, though as $V_{\vect{p}}(\mathcal{T}_h)$ is enriched, we again observe the benefits of employing directional adaptivity. This behaviour is perhaps expected, since for the internal reflection case, no waves are present above the $y=0$ line and moreover it does not possess a dominant wave propagation direction below the $y=0$ line due to the reflected waves, cf.~Figure~\ref{fig:reflection:anal}. In Figures~\ref{fig:reflection:eff:h} and \ref{fig:reflection:eff:hp}, we plot the effectivity indices for both refinement strategies, respectively; here we observe that, apart from an initial pre-asymptotic region, the effectivity indices are roughly constant.

	The corresponding convergence plots for the refraction case, i.e., when $\theta_i=69\degrees$, are presented in Figures~\ref{fig:refraction:error:h} and \ref{fig:refraction:error:hp} when both $h$-- and $hp$--refinement are employed, respectively; in the latter setting, we again observe exponential convergence of the computed relative $L^2$-norm of the error. Moreover, in contrast to the case when there is an internal reflection, here we observe the computational benefits of employing directional adaptivity, in the sense that this typically leads to a reduction in the error, for a given fixed number of degrees of freedom, when compared to the standard refinement strategy; this is particularly evidenced in the $hp$--setting. Indeed, in this case there is a dominant propagation direction throughout the domain, cf.~Figure~\ref{fig:refraction:anal}. Figures~\ref{fig:refraction:eff:h} and~\ref{fig:refraction:eff:hp} show the effectivity indices computed using both $h$-- and $hp$--refinement, respectively; analogous behaviour is observed as for the internal reflection case, i.e., the effectivity indices become roughly constant, after an initial pre-asymptotic region.

Finally, in Figures~\ref{fig:reflection:mesh:h} \&~\ref{fig:reflection:mesh:hp} we show the meshes after 7 $h$-- and $hp$-- adaptive mesh refinements have been performed, respectively, in the case of an internal reflection, i.e., $\theta_i=29\degrees$. Here, the $h$--refinement strategy concentrates most of the elements in the $y<0$ region; although, there is some refinement above $y=0$ to resolve the exponentially decaying solutions present there. Additional mesh smoothing has also been undertaken here to ensure that there is only one hanging node per face, cf.~\cite{Kapita2015}. The $hp$--refinement algorithm also performs some $h$--refinement below the $y=0$ line, though this region is largely $p$--refined; however, most of the refinement occurs around the $y=0$ line to resolve the exponentially decaying solutions. Some $p$--refinement occurs in the rest of the $y>0$ region, which is caused by enforcing the condition that the effective polynomial degree may only vary by one between neighboring elements. In the refraction case, i.e., $\theta_i=69\degrees$, cf. Figures~\ref{fig:refraction:mesh:h} \&~\ref{fig:refraction:mesh:hp}, we note a sharp boundary at the $y=0$, with more refinement undertaken in the $y<0$ region than the region $y>0$.

\subsubsection{Example 4 --- 3D smooth solution (plane wave)}
\begin{figure}[pt]
    \configfigure
    \subfloat[$k=20$; $h$--refinement]{\label{fig:cube:error:20h}\includegraphics[width=0.4\textwidth]{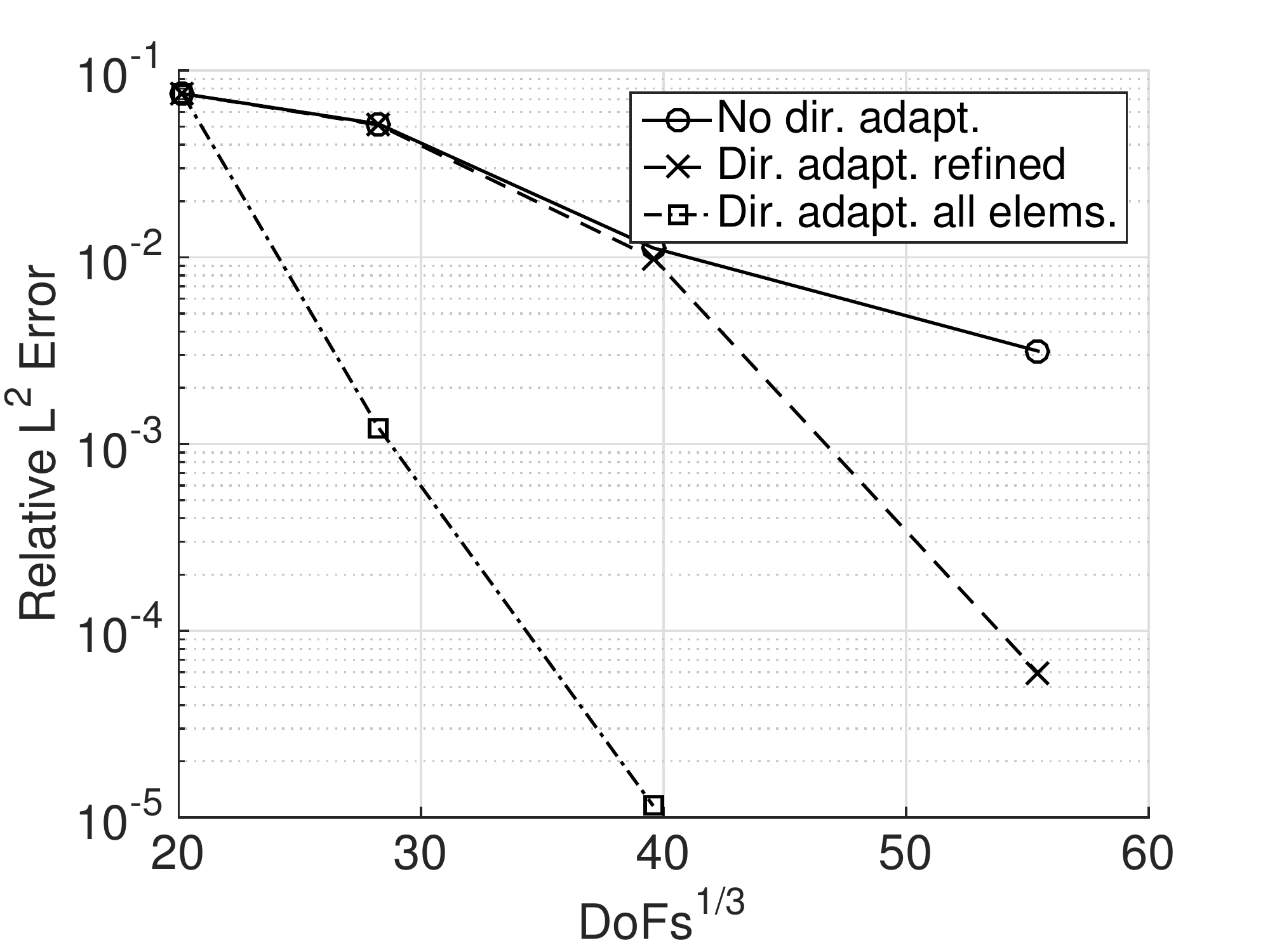}}
    \subfloat[$k=20$; $h$--refinement]{\label{fig:cube:eff:20h}\includegraphics[width=0.4\textwidth]{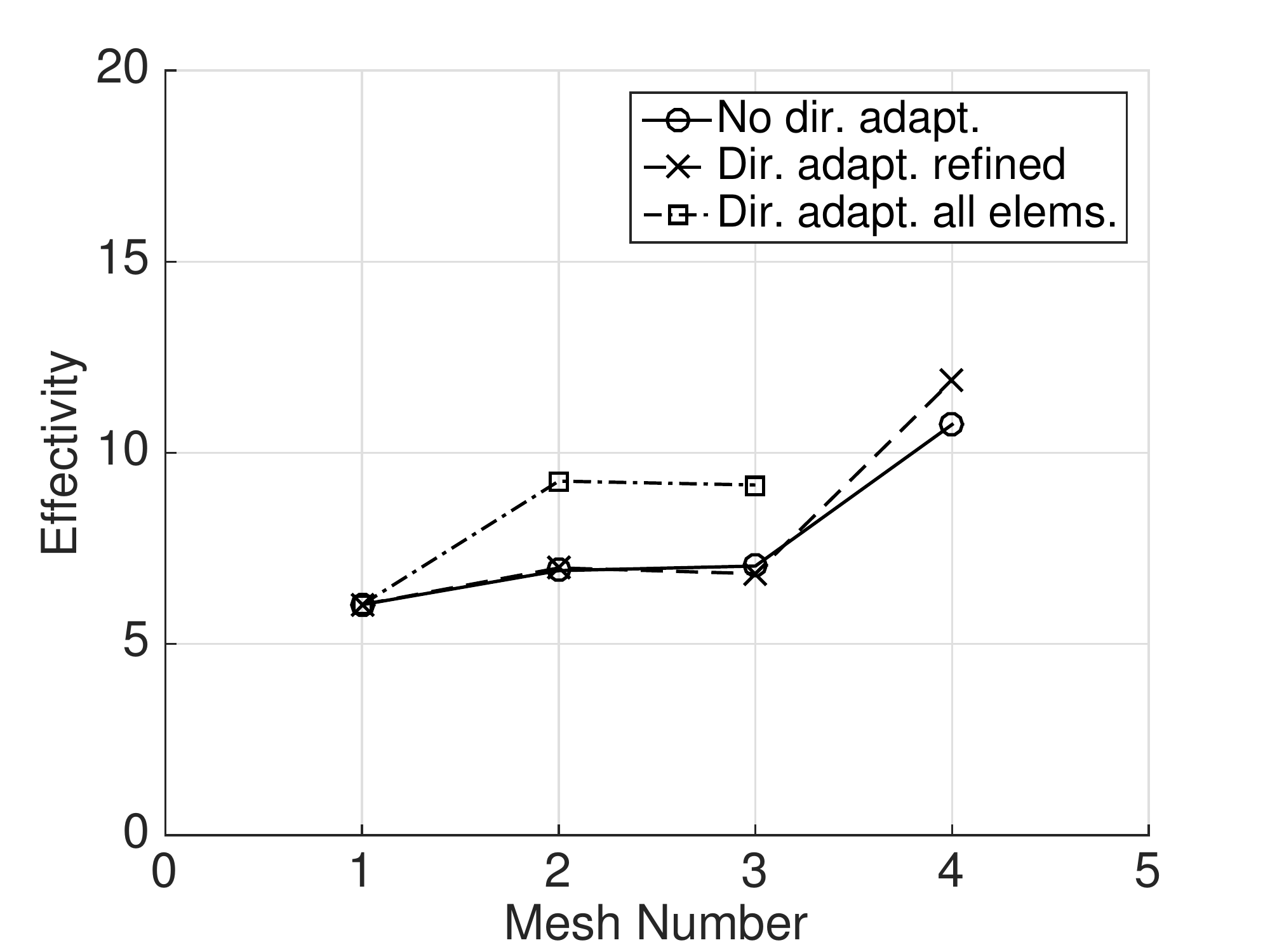}} \\
    \subfloat[$k=20$; $hp$--refinement]{\label{fig:cube:error:20hp}\includegraphics[width=0.4\textwidth]{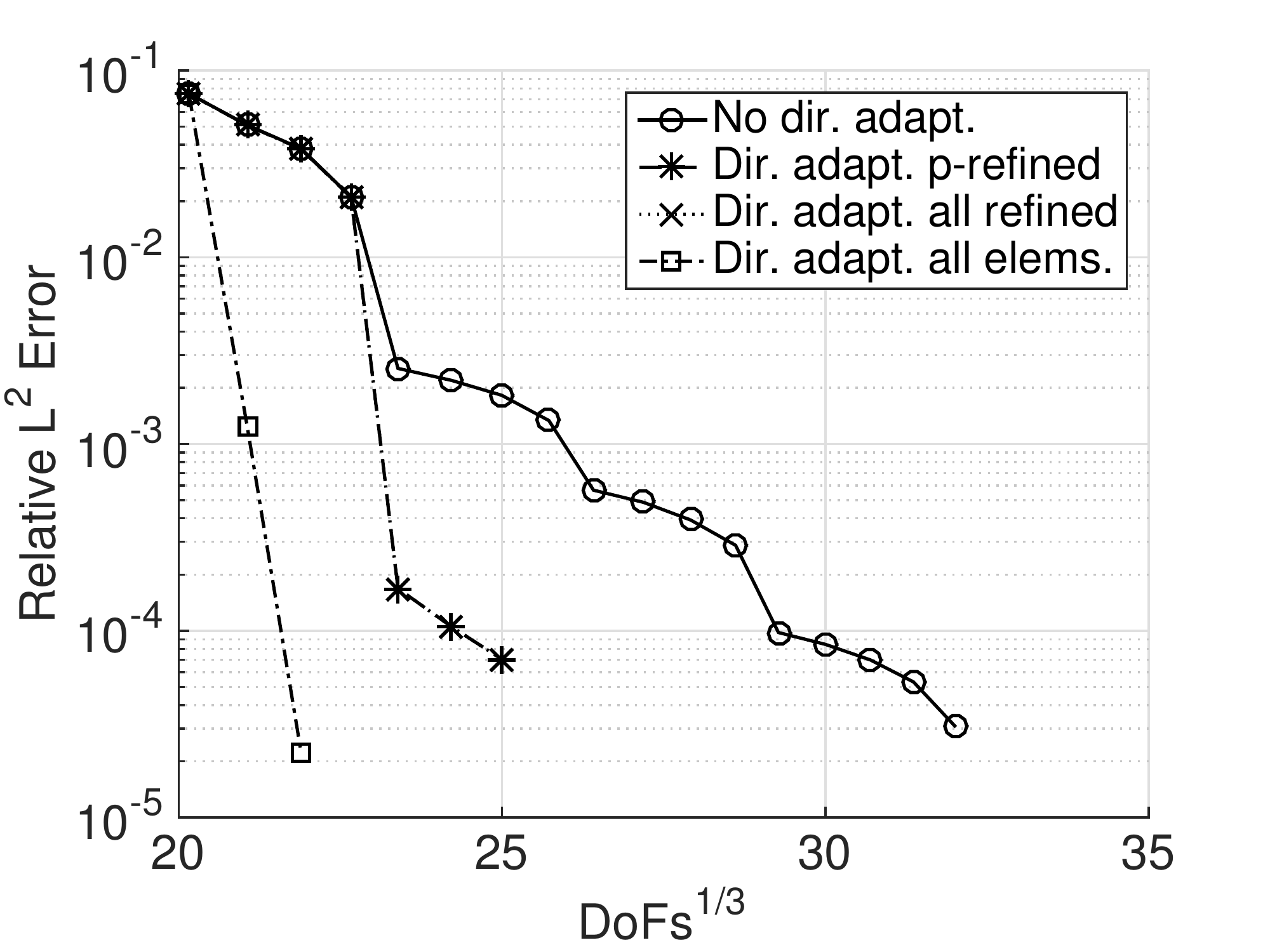}}
    \subfloat[$k=20$; $hp$--refinement]{\label{fig:cube:eff:20hp}\includegraphics[width=0.4\textwidth]{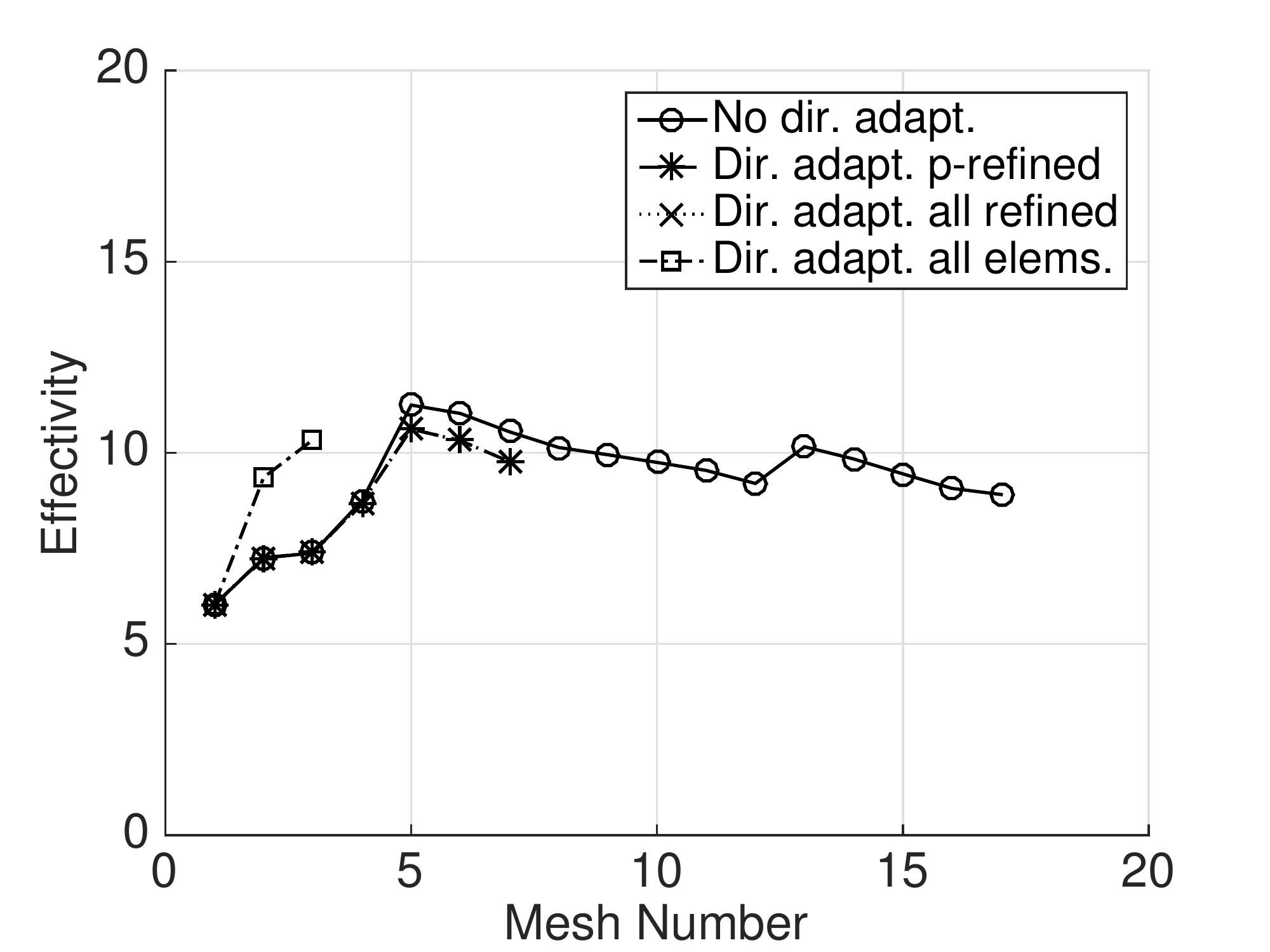}} \\
    \subfloat[$k=50$; $h$--refinement]{\label{fig:cube:error:50h}\includegraphics[width=0.4\textwidth]{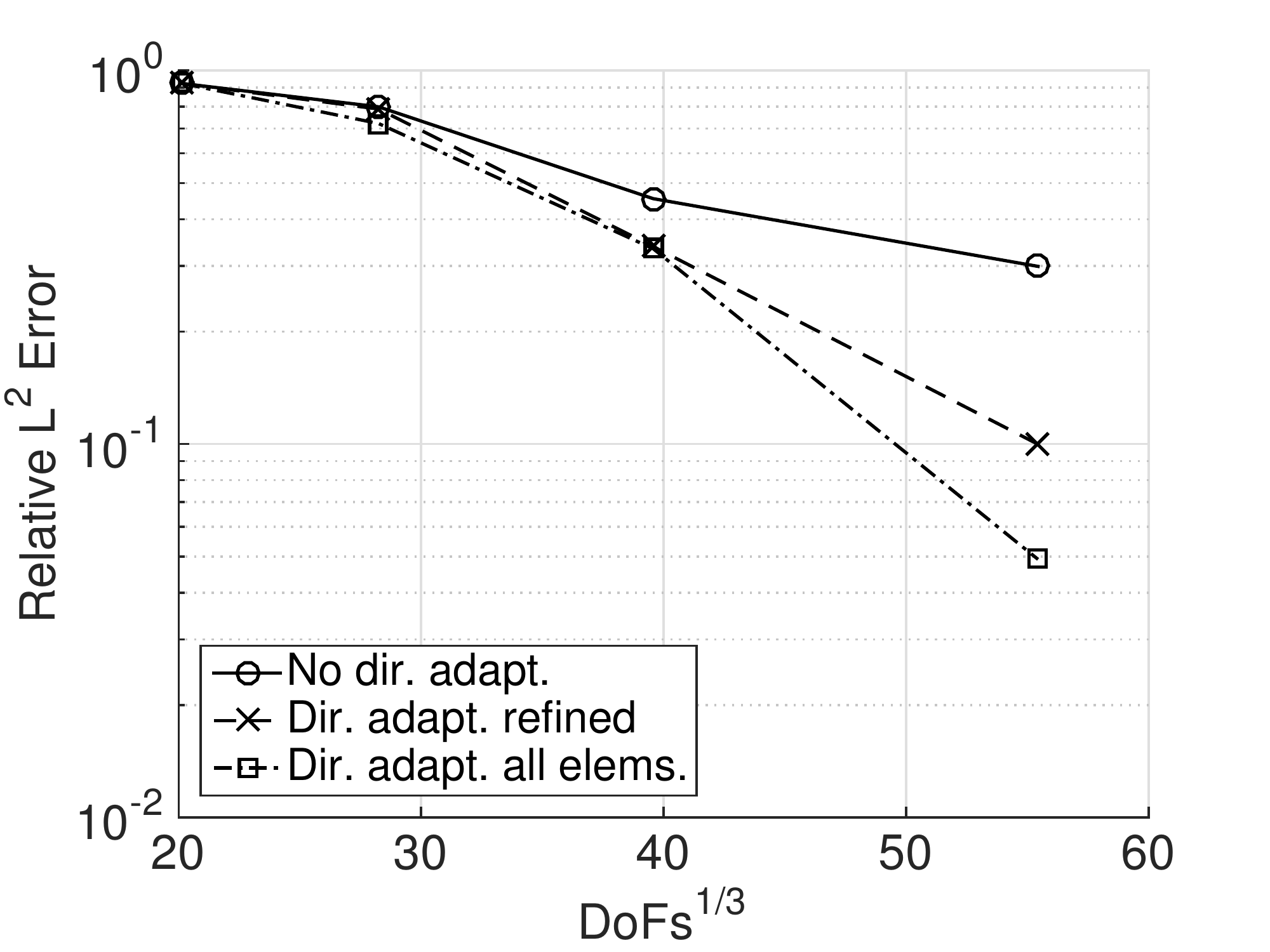}}
    \subfloat[$k=50$; $h$--refinement]{\label{fig:cube:eff:50h}\includegraphics[width=0.4\textwidth]{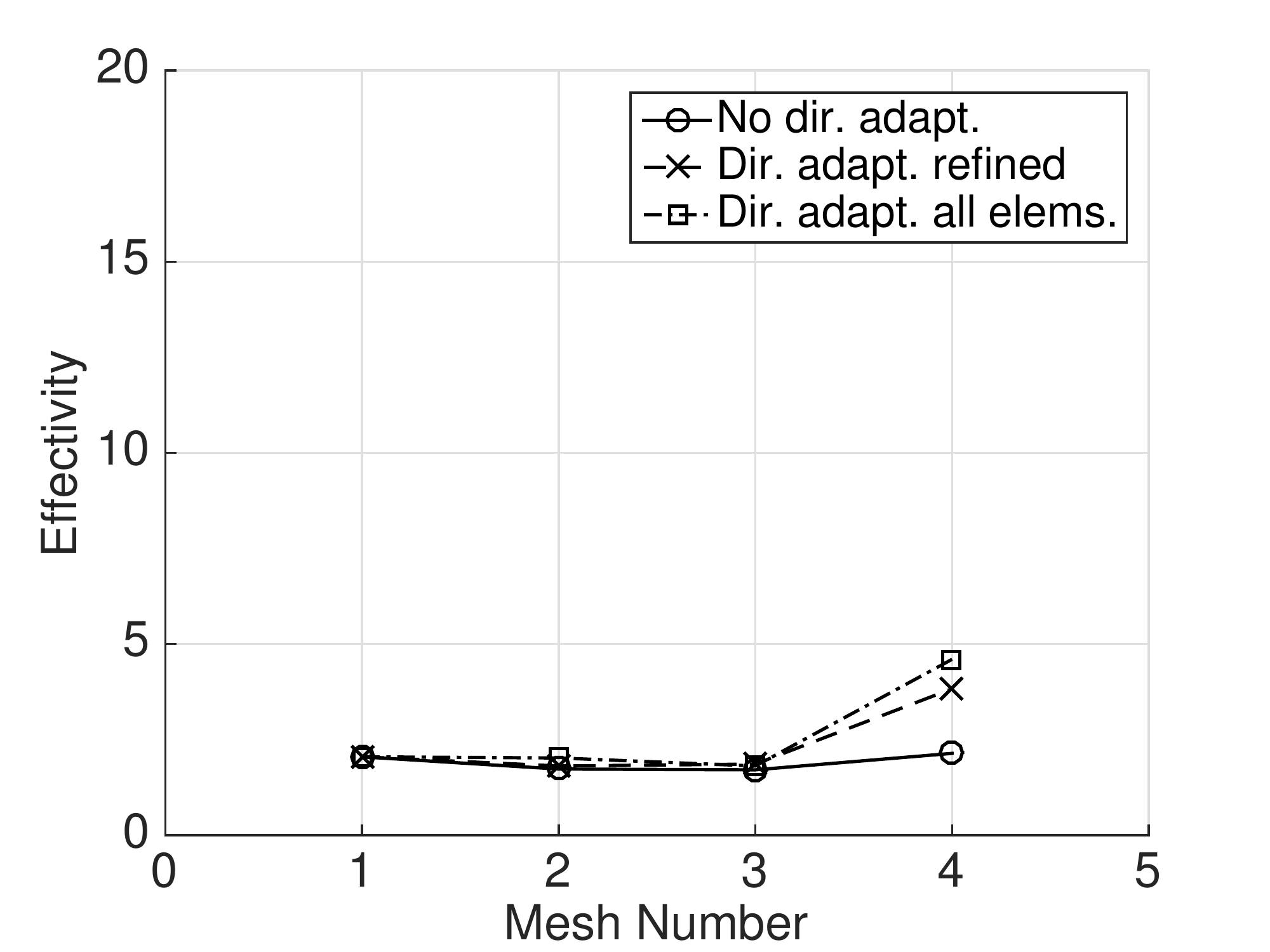}} \\
    \subfloat[$k=50$; $hp$--refinement]{\label{fig:cube:error:50hp}\includegraphics[width=0.4\textwidth]{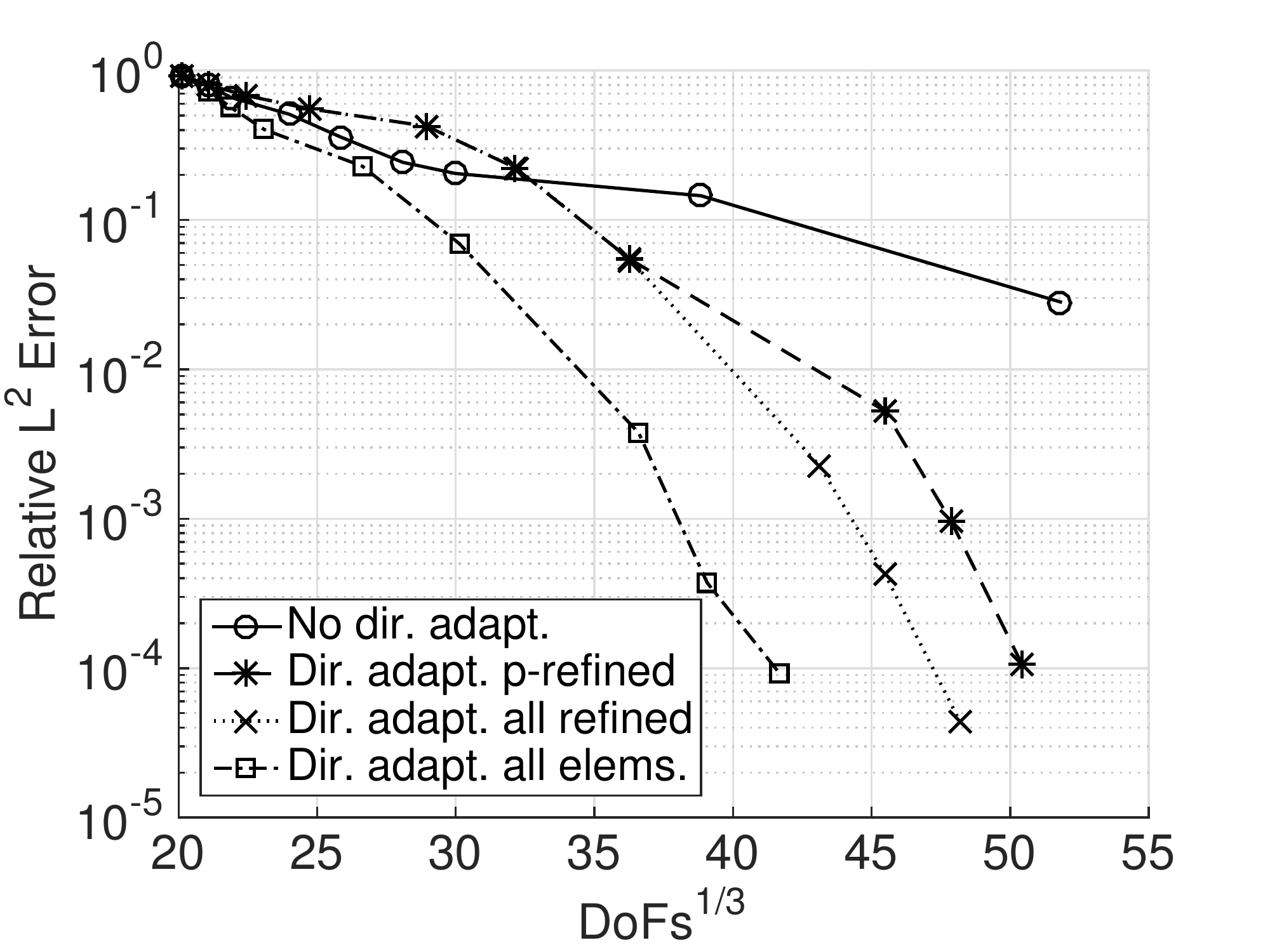}}
    \subfloat[$k=50$; $hp$--refinement]{\label{fig:cube:eff:50hp}\includegraphics[width=0.4\textwidth]{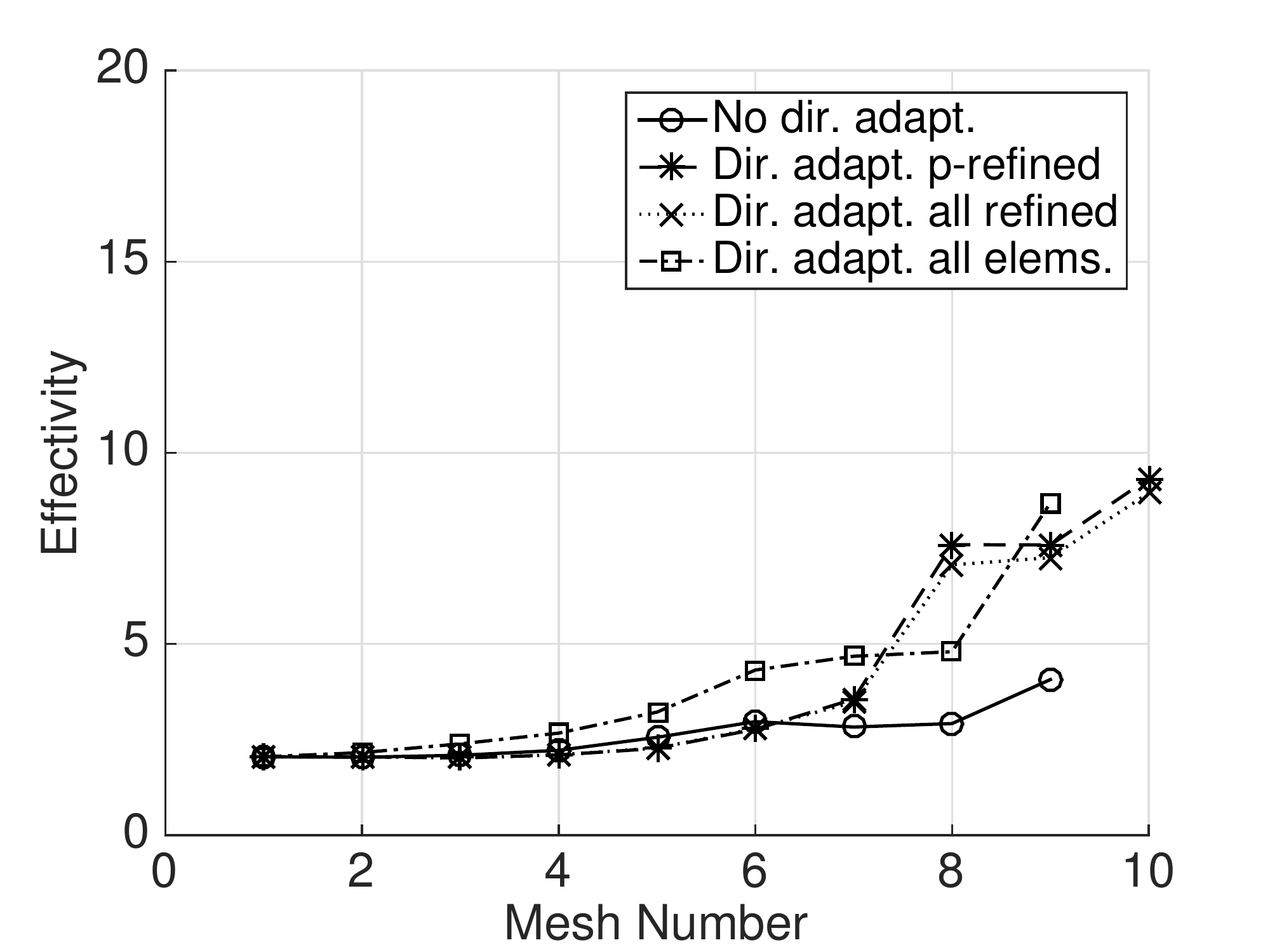}}
    \caption{Example 4: \protect\subref{fig:cube:error:20h} $L^2$-error and \protect\subref{fig:cube:eff:20h} Effectivity index for $h$--refinement with wavenumber $k=20$; \protect\subref{fig:cube:error:20hp} $L^2$-error and \protect\subref{fig:cube:eff:20hp} Effectivity index for $hp$--refinement with $k=20$; \protect\subref{fig:cube:error:50h} $L^2$-error and \protect\subref{fig:cube:eff:50h} Effectivity index for $h$--refinement with $k=50$; \protect\subref{fig:cube:error:50hp} $L^2$-error and \protect\subref{fig:cube:eff:50hp} Effectivity index for $hp$--refinement with $k=50$.}
    \label{fig:cube}
\end{figure}

In this final example, we consider problem \eqref{eqn:helmholtz} posed on the domain $\Omega=(0,1)^3$, $\Gamma_R=\partial\Omega$, and $\Gamma_D\equiv\emptyset$, with Robin boundary condition $g_R$ selected so that the analytical solution $u$ to \eqref{eqn:helmholtz} is given by
\[
    u(\vect{x}) = \e^{ik\vect{\vect{d}\cdot\vect{x}}},
\]
where $\vect{d}_j = \nicefrac{1}{\sqrt{3}}$ for $j=1,2,3$.

In Figures~\ref{fig:cube:error:20h} and \ref{fig:cube:error:50h} we present the performance of the proposed directional adaptivity algorithm employing $h$--refinement with wavenumbers $k=20$ and $k=50$, respectively; the analogous results for $hp$--refinement are given in Figures~\ref{fig:cube:error:20hp} and \ref{fig:cube:error:50hp}, respectively. As in the two--dimensional setting, we observe that selecting more elements for directional adaptivity at each step of the proposed refinement strategy, leads to a greater reduction in the relative $L^2$-norm of the error, for a fixed number of degrees of freedom, when compared to the standard case when directional adaptivity is not employed. Of course, given the simple nature of the analytical solution for this problem, we clearly expect directional adaptivity to be advantageous. In the case when the wavenumber $k=50$ we note that both $h$-- and $hp$--refinement strategies are essentially in the pre-asymptotic region; however, performing directional adaptivity ensures that the method leaves this pre-asymptotic region after only a few mesh refinements. Finally, in Figures~\ref{fig:cube:eff:20h}, \ref{fig:cube:eff:20hp}, \ref{fig:cube:eff:50h}, and \ref{fig:cube:eff:50hp} we plot the effectivity indices of both the $h$-- and $hp$--refinement algorithms for the case when $k=20,50$. We note, especially in the $hp$--refinement case, that the effectivity indices are roughly constant but do slightly rise after the pre-asymptotic region. 

\section{Concluding remarks} \label{sec:conclusions}

In this article we have developed an automatic $hp$--adaptive refinement algorithm
for the TDG approximation of the homogeneous Helmholtz equation. In addition
to employing both local mesh subdivision and local basis enrichment,
we also locally rotate the underlying plane wave basis in such a manner
so that the first basis function is aligned with the dominant wave
direction. The choice to $h$-- or $p$--refine an element is based on
a prediction of how much reduction we expect to observe in the elementwise
error indicator, when a particular refinement is performed. The alignment
of the local basis with the dominant wave direction is undertaken on the basis of an
eigenvalue analysis of the Hessian of the numerical solution, together with
a correction computed from an impedance condition. The computational
efficiency of the proposed adaptive strategy has been studied through a 
series of numerical examples; indeed, the application of $hp$--refinement,
with directional adaptivity, leads to a significant reduction in the
computed error compared to standard refinement strategies. We also note that 
performing directional adaptivity on all elements generally leads to a greater reduction
in the error than the corresponding case when only elements marked for
refinement are directionally adapted; clearly, this error reduction is attained while 
keeping the number of degrees of freedom in the underlying TDG space fixed. Future work
will be devoted to the derivation of robust $hp$--version {\em a posteriori}
error bounds, as well as the application to problems of engineering interest.

\section*{Acknowledgement}
S.~Congreve and I.~Perugia have been funded by the Austrian Science Fund (FWF) through the project P29197-N32. I.~Perugia has also been funded by the FWF through the project F~65.

\bibliographystyle{abbrv}
\bibliography{references_techrep}

\begin{thebibliography}{10}

\bibitem{Amara2014}
M.~Amara, S.~Chaudhry, J.~Diaz, R.~Djellouli, and S.~L. Fiedler.
\newblock A local wave tracking strategy for efficiently solving mid- and
  high-frequency {H}elmholtz problems.
\newblock {\em Comput. Methods Appl. Mech. and Engrg.}, 276:473--508, 2014.

\bibitem{Amara2009}
M.~Amara, R.~Djellouli, and C.~Farhat.
\newblock Convergence analysis of a discontinuous {G}alerkin method with plane
  waves and {L}agrange multipliers for the solution of {H}elmholtz problems.
\newblock {\em {SIAM} J. Numer. Anal.}, 47(2):1038--1066, 2009.

\bibitem{BabI}
F.~Babu\v{s}ka, F.~Ihlenburg, T.~Strouboulis, and S.~K. Gangaraj.
\newblock A posteriori error estimation for finite element solutions of
  {H}elmholtz' equation {I}. {T}he quality of local indicators and estimators.
\newblock {\em Internat. J. Numer. Methods Engrg.}, 40(18):3443--3462, 1997.

\bibitem{BabII}
F.~Babu\v{s}ka, F.~Ihlenburg, T.~Strouboulis, and S.~K. Gangaraj.
\newblock A posteriori error estimation for finite element solutions of
  {H}elmholtz' equation {II}. {E}stimation of the pollution error.
\newblock {\em Internat. J. Numer. Methods Engrg.}, 40(21):3883--3900, 1997.

\bibitem{Betcke2011}
T.~Betcke and J.~Phillips.
\newblock Adaptive plane wave discontinuous {G}alerkin methods for {H}elmholtz
  problems.
\newblock In {\em Proceedings of the 10th International Conference on the
  Mathematical and Numerical Aspects of Waves}, pages 261--264, 2011.

\bibitem{Betcke2012}
T.~Betcke and J.~Phillips.
\newblock Approximation by dominant wave directions in plane wave methods.
\newblock Technical report, UCL, 2012.
\newblock Available at \url{http://discovery.ucl.ac.uk/1342769/}.

\bibitem{Braess}
D.~Braess.
\newblock {\em Finite elements. Theory, fast solvers, and applications in solid
  mechanics}.
\newblock Cambridge University Press, Cambridge, second edition, 2001.

\bibitem{Cessenat1998}
O.~Cessenat and B.~Despr\'es.
\newblock Application of an {U}ltra {W}eak {V}ariational {F}ormulation of
  elliptic {PDE}s to the two-dimensional {H}elmholtz problem.
\newblock {\em SIAM J. Numer. Anal.}, 35(1):255--299, 1998.

\bibitem{SauterDoerfler}
W.~D\"orfler and S.~Sauter.
\newblock A posteriori error estimation for highly indefinite {H}elmholtz
  problems.
\newblock {\em Comp. Meth. Appl. Math.}, 13(3):333--347, 2013.

\bibitem{Formaggia2001}
L.~Formaggia and S.~Perotto.
\newblock New anisotropic a priori error estimates.
\newblock {\em Numer. Math.}, 89(4):641--667, 2001.

\bibitem{Formaggia2003}
L.~Formaggia and S.~Perotto.
\newblock Anisotropic error estimates for elliptic problems.
\newblock {\em Numer. Math.}, 94(1):67--92, 2003.

\bibitem{ghh-paper}
E.~Georgoulis, E.~Hall, and P.~Houston.
\newblock Discontinuous {G}alerkin methods for advection--diffusion--reaction
  problems on anisotropically refined meshes.
\newblock {\em {SIAM} J. Sci. Comput.}, 30(1):246--271, 2007.

\bibitem{Gittelson2008}
C.~J. Gittelson.
\newblock Plane wave discontinuous {G}alerkin methods.
\newblock Master's thesis, ETH Zurich, 2008.
\newblock
  \url{http://www.sam.math.ethz.ch/~hiptmair/StudentProjects/Gittelson/thesis.pdf}.

\bibitem{Gittelson2009}
C.~J. Gittelson, R.~Hiptmair, and I.~Perugia.
\newblock Plane wave discontinuous {G}alerkin methods: Analysis of the
  $h$-version.
\newblock {\em ESAIM Math. Model. Numer. Anal.}, 43(2):297--331, 2009.

\bibitem{hall-thesis}
E.~J.~C. Hall.
\newblock {\em Anisotropic Adaptive Refinement For Discontinuous {G}alerkin
  Methods}.
\newblock PhD thesis, Department of Mathematics, University of Leicester, 2007.

\bibitem{Hiptmair2014}
R.~Hiptmair, A.~Moiola, and I.~Perugia.
\newblock Trefftz discontinuous {G}alerkin methods for acoustic scattering on
  locally refined meshes.
\newblock {\em Appl. Numer. Math.}, 79:79--91, 2014.

\bibitem{Hiptmair2016}
R.~Hiptmair, A.~Moiola, and I.~Perugia.
\newblock A survey of {T}refftz methods for the {H}elmholtz equation.
\newblock In G.~R. Barrenechea, F.~Brezzi, A.~Cangiani, and E.~H. Georgoulis,
  editors, {\em Building Bridges: Connections and Challenges in Modern
  Approaches to Numerical Partial Differential Equations}, pages 237--279.
  Springer, Cham, 2016.

\bibitem{aptofem}
P.~Houston.
\newblock {AptoFEM} finite element analysis software, 2017.
\newblock \url{http://www.aptofem.com} [Online].

\bibitem{Houston2005}
P.~Houston and E.~S{\"u}li.
\newblock A note on the design of $hp$-adaptive finite element methods for
  elliptic partial differential equations.
\newblock {\em Comput. Methods Appl. Mech. Engrg.}, 194(2--5):229--243, 2005.

\bibitem{Huttunen2002}
T.~Huttunen, P.~Monk, and J.~P. Kaipio.
\newblock Computational aspects of the ultra-weak variational formulation.
\newblock {\em J. Comput. Phys.}, 182(1):27--46, 2002.

\bibitem{Kapita2015}
S.~Kapita, P.~Monk, and T.~Warburton.
\newblock Residual-based adaptivity and {PWDG} methods for the {H}elmholtz
  equation.
\newblock {\em {SIAM} J. Sci. Comput.}, 37(3), 2015.

\bibitem{Luostari2013}
T.~Luostari, T.~Huttunen, and P.~Monk.
\newblock Improvements for the ultra weak variational formulation.
\newblock {\em Internat. J. Numer. Methods Engrg.}, 94(6):598--624, 2013.

\bibitem{Melenk2001}
J.~M. Melenk and B.~I. Wohlmuth.
\newblock On residual-based a posteriori error estimation in $hp$-{FEM}.
\newblock {\em Adv. Comp. Math.}, 15(1--4):311--331, 2001.

\bibitem{Mitchell2011}
W.~F. Mitchell and M.~A. McClain.
\newblock A comparison of $hp$-adaptive strategies for elliptic partial
  differential equations.
\newblock Technical Report NISTIR 7824, National Institute of Standards and
  Technology, 2011.

\bibitem{Mitchell2014}
W.~F. Mitchell and M.~A. McClain.
\newblock A comparison of hp-adaptive strategies for elliptic partial
  differential equations.
\newblock {\em ACM Trans. Math. Softw.}, 41(1):2:1--2:39, 2014.

\bibitem{Moiola2011}
A.~Moiola, R.~Hiptmair, and I.~Perugia.
\newblock Plane wave approximation of homogeneous {H}elmholtz solutions.
\newblock {\em Z. Angew. Math. Phys.}, 62(5):809--837, 2011.

\bibitem{SauterZech}
S.~Sauter and J.~Zech.
\newblock A posteriori error estimation of $hp$-{dG} finite element methods for
  highly indefinite {H}elmholtz problems.
\newblock {\em SIAM J. Numer. Anal.}, 53:2414--2440, 2015.

\bibitem{Sloan2004}
I.~H. Sloan and R.~S. Womersley.
\newblock Extremal systems of points and numerical integration on the sphere.
\newblock {\em Adv. Comput. Math.}, 21(1):107--125, 2004.

\bibitem{Womersley2007Online}
R.~S. Womersley.
\newblock Extremal (maximum determinant) points on the sphere ${S}^2$, 2007.
\newblock
  \url{http://web.maths.unsw.edu.au/~rsw/Sphere/Extremal/New/index.html}
  [Online].

\bibitem{Zech}
J.~Zech.
\newblock A posteriori error estimation of $hp$-{DG} finite element methods for
  highly indefinite {H}elmholtz problems.
\newblock Master's thesis, Universit\"at Z\"urich, 2014.
\newblock \url{http://www.math.uzh.ch/compmath/index.php?id=dipl}.

\end{thebibliography}
\end{document}